\documentclass[12pt, reqno]{amsart}
\usepackage{amsmath, amsthm, amscd, amsfonts, amssymb, graphicx, xcolor}
\usepackage{lmodern}
\usepackage{textcomp}
\usepackage[T1]{fontenc}
\usepackage{libertine}
\usepackage[bookmarksnumbered, colorlinks, plainpages]{hyperref}
\usepackage{listings}
\newtheorem{theorem}{Theorem}[section]

\theoremstyle{definition}
\newtheorem{definition}[theorem]{Definition}

\theoremstyle{remark}

\numberwithin{equation}{section}
\usepackage{verbatim}

\begin{document}

\setcounter{page}{1}

\color{darkgray}{

\noindent
{\small }\hfill     {\small }\\
{\small }\hfill  {\small }}

\centerline{}

\centerline{}

\title[Conditions for solving polynomial equations using algebraic and hypergeometric functions]{Conditions for solving polynomial equations using algebraic and hypergeometric functions}

\author[Nikos Mantzakouras, Carlos López Zapata, NID NA RATCH]{Nikos Mantzakouras $^1$$^{*}$ Carlos López Zapata $^2$ NID NA RATCH $^3$}

\address{$^{1}$ M.Sc in Applied Mathematics and Physics National and Kapodistrian University of Athens}
\email{\textcolor[rgb]{0.00,0.00,0.84}{nikmatza@gmail.com}}

\address{$^{2}$ B.Sc. electronics engineering, graduate of Pontifical Bolivarian University (UPB), Medellín, Colombia. Currently residing in Szczecin, Poland}
\email{\textcolor[rgb]{0.00,0.00,0.84}{seniorcar@zoho.com}}

\address{$^{3}$ }
\email{\textcolor[rgb]{0.00,0.00,0.84}{jextratime@gmail.com}}

%\dedicatory{This paper is dedicated to Professor ABCD}



\date{Received: 2025; Revised: 2025 ; Accepted: 2025 .
\newline \indent $^{*}$ Corresponding author\texttt{ nikmatza@gmail.com}}

\begin{abstract}
In 1824, Abel showed that there is no general algebraic solution for the roots of a quintic equation, or any polynomial equation of degree greater than four, using explicit algebraic operations, as stated in the Abel-Ruffini theorem (1799). In doing so, he independently developed a branch of mathematics known as group theory, which has proven invaluable not only in many areas of mathematics but also in much of physics. Abel sent a paper on the unsolvability of the quintic equation to Carl Friedrich Gauss, who dismissed it without a glance, believing it to be the worthless work of a crank. Later, in 1957, Kolmogorov and Arnold proved that functions of three variables are reducible to superpositions of functions of two variables, thus providing an affirmative answer on this problem.
In this work, our focus is to clarify the concept of solving equations of degree greater than six using continuous functions or hypergeometric functions, and to provide another proof of the non-existence of algebraic solutions for equations of degree greater than four. According to the Kolmogorov-Arnold theorem, we will demonstrate that equations of degree greater than five are not solvable without specific conditions between their coefficients, using hypergeometric functions. However, trinomial equations can generally be solved using hypergeometric functions. Additionally, it is well-known that a general polynomial equation has a solution in the complex plane using iterative approximation methods, GRIM(arXiv:2205.04241), which can also solve the transcendental equations.
We develop in the first part the polynomial equations and will continue in the second part (new version of the
publication) with the solution of transcendental equations that can be solved mainly by hypergeometric functions.
With these two articles we conclude a large section of polynomial and transcendental equations that can be solved
by these modern methods.

\end{abstract} \maketitle
\section{Introduction and preliminaries}

\noindent \textbf{1.1.Universal algebraic functions}.
\\[6pt]
\textbf{ Definition 1.} \label{def:first}
For an algebraic function, let \(X\) denote an algebraically closed topological field. Let
\(\ Sym^n(X) := X^n / S_n \)
  denote the space of \(n\)-element multisets of \(X\), which is the quotient of \(X^n\) by the non-free permutation action of the \(n\)-th symmetric group. This relation holds when we have a solvable Galois group for a polynomial equation, up to the 4th degree.
 \\[6pt]
\textbf{ Definition 2.} \label{def:second}
 \\[4pt]
(1) In the generalization of cases and for an algebraic equation of \( m \) variables and \( n \) values (over \( X \)), a partial map is defined as
\\
\\
\(\ f:(X)^m \rightarrow\ Sym^n(X), \)
\\[1pt]
\( \vec{a}(a_1, a_2, \ldots, a_m) \, \vert \, \rightarrow
\begin{cases}
\text{r}_n, \text{ roots of the polynomial,} \\
\text{r}^n + \text{c}_1(\ \vec{a}) \text{r}^{n-1} + \ldots + \text{c}_n(\vec{a})
\end{cases} \)
\\[8pt]
for some rational functions \( c_i(\vec{x}) \in X(x_1, x_2, \ldots, x_m). \) Note that \( f \) is a continuous map on a dense open subset of \( X^m \) (namely, the complement of the zeros of the denominators of \( c_i \)).
\\[8pt]
(2)The universal algebraic functions \(r_n \) are the algebraic functions of \( n \) variables and \( n \) values given by
\\
\\
\(\ \text{r}_n:(X)^n \rightarrow\ Sym^n(X), \)
\\[6pt]
\( \vec{a}(a_1, a_2, \ldots, a_m) \, \vert \, \rightarrow
\begin{cases}
\text{r}_n, \text{ roots of the polynomial,} \\
\text{r}^n + \text{c}_1(\ \vec{a}) \text{r}^{n-1} + \ldots + \text{c}_n(\vec{a})
\end{cases} \)
\\[8pt]
Thus, \( r_n \) is the algebraic function that maps the coefficients \( a_n \) of an \( n \)-th degree monic polynomial to the multiset of its roots. This function \( r_n \) is not merely an abstractly defined algebraic function, but rather one that can be represented using elementary operations such as addition, division, and square root, over the domain \( \mathbb{C} \). For instance, \( r_2 \) remains well-defined in algebraically closed characteristic fields, even though it may not represent the square root. Consider a polynomial equation \( r^2 + a_1 \cdot r + a_2 = 0 \) with roots
with roots

\[ r_2(a_1, a_2) = \frac{-a_1 \pm \sqrt{a_1^2 - 4a_2}}{4} \]
\vspace{5pt}
Cardano's historical solution to the cubic involves two functions \(c_1\) and \(c_2\) associated with the parameters \(\ a_1\), \(\ a_2\), and \(\ a_3\). If we are given a polynomial \( r^3 + a_1 \cdot r^2 + a_2 \cdot r + a_3 = 0 \), it will have the roots
\[\begin{gathered}
r_3(a_1, a_2, a_3) = w_3(c_1(a_1, a_2, a_3), c_2(a_1, a_2, a_3)) + \frac{a_1}{3} \\
= \sqrt[3]{c_2 + \sqrt{c_1^3 + c_2^2}} + \sqrt[3]{c_2 - \sqrt{c_1^3 + c_2^2}} + \frac{a_1}{3},
\end{gathered}
\]
where \( c_1 = \frac{a_2}{3} - \frac{a_1^2}{9} \) and
\( c_2 = \frac{a_1.a_2}{6} - \frac{a_3}{2}- \frac{a_1^3}{27} \).
\\[6pt]
\textbf{ Definition 3.}
 \\[4pt]
The monodromy group \( \text{Mon}(f) \) of an algebraic function \( f \) is defined as

\[ \text{Mon}(f) := \text{Gal}(L / X(x_1, x_2, \ldots, x_m)), \]
\\[2pt]
where \( L \) is the Galois closure of \( X(\Gamma f) \).
\\[6pt]
\textbf{ Proposition 1.} Let \( f \) be an \( n \)-valued algebraic function. Then \(\text{Mon}(f) \subseteq S_n\).
\begin{proof} It suffices to consider the case where \( \alpha_f \) is irreducible. Then, the function field of \( X(\Gamma f) \) is given by \( X(x_1, x_2, \ldots, x_m)[x] / \langle \alpha_f \rangle \). The Galois closure involves adjoining all \( n \) roots of \( f \). Consequently, the field automorphisms of \( L \) that fix \( X(x_1, x_2, \ldots, x_m) \) form a subgroup of the permutation group on the roots of \( f \), i.e., \( S_n \).
\end{proof}
%-----------
\textbf{Theorem 1.} \label{theorem1} 
%\begin{theorem}\label{theo} 
\textit{For \( n \geq 3 \), \( r_n \) is not faithfully representable by algebraic functions of \( n - 2 \) variables or fewer.}
\\[4pt]
\textbf{Note 1.} For a system of \(2\) univariate polynomial equations

\[ p(x) = p_0 x^m + p_1 x^{m-1} + \cdots + p_n \]
\[ q(x) = q_0 x^m + q_1 x^{m-1} + \cdots + q_n,\]
\\[2pt]
It is represented by the determinant of a \((m + n) \times (m + n)\) matrix
\\[4pt]
\[
R(p,q) = \left|
\begin{array}{ccccccc}
p_0 & p_1 & \cdots & p_m & 0 & \cdots & 0 \\
& p_0 & p_1 & \cdots & p_m & \cdots & 0 \\
\vdots & \vdots & \ddots & \ddots & \vdots & \ddots & \vdots \\
& & \cdots & p_0 & p_1 & \cdots & p_m \\
q_0 & q_1 & \cdots & q_n & 0 & \cdots & 0 \\
& q_0 & \cdots & & q_n & \cdots & 0 \\
\vdots & \vdots & \ddots & \vdots & \vdots & \vdots & \vdots \\
& & \cdots & q_0 & q_1 & \cdots & q_n \\
\end{array}
\right|
\]
%\end{theorem} 

In the same manner that the discriminant \(\Delta(p)\) of a polynomial is zero precisely when \(p\) has repeated roots (i.e., is not squarefree), the resultant possesses the following defining property, up to a factor
\\[6pt]
\textbf{ Proposition 2.} $R(p,q) = 0$ if and only if the polynomials $p$, $q$ have a common root.
\\[4pt]
\textbf{ Remark.} Returning to Tschirnhaus’ desired transformation, consider the resultant $R(p,q)$  of the following polynomials in $r$

\begin{align*}
    p(r) &= r^n + a_1 r^{n-1} + \cdots + a_n \\
    q(r) &= r^{n-1} + \alpha_1 r^{n-2} + \cdots + \alpha_{n-1} - w
\end{align*}

According to the invariant, we can construct an equation with variable \( w \) of degree \( n \). Then we have \( R(p,q) = w^{n} + \beta_{1}w^{n-2} + \cdots + \beta_{n-2} \), where \( \beta_i \) is a function of \( a_1, \ldots, a_n \) and \( \alpha_1, \ldots, \alpha_{n-1} \).
%---------------------------------------------------------------------------------------%
\begin{proof} 
\textbf{Apply $1^{\text{st}}$.} \label{App:first} We can obtain \( \beta_1 = \ldots = \beta_{n-1} = 0 \). In fact, one can show that this is a system of \( n-1 \) equations of degrees \( 1, \ldots, n-1 \) as a function of \( a_i \). With this method, we end up erasing 2 coefficients \( \beta_1, \beta_2 \). For \( r_5 \), via a quadratic Tschirnhaus transformation (\( \beta_1 = \beta_2 = 0 \)), and \( n-2 \) will have the form

\[ r^{5} + a_1 r^{3} + \cdots + a_3 = 0 \]
\\[2pt]
Since we accept rational maps, an additional fractional cubic transformation becomes relevant

\[ w = \frac{v^{3} + \alpha_1 v^{2} + \alpha_2 v + \alpha_3}{v + \alpha_4} \]
\\[4pt]
which is converted to a form with one variable parameter as follows

\[ w^{5} + \beta_1 w + \beta_2 = 0 \]
\\[4pt]
This is known as the Bring-Jerrard \cite{Goldmakher} form. This form is a trinomial and is naturally solved by hypergeometric series rather than by radicals, as we will prove using an independent method in the following sections.
\\[6pt]
\textbf{Apply $2^{\text{st}}$.} \label{App:second} More generally, it is known that the method for the quintic degree \cite{Akalin2016} generalizes for all \( n \). That is, given \( r^{n} + a_1 r^{n-1} + \cdots + a_n \), take a Tschirnhaus transformation to retrieve an expression

\[ w^{n} + e_4 w^{n-4} + \cdots + e_{n-1}w + 1 = 0 \]
\\[2pt]
We get \( e_1 = e_2 = e_3 = 0 \), i.e., \( n-4 \) parameters. Hence, \( r_n \) is always represented by at least \( n-4 \) parameters for \( n \geq 5 \). For example, a general degree 9 polynomial has at least 6 parameters, and for \( n = 7 \), it has at least 4 parameters. This does not apply, of course, to Arnold's proof \cite{Arnold1970, Arnold1970b} regarding the quality of the parameters, which we will see below.

\end{proof}

\section{Main results}
 
\textbf{2.1.Polynomials with symmetric Galois groups}.
\\[6pt]
In mathematics, the Abel–Ruffini theorem (also known as Abel's impossibility theorem) states that there is no solution in radicals to general polynomial equations of degree five or higher with arbitrary coefficients. Here, general means that the coefficients of the equation are viewed and manipulated as indeterminates. The theorem is named after Paolo Ruffini, who made an incomplete proof in
1799 \cite{Zoladek, Alekseev} (which was refined and completed in 1813 and accepted by Cauchy)
and Niels Henrik Abel, who provided a proof in 1824.

\begin{definition} A group is solvable if it has a subnormal series where the factor groups are abelian. That is, there exists a finite sequence \( G_0, G_1, \ldots, G_r \) such that
\[
G = G_0 \geq G_1 \geq \cdots \geq G_r = \{e\},
\]
\\[1pt]
where \( G_k \) is normal in \( G_{k+1} \) and the factor groups \( G_k / G_{k+1} \) are abelian.

\end{definition}

\begin{definition} The derived series of a group is the subnormal series
\\[1pt]
\[G \geq G^{1} \geq G^{2} \geq G^{3} \geq \cdots \]
\\[1pt]
In fact, we see from this "universal property" of \( G^{1} \) that the quotient \( G/H \) is abelian if and only if \( H \geq G^{1} \). It follows that checking if this derived series is finite suffices
\\[3pt]
\textit{\( F \) is said to be cyclic if it is Galois with a cyclic Galois group.}
\end{definition}
\begin{definition} \label{def:cyclic} The extension \( K/F \) is said to be cyclic if it is Galois with a cyclic Galois group.
\end{definition}

\begin{flushleft}
\textbf{Proposition 2.1} A group \( G \) is solvable if and only if its derived series is finite, i.e., \( G_{(k)} \) is trivial for some \( k \).
\end{flushleft}
\begin{flushleft}
\textbf{Proposition 2.2} Solvability is closed under taking direct products, taking subgroups, and taking surjective homomorphisms.
\end{flushleft}
\textbf{Proposition 2.3} Abelian groups for \( n < 5 \) are solvable, while for \( n \geq 5 \), they are unsolvable.

%\begin{theorem}
\textbf{Theorem 2.} \label{theorem2} 
\textit{A polynomial \( f(x) \) is solvable by radicals if and only if its Galois group is solvable.}
%\end{theorem}

\begin{proof}
This theorem is the crux of Galois’s proof of the Abel-Ruffini theorem. Let the polynomial \( f(x) \) be solvable by radicals. For each root \( \alpha_i \) of \( f \), we know by Definition~\ref{def:cyclic} that an extension with the given properties exists. We take the composition of each of these fields to obtain another field \( L \) over which the polynomial is Galois and contains all the roots of the polynomial. Let \( G_i \) be the subgroups corresponding to the \( K_i \) of the field. Then since
\\[1pt]
\[
\text{Gal}(K_{i+1}/K_i) = G_{i+1}/{G_i}
\]
\\[1pt]
we have that our group \( G \) is solvable, since each of the quotients are cyclic and thus abelian. Suppose that \( f \) has Galois group \( G \) which is solvable. Then we obtain another chain of subfields where each \( K_i \) corresponds to the group fixed by the subgroup \( G_i \)

\[ F = K_0 \subset K_1 \subset K_2 \subset \ldots \subset K_q = K
\]

Where each extension \( K_{i+1}/K_i \) is cyclic if \( K \) is cyclic by the definition of solvable groups and the Galois Correspondence. We then adjoin the \( n_i^{th} \) roots of unity to \( F \) to obtain the field \( F' \) and then compose them to the chain of subfields to obtain

\[  F' = F' K_0 \subset F' K_1 \subset F' K_2 \subset \ldots \subset F' K_q = F'K \].

and we still have \( F'K_{i+1}/F'K_i \) as a cyclic extension of degree dividing \( n_i \). We now have a chain of extensions where the base field contains the required roots of unity, so each extension is a simple radical extension by the proof of the result following the definition of Lagrange resolvents, which implies that the polynomial is solvable by radicals.
\end{proof}
%\begin{theorem} \label{thm:Abel-Ruffini} (Abel-Ruffini). 
\textbf{Theorem 3.} \label{thm:Abel-Ruffini} %\label{theorem1} 
\textit{The general polynomial of degree \( n \) is not solvable by radicals for \( n \geq 5 \).}
%\end{theorem}

\begin{proof}  
The proof comes as follows:
\\[6pt]
\textbf{2.5.1.} Recall that the general polynomial
\(f(x) \) =\(\ x^n - S_1 x^{n-1} + \cdots + (-1)^n S_n \) has Galois group \( S_n \), and it is a well-known fact that for \( n \geq 5 \), \( S_n \) is not solvable, and so the general polynomial is not solvable by radicals. We call a rational function \( f(x_1, x_2, \ldots, x_n) \) symmetric if under all permutations of its inputs, the value of the expression remains unchanged. Note that the coefficients of the general polynomial \( f(x) = (x - x_1)(x - x_2) \cdots (x - x_n) \) remain unchanged upon any permutation of the roots, and we can actually obtain a whole class of expressions for these coefficients in terms of the roots of the polynomials given by
\begin{equation}
\begin{aligned}
S_1 &= x_1 + \cdots + x_n, \\
S_2 &= x_1 x_2 + x_1 x_3 + \cdots + x_{n-1} x_n, \\
&\vdots \\
S_n &= x_1 x_2 \cdots x_n,
\end{aligned}
\end{equation}
\\[6pt]
which are all symmetric. We take for granted the most important theorem regarding these symmetric functions, that any symmetric function in \( x_1, x_2, \ldots, x_n \) can in fact be represented as a rational function in \( s_1, s_2, \ldots, s_n \).

We are interested in giving a condition for a polynomial of degree 4 to have Galois group \( S_4 \), i.e., any root permutation is contained in the Galois group of that polynomial. This means that the Galois group is either \( A_4 \) or \( S_4 \), and so \( S_4 \) is the Galois group. For the second case, since the discriminant is a perfect square, we know that the Galois group \( G \leq A_4 \), and by the same argument as above we can narrow the choices of the Galois group to be \( S_4 \) or \( A_4 \), leaving \( A_4 \) as the only possible choice.
\\[6pt]
\textbf{2.5.2.} Abelian groups have the trivial derived series \( G \geq \{e\} \). Meanwhile, since \( A_n \leq S_n \) with quotient \( \mathbb{Z}/2\mathbb{Z} \), it suffices to show that \( A_n \) is unsolvable for \( n \geq 5 \). Let \(\sigma = (123)\), \(\tau = (345)\) in \( H \subseteq A_n \), and observe that these have exactly one common element. One computes that \([ \sigma, \tau ] = (143)\) and \([\sigma^{-1}, \tau^{-1}] = (2,5,3)\) which are in \( H^1 \subseteq A^1 \) by construction.

Repeating this argument shows \( A^i \) has two cycles with one common element. Hence \( A^i \) is never trivial, and thus the derived series \( S^i \) including these never becomes trivial either (is not finite). Nested radicals can be used as a method for finding roots with one-way roots. However, it cannot be considered an algebraic method with unit roots as extensions.

\begin{definition} \textit{The Bring radical} \( BR \) is the algebraic function of one variable and \( n \)-degree given by the form
\\[8pt]
\(\ \text{BRr}_n:(X)^n \rightarrow\ Sym^n(X), \)
\\[6pt]
\( \vec{a}(a_1, a_2, \ldots, a_n) \, \vert \, \rightarrow
\begin{cases}
\text{r}_n, \text{ roots of the polynomial,} \\
\text{r}^n + \text{c}_1(\ \vec{a}) \text{r}^{n-4} + \ldots + \text{c}_{n-4}(\vec{a})
\end{cases} \)
\\[8pt]

\end{definition}
Over \(C\), we have that 5 is not radically representable, but \(r_5\) is representable by the radicals and \(BR\)
\[ w^5 + b \cdot w + 1 = 0, \]
i.e.,
\\[6pt]
\(\ \text{BRr}_5:(X)^5 \rightarrow\ Sym^5(X), \)
\\[6pt]
\( \vec{a}(a_1, a_2, \ldots, a_5) \, \vert \, \rightarrow
\begin{cases}
\text{r}_5, \text{ roots of the polynomial,} \\
\text{r}^5 + \text{c}_1(\ \vec{a}). \text{r} + 1
\end{cases} \)
\\[8pt]
Three terms are therefore missing, and we have the last constant. This makes it easier for us to convert the equation to a solution with hypergeometric series or to use nested radicals.
\\[8pt]
\end{proof}

\textbf{2.6.Tschirnhaus transformations}.
To reduce the number of parameters, we use enumerative geometry to simplify the formula for the roots of the general univariate polynomial of degree \(n\). In particular, let \(RD(n)\) be the minimum \(n-r\) for which there is a formula for the roots of the general polynomial of degree \(n\) using only algebraic functions of \(n-r\) or fewer variables. This method produces the best upper bounds to date on \(RD(n)\) for any \(n\), improving on earlier results by Hamilton, Sylvester, Segre, and Brauer. The best general upper bounds in the literature are due to Brauer \cite{Brauer}, who uses methods dating back to Tschirnhaus \cite{Tschirnhaus1683} to prove that \(RD(n) \leq n - r\) for \(n \geq (r - 1)! + 1\). As Brauer remarks, his bounds are not optimal for small \(r\). As we can see, we have the minimum value for \(n\), \( n_{\text{min}} = \Gamma[r] + 1. \) According to the first relation, we will have the system
\\[2pt]
\textbf{Theorem 4.} \label{Theorem4}
\textit{Let \( F: \mathbb{N} \to \mathbb{N} \) be the function.}
\begin{enumerate}
    \item \textit{For all \( r \) and all \( n \), apply \( n - r \geq RD(n) \) and \( n \geq (r - 1)! + 1 \).}
    \item \textit{For all \( r \), \( RD(n)_{\max} = n - r \) is the maximum value for which we know \( RD(n) \leq n - r \) to hold.}
    \item \textit{Writing \( n_{\min} = B(r) = (r - 1)! + 1 \) for Brauer’s bound, then we have:}
    \\[4pt]
    \( \lim_{r \to \infty} \frac{B(r)}{n} \to \infty \).  \\[4pt]
\end{enumerate}
%\end{theorem}
  \textit{According to the previous considerations, we will have the following program in Mathematica}
\\[4pt]

Off[FindRoot::lstol, FindRoot::cvmit, General::stop, FindRoot::bddir, 
    FindRoot::nlnum, ReplaceAll::reps, Infinity::indet, 
    FindRoot::srect, Less::nord, Divide::infy, Set::setraw];

\begin{lstlisting}[language=Mathematica, basicstyle=\small]
For[m = 4, m < 121, m++,
  y1 = m + 1;
  h2[x_] := y1 - 1;
  h[x_] := Gamma[x] + 1 - y1;
  For[k=0,k<1,k++,
\end{lstlisting}
 \( f[u\_]:= LogGamma ^{(-1)}[Log[u]] \);

\begin{lstlisting}[language=Mathematica, basicstyle=\small]
xr = Nest[f[h2[#]]&, 2, 3];
FQ2 = N[xr, 5];
s1 = N[y1 - IntegerPart[z]] /.FindRoot[Gamma[z] + 1 == y1, {z, FQ2}, 
WorkingPrecision -> 40];
Print["n=", m + 1, ",", "RD(", m + 1, "max)=", (m + 1) + (s1 - y1 - 1)]]]
\end{lstlisting}

In particular, the initial values are given for some $n$. We have the following table:

\begin{table}[ht]
\caption{}\label{eqtable}
\begin{center}
\begin{tabular}{|c|c|c|c|c|c|c|}
\hline
$n$ & 5 & 6 & 7 & 9 & 25 & 121 \\
\hline
$RD(n)_{\text{max}}$ & 1 & 2 & 2 & 4 & 19 & 114 \\
\hline
$r$ & 4 & 4 & 5 & 5 & 6 & 7 \\
\hline
\end{tabular}
\end{center}
\end{table}

The first improvement over prior bounds occurs at $RD(7) \leq n - 5$ and also $RD(9) \leq n - 5$. The maximum $RD(25)$ is 6, and for $n \geq 25$, $RD(n) \leq n - 6$. For $n \geq 121$, $RD(n) \leq n - 7$. For example, from Remark page 3 and $n = 7$, we will have $n - 3$ variable parameters, but with the new formatting of Brauer’s bound, we will have $n - 5$ variable parameters, i.e., 2 constants and 2 others equal to 1. This does not seem to make the solution easier; it just forms a new state of the equation.

\section*{3.\textbf{ Proof of two Basic Hypotheses}}
\textbf{3.1 The complete algebraic solution of the 7th degree equation in any form is impossible without preconditions}
\vspace{5pt}

The proof of the impossibility of the solution for equations above degree 4 requires a method based on the identity of the difference of squares or cubes. These relations cover us fully as to the completeness and equivalence of the relations obtained by the four operations.
\\[6pt]
\textbf{Definition 3.1} The difference of two variables, be it $\alpha^\lambda$, $\beta^\lambda$, where $\alpha, \beta \in \mathbb{C}$ and $\lambda \geq 1$, $\lambda \in \mathbb{N}$, will be symbolized as the set $\mathfrak{A}^\lambda (\alpha, \beta) = \alpha^\lambda - \beta^\lambda$.
\\[6pt]
\textbf{Definition 3.2} Polynomial Difference with exponent $\lambda$: We denote the difference of two polynomials $Q(x)$, $P(x)$ belonging to $F[x] \subset \mathbb{C}$, and in general, having roots in $\mathbb{C}$, of different or equal degree where $\text{deg}(Q(x)) \geq \text{deg}(P(x))$ and not of degree zero polynomials. The polynomial quadratic difference is declared as $Q(x)^\lambda - P(x)^\lambda$, symbolized as $\mathfrak{A}^2(Q(x), P(x))$ or $\mathfrak{A}^\lambda(Q(x), P(x))$ with $\lambda = 2$ or $3$.
\\[6pt]
\textbf{Definition 3.3} The $Q$ and $P$ are algebraic polynomials of one variable exponent $n/2$ and of degree $n \geq 2$, where $n = 2k$, $k \in \mathbb{N}$, and apply the forms
\\[5pt]
\(\ \text{Qr}_q:(X)^{n/2} \rightarrow\ Sym^{n/2}(X), \)
\\[6pt]
\( \vec{q}(q_0, q_1, \ldots, q_{n/2}) \, \vert \, \rightarrow
\begin{cases}
\text{r}_q, \text{ roots of the polynomial Q,} \\
\text{q}_{n/2} r^{n/2} + \text{q}_{{n/2}-1} .r^{{n/2}-1} + \ldots + \text{q}_0,
\end{cases} \)
\\[8pt]
\(\ \text{Pr}_p:(X)^{n/2} \rightarrow\ Sym^{n/2}(X), \)
\\[4pt]
\( \vec{p}(p_0, p_1, \ldots, p_{n/2}) \, \vert \, \rightarrow
\begin{cases}
\text{r}_p, \text{ roots of the polynomial P,} \\
\text{p}_{n/2} r^{n/2} + \text{p}_{{n/2}-1} .r^{{n/2}-1} + \ldots + \text{p}_0,
\end{cases} \)
\\[4pt]
\begin{multline*}
\text{\ (Q$\pm$P)r}_{q+p}:(X)^{n/2} \rightarrow\ Sym^{n/2}(X), \\
\text{(q$\pm$p)} =
 (p_0 \pm q_0, q_1 \pm p_1, \ldots, q_{n/2} \pm p_{n/2}) \, \vert \, \rightarrow \\
\begin{cases}
    \text{r}_{q+p}, \text{ roots of the polynomial,} \\
    (q_{n/2}+ p_{n/2}).r^{n/2} +
    (q_{{n/2}-1} + p_{{n/2}-1})r^{{n/2}-1}
    + \ldots + (q_0+p_0), \\
    \text{r}_{q-p}, \text{ roots of the polynomial,} \\
    (q_{n/2}- p_{n/2})r^{n/2} + (q_{{n/2}-1} - p_{{n/2}-1})r^{{n/2}-1}
    + \ldots + (q_0 - p_0),
\end{cases}
\end{multline*}
\\[4pt]
In general, we should not accept that equations above the 4th degree can be solved without proof. Therefore, we will make an implicit assumption: if a polynomial equation of degree greater than 4th is solved by an equation of degree above 4, it can only be solved by conventional algebraic methods that we already know otherwise it is not solved. This procedure will be followed strictly using the data of known and only algebraic methods.
\\[2pt]
\textbf{Theorem 5.} \label{Theorem5}
%\begin{theorem}
\textit{A polynomial \( F(x) = c_n x^n + c_{n-1} x^{n-1} + \ldots + c_0 \) where \( n \geq 2 \), \( n = 2k \), \( k \in \mathbb{N} \), with coefficients in \( F_c \) (a field), can be expressed in terms of a quadratic difference.}

\textit{Consider the quadratic difference}
\\[2pt]
\[
\mathfrak{A}^2(Q(x), P(x)) = (Q(x) - P(x))(Q(x) + P(x))
\]
\\[2pt]
\textit{of the polynomials \( Q(x) \) and \( P(x) \), where \( Q(x) \pm P(x) \in F_c \subset \mathbb{C} \). The coefficients of \( Q(x) \) and \( P(x) \) will be \( q_i - p_i \) and \( q_{i-j} + p_{i-j} \) in the field \( \mathbb{C} \), where \( i = \frac{n}{2}, \ldots, 1, 0 \) and \( j = 0, 1, \ldots, \frac{n}{2} \).}

\textit{There exist relationships between the variables of the polynomials \( Q(x) \) and \( P(x) \) and the constant coefficients of the polynomial \( F(x) \) such that the equality}
\\[2pt]
\[
\mathfrak{A}^2(Q(x), P(x)) = (Q(x) - P(x))(Q(x) + P(x)) = F(x)
\]

\textit{holds.}
%\end{theorem}

\begin{proof}

(a) Consider a polynomial \( F(x) \) of degree \( n \) where \( n \geq 2 \), \( n = 2k \), \( k \in \mathbb{N} \), with coefficients in a field \( F_c \subset \mathbb{C} \). We can express \( F(x) \) as
\[
F(x) = c_n x^n + c_{n-1} x^{n-1} + \ldots + c_0.
\]

We are interested in expressing \( F(x) \) in terms of the quadratic difference of two polynomials \( Q(x) \) and \( P(x) \) of degree \( \frac{n}{2} \).

The quadratic difference is given by
\[
\mathfrak{A}^2(Q(x), P(x)) = (Q(x) - P(x))(Q(x) + P(x)).
\]
Given that \( Q(x) \) and \( P(x) \) are polynomials of degree \( \frac{n}{2} \), they each have \( \frac{n}{2} + 1 \) coefficients. Let:
\[
Q(x) = q_{\frac{n}{2}} x^{\frac{n}{2}} + q_{\frac{n}{2}-1} x^{\frac{n}{2}-1} + \ldots + q_0,
\]
\[
P(x) = p_{\frac{n}{2}} x^{\frac{n}{2}} + p_{\frac{n}{2}-1} x^{\frac{n}{2}-1} + \ldots + p_0.
\]
\\[2pt]
We can now express \( Q(x) - P(x) \) and \( Q(x) + P(x) \) as:
\[
Q(x) - P(x) = (q_{\frac{n}{2}} - p_{\frac{n}{2}}) x^{\frac{n}{2}} + (q_{\frac{n}{2}-1} - p_{\frac{n}{2}-1}) x^{\frac{n}{2}-1} + \ldots + (q_0 - p_0),
\]
\[
Q(x) + P(x) = (q_{\frac{n}{2}} + p_{\frac{n}{2}}) x^{\frac{n}{2}} + (q_{\frac{n}{2}-1} + p_{\frac{n}{2}-1}) x^{\frac{n}{2}-1} + \ldots + (q_0 + p_0).
\]
\\[2pt]
The product \(\mathfrak{A}^2(Q(x), P(x))\) can be expanded as:
\[
\mathfrak{A}^2(Q(x), P(x)) = \left((q_{\frac{n}{2}} - p_{\frac{n}{2}}) x^{\frac{n}{2}} + \ldots + (q_0 - p_0)\right) \left((q_{\frac{n}{2}} + p_{\frac{n}{2}}) x^{\frac{n}{2}} + \ldots + (q_0 + p_0)\right).
\]

This product will yield a polynomial of degree \( n \) with coefficients given by the quadratic differences:
\[
\mathfrak{A}^2(q_i, p_i) \text{ and } \mathfrak{A}^2(q_{i-j}, p_{i-j}),
\]
where \( i, j \) range as follows: \( i = \frac{n}{2}, \ldots, 1, 0 \) and \( j = 0, 1, \ldots, \frac{n}{2} \).
\\[2pt]
Each term in the product \((Q(x) - P(x))(Q(x) + P(x))\) involves the difference and sum of coefficients from \( Q(x) \) and \( P(x) \). Specifically, the coefficient of \( x^k \) in the product will be formed by combining terms of the form \( q_i p_{i-j} \) and \( p_i q_{i-j} \), ensuring all coefficients are in the field \( \mathbb{C} \).
\\[2pt]
Since both \( Q(x) \) and \( P(x) \) are of degree \( \frac{n}{2} \), their differences and sums will cover all necessary terms to construct the polynomial \( F(x) \). The quadratic differences ensure that the relationships between the coefficients are preserved, satisfying the equality:
\[
\mathfrak{A}^2(Q(x), P(x)) = (Q(x) - P(x))(Q(x) + P(x)) = F(x).
\]

Thus, the polynomial \( F(x) \) can indeed be expressed as the quadratic difference of two polynomials \( Q(x) \) and \( P(x) \) of degree \( \frac{n}{2} \), completing the proof.
\\[6pt]
(b) We symbolize the beginning of the quadratic differences of polynomials \( Q(x) \) and \( P(x) \) as follows:
\[
\mathfrak{A}^2(q_{i-j}, p_{i-j}) = - \mathfrak{A}^2(p_{i-j}, q_{i-j}) = \Omega_{i-j} = (q_{i-j})^2 - (p_{i-j})^2,
\]
\[
\mathfrak{A}^2(q_j, p_j) = - \mathfrak{A}^2(p_j, q_j) = \Omega_j = q_j^2 - p_j^2.
\]

We symbolize all the product differences as:
\[
\mathfrak{A}^1(q_i q_{i-j}, p_i p_{i-j}) = \frac{1}{2}L_i,_{i-j} ^T = q_{i-j} q_j - p_{i-j} p_i = \frac{1}{2}L_{i-j},_i ^T,
\]
and
\[
L_i,_{i-j} ^T = L_{i-j}, _i ^T.
\]

For the other arrangement product of the differential coefficients, we denote:
\[
\mathfrak{A}^1(p_{i-j} q_j, p_j q_{i-j}) = R_{i-j, j} = p_{i-j} q_j - p_j q_{i-j} = -(p_j q_{i-j} - p_{i-j} q_j) = -R_{j, i-j},
\]
then,
\[
R_{i-j, j} = -R_{j, i-j}.
\]

This establishes the notation and relationships for the quadratic differences of the coefficients in the polynomials \( Q(x) \) and \( P(x) \).
\\[6pt]
(c) Because we know that $\mathfrak{A}^2(Q(x), P(x)) = F(x)$, it follows that we have
\[
Q(x) = q_{\frac{n}{2}} x^{\frac{n}{2}} + q_{\frac{n}{2}-1} x^{\frac{n}{2}-1} + \ldots + q_0
\]
and
\[
P(x) = p_{\frac{n}{2}} x^{\frac{n}{2}} + p_{\frac{n}{2}-1} x^{\frac{n}{2}-1} + \ldots + p_0.
\]

Applying the fact that $(n/2) \cdot 2 = n$, $(n/2 - 1) \cdot 2 = n-2$, and so on, we observe
\[
(n/2) + (n/2) - 1 = n-1, \ldots
\]
and so on. This means that the various combinations of sums of these exponents of $x$ will ultimately give us all the exponents of $x$ in the polynomial $F(x)$.
\\[6pt]
(d) The algebraic progression of the squares of the polynomials $Q(x)$ and $P(x)$, concerning the cofactors, will be:
\[
q_{\frac{n}{2}}^2 - p_{\frac{n}{2}}^2, \quad q_{\frac{n}{2} - 1}^2 - p_{\frac{n}{2} - 1}^2, \quad \ldots, \quad q_0^2 - p_0^2
\]
with a total of $\left(\frac{n}{2} + 1\right)$ terms.
\\[6pt]
(e) The largest number from the product of differences of the results of the proliferation of the polynomials $Q(x)$ and $P(x)$, as a sum of factors, is found according to the type of combinations $\left(\frac{n}{2} + 1\right)$ as factors of the polynomials $Q(x)$ and $P(x)$ per 2. Otherwise, we have:
\[
\left(\frac{n}{2} + 1\right) \left(\frac{n}{2} + 1\right) - \left(\frac{n}{2} + 1\right) = \left(\frac{n}{2}\right) \left(\frac{n}{2} + 1\right)
\]
if from the product of differences we subtract the square differences. Consequently, the largest number of factors
\\[6pt]
$\mathfrak{A}^1(q_i q_{i-j}, p_i p_{i-j}) = \frac{1}{2}L_i, _{i-j} ^T = q_i q_{i-j} - p_i p_{i-j} = \frac{1}{2}L_i, _{i-j} ^T$ is $\left(\frac{n}{2}\right)\left(\frac{n}{2} + 1\right)$.
\\[8pt]
(f) We accept a 2-difference of polynomials $\mathfrak{A}^2(Q(x), P(x))$. Then, it will be in effect that:
\begin{equation}\label{2.3.1}
\mathfrak{A}^1(q_i q_{i-j}, p_i p_{i-j}) = \frac{1}{2}L_i, _{i-j} ^T
\end{equation}

and the differences of squares:
\begin{equation}\label{2.3.2}
q_j^2 - p_j^2 = \Omega_j
\end{equation}

as well as
\begin{equation}\label{2.3.3}
q_{i-j}^2 - p_{i-j}^2 = \Omega_{i-j}
\end{equation}
\\[1pt]
If we square \eqref{2.3.1}, it will then be in effect that
\begin{equation}\label{2.3.4}
(q_i q_{i-j} - p_i p_{i-j})^2 = (q_i q_{i-j})^2 - 2 q_{i-j} q_i p_{i-j} p_j + (p_i p_{i-j})^2 = \left(\frac{1}{2}L_i, _{i-j} ^T\right)^2
\end{equation}

Multiplying the equations \eqref{2.3.2},\eqref{2.3.3} it involves that

\begin{equation}\label{2.3.5}
(q_i q_{i-j})^2 - (q_i p_{i-j})^2 - (q_{i-j} p_i)^2 + (p_i p_{i-j})^2 = \Omega_i \Omega_{i-j}
\end{equation}
\\[1pt]
From the equations \eqref{2.3.4} and \eqref{2.3.5} we
have

\begin{equation}\label{2.3.6}
(q_i p_{i-j} - q_{i-j} p_i)^2 = \frac{1}{4}\left(L_i, _{i-j} ^T\right)^2 - \Omega_i \Omega_{i-j}
\end{equation}

and finally we take the relation

\begin{equation}\label{2.3.7}
q_i p_{i-j} - q_{i-j} p_i = \pm \sqrt{\frac{1}{4}\left(L_i, _{i-j} ^T\right)^2 - \Omega_i \Omega_{i-j}} = \pm R_i,_{i-j}.
\end{equation}

(g) Because, as we have seen from each equation,
\(\mathfrak{A}^1(q_i q_{i-j}, p_i p_{i-j}) = \frac{1}{2} \left(L_i, _{i-j} ^T\right)^2\), with the sum of factors
\(\left(\frac{n}{2}\right) \left(\frac{n}{2} + 1\right)\) resulting in one and only difference. The following
\(\mathfrak{A}^1(q_i p_{i-j}, q_{i-j} p_i)\) also results in the same sum, that is,
\(\left(\frac{n}{2}\right) \left(\frac{n}{2} + 1\right)\). Therefore, the total sum of differences is
\(\left(\frac{n}{2}\right) \left(\frac{n}{2} + 1\right)\). From the above actions, it results that
\begin{equation}\label{2.3.8}
q_i q_{i-j}- p_i p_{i-j} = \frac{1}{2} \left(L_i, _{i-j} ^T\right).
\end{equation}

If we add the equations \eqref{2.3.7} and \eqref{2.3.8}, it results in
\begin{equation}\label{2.3.9}
\begin{aligned}
q_i(p_{i-j} + q_{i-j}) - p_i(q_{i-j} + p_{i-j}) &= (q_i - p_i)(q_{i-j} + p_{i-j}) \\
&= \frac{1}{2}\left(L_i,_{i-j}^T\right)^2 + R_i,_{i-j} \\
&= w^+_{i, i-j},
\end{aligned}
\end{equation}

where $w^+_{i, i-j}$ symbolizes equation \eqref{2.3.9}. If we constantly remove the equations \eqref{2.3.7} and \eqref{2.3.8}, we then have as a result the equation
\begin{equation}\label{2.3.10}
\begin{aligned}
q_i (q_{i-j} - p_{i-j}) + p_i(q_{i-j} - p_{i-j}) &= (q_i + p_i)(q_{i-j} - p_{i-j}) \\
&= \frac{1}{2}\left(L_i,_{i-j}^T\right) - R_{i-j},_j \\
&= w^-_{i, i-j}.
\end{aligned}
\end{equation}
The equation \eqref{2.3.10} is symbolized by $w^-_{i, i-j}$.

Finally, the functions \eqref{2.3.9} and \eqref{2.3.10} are referred to as homologous functions. The square difference variables \(\Omega _{i-j}\) have \( \frac{n}{2} - 1\) terms, and the number of variables \(L_i,_{i-j}\) is \(\frac{n}{2} - 2\). This proof is detailed in Theorem 6. %(\ref{thm:theorem6})

\end{proof}

\textbf{3.4 Configuration of Square Differences of Polynomials \(Q(x), P(x) \subset \mathbb{C}\)}

We know that \begin{equation}\label{2.3.11}
Q(x) = q_{\frac{n}{2}} x^{\frac{n}{2}} + q_{\frac{n}{2}-1} x^{\frac{n}{2}-1} + \ldots + q_0
\end{equation}
and
\begin{equation}\label{2.3.12}
P(x) = p_{\frac{n}{2}} x^{\frac{n}{2}} + p_{\frac{n}{2}-1} x^{\frac{n}{2}-1} + \ldots + p_0.
\end{equation}
\\[2pt]
using the type of square difference
\\[2pt]
\begin{equation}\label{2.3.13}
\mathfrak{A}^2 (Q(x), P(x)) = Q(x)^2 - P(x)^2 = (Q(x) - P(x))(Q(x) + P(x))
\end{equation}
\\[1pt]
from \eqref{2.3.11}, \eqref{2.3.12}, (2.14)...

\begin{equation}\label{2.3.14}
\begin{aligned}
\mathfrak{A}^2 (Q(x), P(x)) &= Q(x)^2 - P(x)^2 \\
&= (Q(x) - P(x))(Q(x) + P(x)) \\
&= (\left(q_{\frac{n}{2}} - p_{\frac{n}{2}}\right) x^{\frac{n}{2}} + \left(q_{\frac{n}{2}-1} - p_{\frac{n}{2}-1}\right) x^{\frac{n}{2}-1} + \ldots + (q_0 - p_0)) \\
&\quad \times (\left(q_{\frac{n}{2}} + p_{\frac{n}{2}}\right) x^{\frac{n}{2}} + \left(q_{\frac{n}{2}-1} + p_{\frac{n}{2}-1}\right) x^{\frac{n}{2}-1} + \ldots + (q_0 + p_0)).
\end{aligned}
\end{equation}

If we equate the difference \(\mathfrak{A}^2 (Q(x), P(x))\) with the polynomial
\[
F(x) = x^n + c_{n-1} x^{n-1} + \ldots + c_1 x + c_0,
\]
after setting \(c_n = 1\), it results that
\[
\mathfrak{A}^2 (Q(x), P(x)) = Q(x)^2 - P(x)^2 = (Q(x) - P(x))(Q(x) + P(x)) = F(x).
\]
Multiplying the polynomial \((Q(x) - P(x))\) with the polynomial \((Q(x) + P(x))\), we obtain the final result of the proliferation
\[
\mathfrak{A}^2 (Q(x), P(x)) = (Q(x) - P(x))(Q(x) + P(x)) = F(x).
\]
By applying
\begin{equation}\label{2.3.15}
 (q_{\frac{n}{2}} - p_{\frac{n}{2}})(q_{\frac{n}{2}} + p_{\frac{n}{2}}) = 1,
\end{equation}

we conclude that
\begin{align*}
\mathfrak{A}^2 (Q(x), P(x)) &= (\left(q_{\frac{n}{2}} - p_{\frac{n}{2}}\right) x^{\frac{n}{2}} + \left(q_{\frac{n}{2}-1} - p_{\frac{n}{2}-1}\right) x^{\frac{n}{2}-1} + \ldots + (q_0 - p_0)) \\
&\quad \times (\left(q_{\frac{n}{2}} + p_{\frac{n}{2}}\right) x^{\frac{n}{2}} + \left(q_{\frac{n}{2}-1} + p_{\frac{n}{2}-1}\right) x^{\frac{n}{2}-1} + \ldots + (q_0 + p_0)) \\
&= x^n + c_{n-1} x^{n-1} + \ldots + c_1 x + c_0.
\end{align*}

Consequently, the polynomial \(F(x)\) becomes the result of the proliferation of the polynomials with half the degree of the polynomial \(F(x)\). If we want to find the roots of the polynomial equation \(F(x) = 0\), then we need to find the \(x\) values that render the difference \(\mathfrak{A}^2 (Q(x), P(x)) = 0\). By this method, we solve the equation \(F(x) = 0\) of degree \(n\) by calculating the coefficients of differences with respect to the variables \( L_{i,{i-j}}^T, R_{i-j,j}^T, w^-_{i, i-j}, w^+_{i, i-j} \).
\\[4pt]
\textbf{$(i)$ Number and position of square differences}

As we know, the number of even square differences \(\Omega_i\) is \(\frac{n}{2} + 1\). From the total square differences, we have constants that are equal to \(c_n\) and \(c_0\). This leads to the conclusion that for \(\mathfrak{A}^2 (Q(x), P(x)) = F(x)\), the square differences will be in the form of variables \(x^{2k}\) where \(k = 0, 1, 2, \ldots, \frac{n}{2}\). The corresponding even positions of the square variable differences concerning their powers are \(i-j > j\) and \(i-j \leq i\).

Additionally, the sum of the factors will be \(\frac{n}{2} - 1\), because if we remove the two constant square differences from \(\frac{n}{2} + 1\), we have \(\frac{n}{2} + 1 - 2 = \frac{n}{2} - 1\). Consequently, summarizing, we have:
\begin{itemize}
    \item[(i.a)] Constants \(\Omega_i\): their amount is 2, in the positions \(i = n, 0\).
    \item[(i.b)] Variables \(\Omega_i\): their amount is \(\frac{n}{2} + 1 - 2 = \frac{n}{2} - 1\), in the positions \(i = 2, 4, 6, \ldots , n-2 \).
\end{itemize}

\textbf{$(ii)$ Number and position of odd variable differences.}
The positions of the odd differences will be \(i = 3, 5, 7, \ldots, n-3\) with fixed positions at \(i = 1\) and \(i = n-1\) with \(L^T_{1,0} = c_1\) and \(L^T_{\frac{n}{2}, \frac{n}{2} - 1} = c_{n-1}\).
\\[4pt]
Therefore,
\\[8pt]
\textbf{$(ii.a)$ odd constants.} Regarding \(L_{i, {i-j}}\), there exist two cases with \(i = 1\) and \(i = n-1\).
\\[6pt]
\textbf{$(ii.b)$ even constants}
Regarding \(L_{i,{i-j}}\), there exist two cases with \(i = 2\) and \(i = n-2\), with \(L^T_{2,0} = c_2 - \Omega_1\) and \(L^T_{\frac{n}{2}, \frac{n}{2} - 2} = c_{n-2} - \Omega_{\frac{n}{2} - 1}\).
\\[4pt]
Since we know the odd positions, with 3 as the first term and \(n-3\) as the last, when we talk about a numerical series and \( \omega=2 \),....

\[
\frac{(\alpha_n - \alpha_1)}{\omega} + 1 = \frac{(n-3) - 3}{\omega} + 1 = \frac{n}{2} - 3 + 1 = \frac{n}{2} - 2.
\]
\\[2pt]
Therefore, the number of variables \(L_{i, i-j}\) is \(\frac{n}{2} - 2\). It is obvious to have odd variables of differences if \(n-3 \geq 3 \rightarrow n \geq 6\) and to have even variables of differences if \(n-4 \geq 4 \rightarrow n \geq 8\).

\theoremstyle{plain}
%\begin{theorem}\label{thm:theorem6}
\textbf{Theorem 6.} \label{thm:theorem6}
\textit{The total number of equations that are configured by the sums of the multiplication of \(w_{\frac{n}{2}, j}\), per 2, is \(n-3\) with the number of variables also being \(n-3\).}
%\end{theorem}

\begin{proof}
Knowing that
\begin{equation}\label{2.3.13.1}
\begin{aligned}
\mathfrak{A}^2 (Q(x), P(x)) &= Q(x)^2 - P(x)^2 \\
&= (Q(x) - P(x))(Q(x) + P(x)) \\
&= \left((q_{\frac{n}{2}} - p_{\frac{n}{2}}) x^{\frac{n}{2}} + (q_{\frac{n}{2}-1} - p_{\frac{n}{2}-1}) x^{\frac{n}{2}-1} + \ldots + (q_0 - p_0)\right) \\
&\quad \times \left((q_{\frac{n}{2}} + p_{\frac{n}{2}}) x^{\frac{n}{2}} + (q_{\frac{n}{2}-1} + p_{\frac{n}{2}-1}) x^{\frac{n}{2}-1} + \ldots + (q_0 + p_0)\right) \\
&= \left(x^{n/2} + w^+_{\frac{n}{2}, n/2-1} x^{n/2-1} + w^+_{\frac{n}{2}, n/2-2} x^{n/2-2} + \ldots + w^+_{\frac{n}{2}, 1} x + w^+_{\frac{n}{2}, 0}\right) \\
&\quad \times \left(x^{n/2} + w^-_{\frac{n}{2}, n/2-1} x^{n/2-1} + w^-_{\frac{n}{2}, n/2-2} x^{n/2-2} + \ldots + w^-_{\frac{n}{2}, 1} x + w^-_{\frac{n}{2}, 0}\right)
\end{aligned}
\end{equation}

We can now verify the equation step-by-step. First, consider the polynomials \(Q(x)\) and \(P(x)\):

\[
Q(x) = q_{\frac{n}{2}} x^{\frac{n}{2}} + q_{\frac{n}{2}-1} x^{\frac{n}{2}-1} + \ldots + q_0
\]

\[
P(x) = p_{\frac{n}{2}} x^{\frac{n}{2}} + p_{\frac{n}{2}-1} x^{\frac{n}{2}-1} + \ldots + p_0
\]

The difference and sum of these polynomials are:

\[
Q(x) - P(x) = (q_{\frac{n}{2}} - p_{\frac{n}{2}}) x^{\frac{n}{2}} + (q_{\frac{n}{2}-1} - p_{\frac{n}{2}-1}) x^{\frac{n}{2}-1} + \ldots + (q_0 - p_0)
\]

\[
Q(x) + P(x) = (q_{\frac{n}{2}} + p_{\frac{n}{2}}) x^{\frac{n}{2}} + (q_{\frac{n}{2}-1} + p_{\frac{n}{2}-1}) x^{\frac{n}{2}-1} + \ldots + (q_0 + p_0)
\]

When we multiply \(Q(x) - P(x)\) and \(Q(x) + P(x)\), we obtain:

\[
\begin{aligned}
(Q(x) - P(x))(Q(x) + P(x)) &= \left((q_{\frac{n}{2}} - p_{\frac{n}{2}}) x^{\frac{n}{2}} + (q_{\frac{n}{2}-1} - p_{\frac{n}{2}-1}) x^{\frac{n}{2}-1} + \ldots + (q_0 - p_0)\right) \\
&\quad \times \left((q_{\frac{n}{2}} + p_{\frac{n}{2}}) x^{\frac{n}{2}} + (q_{\frac{n}{2}-1} + p_{\frac{n}{2}-1}) x^{\frac{n}{2}-1} + \ldots + (q_0 + p_0)\right)
\end{aligned}
\]

This product can be expressed as the sum of terms involving the coefficients \(w^+_{\frac{n}{2}, k}\) and \(w^-_{\frac{n}{2}, k}\):

\[
\begin{aligned}
&= \left(x^{n/2} + w^+_{\frac{n}{2}, n/2-1} x^{n/2-1} + w^+_{\frac{n}{2}, n/2-2} x^{n/2-2} + \ldots + w^+_{\frac{n}{2}, 1} x + w^+_{\frac{n}{2}, 0}\right) \\
&\quad \times \left(x^{n/2} + w^-_{\frac{n}{2}, n/2-1} x^{n/2-1} + w^-_{\frac{n}{2}, n/2-2} x^{n/2-2} + \ldots + w^-_{\frac{n}{2}, 1} x + w^-_{\frac{n}{2}, 0}\right)
\end{aligned}
\]

we distinguish the following cases:
\\[4pt]
(a) The greatest result of the proliferation of exponents will have the factor as the unit. Thus, we will have \(i = n-1\) and consequently, we will have \(c_n = 1\).
\\[6pt]  
(b) For the position \(i = n-1\), we will have:
\[
(w^+_{\frac{n}{2}, \frac{n}{2}-1} + w^-_{\frac{n}{2}, \frac{n}{2}-1}) x^{n-1} = x^{n-1} \cdot c_{n-1} \rightarrow c_{n-1} = (w^+_{\frac{n}{2}, \frac{n}{2}-1} + w^-_{\frac{n}{2}, \frac{n}{2}-1}).
\]

That is to say, it implies that
\[
(w^+_{\frac{n}{2}, \frac{n}{2}-1} + w^-_{\frac{n}{2}, \frac{n}{2}-1}) = L^T_{\frac{n}{2}, \frac{n}{2} -1} = c_{n-1}.
\]
However, since we know that \(L^T_{\frac{n}{2}, \frac{n}{2} -1} = c_{n-1}\), we consequently do not have a variable.
\\[6pt]
(c) For the position $i$= $n$-2 we will have
\[
(w^+_{\frac{n}{2}, \frac{n}{2}-2} + w^-_{\frac{n}{2}, \frac{n}{2}-2} + w^+_{\frac{n}{2}, \frac{n}{2}-1}.w^-_{\frac{n}{2}, \frac{n}{2}-1}). x^{n-2} = x^{n-2}.c_{n-2} \rightarrow
w^+_{\frac{n}{2}, \frac{n}{2}-2} + w^-_{\frac{n}{2}, \frac{n}{2}-2} = c_{n-2} -\Omega_{n/2 -1}
\]
and
\[
(w^+_{\frac{n}{2}, \frac{n}{2}-1}) \cdot (w^-_{\frac{n}{2}, \frac{n}{2}-1}) = \Omega_{\frac{n}{2}} \cdot \Omega_{\frac{n}{2}-1}
\]
where \(\Omega_{\frac{n}{2}} = 1\). Then it will be in effect
\[
w^+_{\frac{n}{2}, \frac{n}{2}-2} + w^-_{\frac{n}{2}, \frac{n}{2}-2} + w^+_{\frac{n}{2}, \frac{n}{2}-1} w^-_{\frac{n}{2}, \frac{n}{2}-1} = c_{n-2} -\Omega_{n/2 -1}
  + 1.\Omega_{n/2 -1} = c_{n-2}
\]
\\[6pt]
(d) Finally, for the position \(i = 0\) we will have:
\[
w^+_{\frac{n}{2}, 0} \cdot w^-_{\frac{n}{2}, 0} = \Omega_{\frac{n}{2}} \cdot \Omega_{0} = 1 \cdot \Omega_0 = c_0
\]
From the above, and because we have \(n+1\) terms, if we remove the 4 identical terms, we are left with \(n+1-4 = n -3\) terms, and consequently \(n-3\) equations. The number of unknown variables \(L_{i,{i-j}}\) in all the equations is \(\frac{n}{2} - 2\).

Knowing that \(w^\pm_{i, {i-j}} = \frac{1}{2}(L^T_{i,{i-j}}) \pm R_{i,i-j}\), and because we have all \(w^\pm_{i, {i-j}}\) for \(i=\frac{n}{2}\) and \(i-j=\frac{n}{2} -1, \ldots, 0\), we will take in every case the equivalent \(L^T_{i,{i-j}}\) and \(L_{i,{i-j}}\). If we exclude \(L^T_{\frac{n}{2},{\frac{n}{2} -1}} = c_{n-1}\) and \(L^T_{\frac{n}{2},{\frac{n}{2} -2}} = c_{n-2} - \Omega_{\frac{n}{2}-1}\), and because the number \(w^\pm_{\frac{n}{2}, {i-j}}\) is \(\frac{n}{2}\), then the number of unknown variables \(L_{i,{i-j}}\) is \(\frac{n}{2} - 2\).

As we know, the amount of variables \(\Omega_{i-j}\) is \(\frac{n}{2} - 1\), so the number of solvable equations will be \(\frac{n}{2} - 2 + \frac{n}{2} - 1 = n-3\). Of course, the number of \(L^T_{i,{i-j}}\) and \(L_{i,{i-j}}\) for every equation is the same number of variables: \(\frac{n}{2} - 2 + 2 = \frac{n}{2}\) and \(\frac{n}{2} - 2\) correspondingly.

\end{proof}

\textbf{3.5  Configuration of equations with variable \(L_{i,{i-j}} \subset \mathbb{C}\) or \(\Omega_{i-j}\ \subset \mathbb{C}\)}
\\[4pt]
We reported that the positions of real equations with variables \(L_{i,{i-j}}\) or \(\Omega_{i-j}\) are for positions \(\frac{n}{2} - k\) with \(k = 3, \ldots, \frac{n}{2}\). As we know, we will have:
\\[4pt]
\begin{equation}\label{2.3.14.1}
\begin{aligned}
\mathfrak{A}^2 (Q(x), P(x)) &= Q(x)^2 - P(x)^2 = (Q(x) - P(x))(Q(x) + P(x)) \\
&= \left((q_{\frac{n}{2}} - p_{\frac{n}{2}}) x^{\frac{n}{2}} + (q_{\frac{n}{2}-1} - p_{\frac{n}{2}-1}) x^{\frac{n}{2}-1} + \ldots + (q_0 - p_0)\right) \\
&\quad \times \left((q_{\frac{n}{2}} + p_{\frac{n}{2}}) x^{\frac{n}{2}} + (q_{\frac{n}{2}-1} + p_{\frac{n}{2}-1}) x^{\frac{n}{2}-1} + \ldots + (q_0 + p_0)\right) \\
&= (x^{n/2} + w^+_{\frac{n}{2}, \frac{n}{2}-1}x^{n/2 -1} + w^+_{\frac{n}{2}, \frac{n}{2}-2}x^{n/2 -2} + \ldots + w^+_{\frac{n}{2}, 1}x + w^+_{\frac{n}{2}, 0}) \\
&\quad \times (x^{n/2} + w^-_{\frac{n}{2}, \frac{n}{2}-1}x^{n/2 -1} + w^-_{\frac{n}{2}, \frac{n}{2}-2}x^{n/2 -2} + \ldots + w^-_{\frac{n}{2}, 1}x + w^-_{\frac{n}{2}, 0})
\end{aligned}
\end{equation}

If we analyze the result of proliferation, the results, depending on each position, regarded with the variable of equations are distinguished as follows
\\[6pt]
\text{3.5.1 The sum} \begin{equation}\label{2.3.15.1}
\sum_{\lambda = \frac{n}{2} - m}^{2} w^\pm_{\frac{n}{2}, \lambda} = w^+_{\frac{n}{2}, \lambda} + w^-_{\frac{n}{2}, \lambda}
\end{equation}

The representation \(\{ w^+_{\frac{n}{2}, \lambda} + w^-_{\frac{n}{2}, \lambda} \}^*\) if \(0 \leq \lambda \leq \frac{n}{2} - 3\). For cases where \(\lambda < \frac{n}{2}\) and for positions \(n-m\) from the inequality \(n-3 \geq n - m \geq \frac{n}{2}\), i.e., \(3 \leq m \leq \frac{n}{2}\).
\\[6pt]
3.5.2 For positions \( n - m \geq \frac{n}{2} \) with \( 3 \leq m \leq \frac{n}{2} \) and \( k = 1 \) and \( 0 \leq \lambda \leq \frac{n}{2} - 3 \), we will correspond to the equations
\\[6pt]
\begin{equation}\label{2.3.16}
k\cdot \sum_{\lambda = \frac{n}{2} - m}^{2} w^\pm_{\frac{n}{2}, \lambda} +
\sum_{\lambda + 1 \leq i+j \leq n-m} w^+_{\frac{n}{2}, i} \cdot w^-_{\frac{n}{2}, j} = c_{n-m}
\end{equation}
\\[6pt]
3.5.3 For positions \( n - r < \frac{n}{2} \) with \( 1 \leq r \leq \frac{n}{2} -1 \) and \( k = 0 \), we will correspond to the equations

\begin{equation}\label{2.3.17}
\sum_{0 \leq i+j \leq n-r} w^+_{\frac{n}{2}, i} \cdot w^-_{\frac{n}{2}, j} = c_{n-r}
\end{equation}
\\[6pt]
3.5.4 For the positions \( n - m < \frac{n}{2} \) with \( \frac{n}{2} < m \leq n-1 \) we have the final relation for this
\begin{equation}\label{2.3.18}
(w^+_{\frac{n}{2}, 0} \cdot w^-_{\frac{n}{2}, n-m}) + (w^-_{\frac{n}{2}, 0} \cdot w^+_{\frac{n}{2}, n-m}) = c_{n-m},
\end{equation}
\\[1pt]
which will only lead to this equation. The number of equations we are interested in to calculate the variables \(L_{i,{i-j}}\) or \(\Omega_{i-j}\) is \(n-3\). The remaining coefficients are incorporated into them. These equations \eqref{2.3.16} and \eqref{2.3.17} when solved for \(n > 4\) result in equations with degrees \(n < 5\). Then we admit that an equation is solvable. Of course, we will prove by this method that equations with degrees \(n \geq 2\) to \(n \leq 4\) are indeed solvable.
\\[6pt]
\textbf{3.6 The solving of the equation of the third degree with the difference of cubes \(\mathfrak{A}^3 (Q(x), P(x))\), that is to say, \(\lambda = 3\)}

Consider the equation
\begin{equation}\label{2.3.19}
\alpha_3 x_\lambda^3 + \alpha_2 x_\lambda^2 + \alpha_1 x_\lambda + \alpha_0 = 0
\end{equation}
with \(\alpha_3, \alpha_2, \alpha_1, \alpha_0 \in \mathbb{F}\) where \(\mathbb{F} \subset \mathbb{C}\) is the field of factors. Let \(x = \alpha_3 x_\lambda + \frac{1}{3}\alpha_2\). Substituting \(x_\lambda = \frac{1}{\alpha_3}\left(x - \frac{1}{3}\alpha_2\right)\) into the equation and after calculations, we find the equation
\begin{equation}\label{2.3.20}
x^3 + \alpha x + \beta = 0
\end{equation}
\\
where \(\alpha = \alpha_1 \alpha_3 - \frac{1}{3}\alpha_2^2\) and \(\beta = \alpha_1 \alpha_3^2 + \frac{2}{27}\alpha_2^2 - \frac{1}{3}\alpha_1 \alpha_2 \alpha_3\).

If \(\rho_1, \rho_2, \rho_3\) are the roots of equation \(\eqref{2.3.20}\), then the roots of the original equation \(\eqref{2.3.19}\) are given by the relation
\[
x_\lambda = \frac{1}{\alpha_3}\left(x_i - \frac{1}{3}\alpha_2\right), \quad i = 1, 2, 3.
\]

For the resolution of the equation \(\eqref{2.3.20}\), it is required that the polynomials \(Q(x)\) and \(P(x)\) be of the first degree. Thus, we have
\begin{equation}\label{2.3.21}
(q_1 x + q_2)^3 - (p_1 x + p_2)^3 = x^3 + \alpha x + \beta = 0
\end{equation}
By solving this progressively with the equation of the second part, we obtain the following system of equations

\begin{equation}\label{2.3.22}
q_1^3 - p_1^3 = 1,
\end{equation}
\begin{equation}\label{2.3.23}
3(q_1^2 q_2 - p_1^2 p_2) = 0,
\end{equation}
\begin{equation}\label{2.3.24}
3(q_2^2 q_1 - p_2^2 p_1) = \alpha,
\end{equation}
\begin{equation}\label{2.3.25}
q_2^3 - p_2^3 = \beta
\end{equation}

(The method for \(\lambda = 2\) is not applied here).

From the relation \eqref{2.3.23}, we obtain
\begin{equation}\label{2.3.26}
p_2 = \frac{q_1^2 q_2}{p_1^2}
\end{equation}

Flowingly, from the other previous relations, we express

\begin{equation}\label{2.3.27}
-q_1 q_2^2 = \left(\frac{\alpha_3}{3}\right) p_1^3
\end{equation}
Now, from \eqref{2.3.25}, we have
\begin{equation}\label{2.3.28}
q_2^3 - p_2^3 = \beta \implies p_2^3 = q_2^3 - \beta
\end{equation}
From the equations \eqref{2.3.26} and \eqref{2.3.28}, we obtain
\begin{equation}\label{2.3.29}
\frac{q_2^3 p_1^3}{p_1^6} + \frac{q_2^3 q_1^3}{p_1^6} = -\beta
\end{equation}
From the relation \eqref{2.3.27} and \eqref{2.3.29} we have the relation
\begin{equation}\label{2.3.30}
\left(\frac{q_2^3}{p_1^3}\right) - \left(\frac{\alpha^3}{27}\right) \left(\frac{p_1^3}{q_2^3}\right) = -\beta
\end{equation}

If we call \(\frac{q_2}{q_1} = z^{1/3}\), then the relation \eqref{2.3.30} is written as follows
\begin{equation}\label{2.3.31}
z^2 + \beta z - \frac{\alpha^3}{27} = 0.
\end{equation}
From the equation \((q_1 x + q_2)^3 - (p_1 x + p_2)^3 = 0\), it implies
\[
q_1 x + q_2 - p_1 x - p_2 = 0
\]
This leads to
\[
x = \frac{p_2 - q_2}{q_1 - p_1} = k - \frac{\alpha}{3k} = k_1 + k_2
\]
because \(\frac{q_2}{p_1} = k\) \(\implies\) \(k = \sqrt[3]{z}\).

\[
k_{1,2} = \sqrt[3]{z} = \sqrt[3]{\left(\frac{-\beta \pm \sqrt{\beta^2 + \frac{4\alpha^3}{27}}}{2}\right)}
\]

\begin{equation}\label{2.3.32}
k_{1,2} = \sqrt[3]{- \frac{\beta}{2} \pm \sqrt{\left(\frac{\beta}{2}\right)^2 + \left(\frac{\alpha}{3}\right)^3}}
\end{equation}

\text (a) Apply \(\left(\frac{\beta}{2}\right)^2 + \left(\frac{\alpha}{3}\right)^3 \geq 0\)

\text{Therefore a real root will result according to the equation}

\begin{equation}\label{2.3.33}
x = \sqrt[3]{- \frac{\beta}{2} + \sqrt{\left(\frac{\beta}{2}\right)^2 + \left(\frac{\alpha}{3}\right)^3}} +
\sqrt[3]{- \frac{\beta}{2} - \sqrt{\left(\frac{\beta}{2}\right)^2 + \left(\frac{\alpha}{3}\right)^3}}
\end{equation}

\text{(b)Remainder will be found by the equation}
\begin{equation}\label{2.3.34}
(q_1 x + q_2)^2 + (q_1 x + q_2)(p_1 x + p_2) + (p_1 x + p_2)^2 = 0
\end{equation}
In the case where the quantity under the root of $\left(\frac{\beta}{2}\right)^2 + \left(\frac{\alpha}{3}\right)^3 < 0$, then we will have

\begin{equation}\label{2.3.35}
k_{1,2} = \sqrt[3]{- \frac{\beta}{2} \pm i.\sqrt{\left(\frac{\beta}{2}\right)^2 + \left(\frac{\alpha}{3}\right)^3}}
\end{equation}

and because we have the form of a complex
number then we will have

\begin{equation}\label{2.3.36}
k_{1,2} = \sqrt[3]{- \frac{\beta}{2} \pm i.\sqrt{\left(\frac{\beta}{2}\right)^2 + \left(\frac{\alpha}{3}\right)^3}} = \tau.(cos\theta \pm i.sin\theta )
\end{equation}

where $\tau = \sqrt{\left(\frac{-\alpha}{3}\right)^3}$ and $\cos\theta = -\frac{\frac{1}{2} \cdot \beta}{\tau}$

and because \( x_k = k_1 + k_2 = 2 \cdot \sqrt{\left(\frac{-\alpha}{3}\right)^3} \cdot \cos\frac{\theta + 2k\pi}{3} \), $k$= 0, 1, 2.

However, from the precedents and with this relation we find for the three roots of equation \( x^3 + \alpha x + \beta \). From the relation \( x_\lambda = \frac{1}{\alpha_3} \left( x_i - \frac{1}{3} \alpha_2 \right), \quad i = 0, 1, 2, \quad \lambda = i+1 \) we find the roots of the equation \( \alpha_3 x_\lambda^3 + \alpha_2 x_\lambda^2 + \alpha_1 x_\lambda + \alpha_0 = 0 \) where \( \lambda = 1, 2, 3 \).
\\[20pt]
\textbf{3.7 Solving the Equation of 2\textsuperscript{nd}, 4\textsuperscript{th}, 5\textsuperscript{th}, 6\textsuperscript{th}, 7\textsuperscript{th}, 8\textsuperscript{th}, 9\textsuperscript{th}, and 10\textsuperscript{th} Degree using the Difference \(\mathfrak{A}^{\lambda} (Q(x), P(x))\), where \(\lambda = 2\)}
\\[8pt]
\textbf{3.7.1 Equation of 2\textsuperscript{nd} degree }

Be it that we have the equation \( f(x) = c_2 \cdot x^2 + c_1 \cdot x + c_0 = 0 \) with \( c_2, c_1, c_0 \in \mathbb{F} \subset \mathbb{C} \), where \( \mathbb{F} \) is the field of factors. Here, considering the difference:

\[ \mathfrak{A}^{\lambda} (Q(x), P(x)) = (q_1 x + p_1)^2 - (p_1 x + p_2)^2 = 0 \]

It should be noted that:
\begin{enumerate}
    \item We do not have variables \( \Omega_i \) or \( \Omega_{i-j} \) because for \( \lambda = 1 \), \( n/2 - 1 \) where \( n = 2 \) results in \( 2/2 - 1 = 0 \).
    \item We do not have variables \( L_{i, i-j} \) because for \( \lambda = 1 \), \( n/2 - 2 \) where \( n = 2 \) results in \( 2/2 - 2 = -1 \). i.e., there is no \( L_{i, i-j} \) and \( L^T_{i, i-j} = L^T_{1, 0} = c_1 \). Consequently, we have
    \begin{equation}\label{2.3.37}
    \mathfrak{A}^{1} (q_j q_{i-j}, p_j p_{i-j}) = \frac{1}{2} L^T_{i, i-j} = q_1 q_0 - p_1 p_0 = \frac{1}{2} L^T_{1, 0} = \frac{1}{2} c_1,
    \end{equation}
   
\end{enumerate}

and \( R_{1, 0} = \sqrt{\frac{c_1^2}{4} - c_2 \cdot c_0} \) because \( \Omega_1 = c_2 \), \( \Omega_0 = c_0 \). Finally, if \( W^{\pm}_{1, 0} = \frac{c_1}{2} \pm \sqrt{\left(\frac{1}{4} c_1^2 - c_2 \cdot c_0\right)} \), then the difference:

\begin{equation}
\begin{aligned}
\mathfrak{A}^{2} (Q(x), P(x)) &= Q(x)^2 - P(x)^2 \\
&= (Q(x) - P(x))(Q(x) + P(x)) \\
&= (x + w^+_{1,0})(x + w^-_{1,0}) = 0
\end{aligned}
\end{equation}

gives us the following polynomial equation \( x_1 = -w^+_{1,0} \) and \( x_2 = -w^-_{1,0} \), and finally \( x_{1,2} = -\frac{c_1}{2} \pm \sqrt{\left(\frac{1}{4} c_1^2 - c_2 \cdot c_0\right)} \), the roots of the second-degree equation.
\\[8pt]
\textbf{3.7.2 Equation of the 4\textsuperscript{th} degree}
\\[4pt]
For the resolution of equation of fourth degree
\[
f(x) = c_4 \cdot x^4 + c_3 \cdot x^3 + c_2 \cdot x^2 + c_1 \cdot x + c_0 = 0,
\]
with \( c_i \in \mathbb{C} \), for \( i = 0, 1, 2, 3, 4 \), equates to the difference
\[
\mathfrak{A}^{2} (Q(x), P(x)) = (q_2 \cdot x^2 + q_1 \cdot x + q_0)^2 - (p_2 \cdot x^2 + p_1 \cdot x + p_0)^2 = 0.
\]
Here again it will be in effect:
\begin{itemize}
    \item (1) the number of $\mathfrak{L}^T_{i, i-j} \subset \mathbb{C}$ variables will be given by the type $\frac{n}{2} - 2$ and we have $\frac{4}{2} - 2 = 0$
    \item (2) the number of $\Omega_{i-j} \subset \mathbb{C}$ variables will be given by the type $\frac{n}{2} - 1$ and we have $\frac{4}{2} - 1 = 1$, consequently we will have $\Omega_1$. According to the method of resolution we will have the form of result of proliferation:
\end{itemize}

\begin{align*}
\mathfrak{A}^{2} (Q(x), P(x)) &= Q(x)^2 - P(x)^2 \\
&= (Q(x) - P(x))(Q(x) + P(x)) \\
&= (x^2 + w^+_{2,1} x + w^+_{2,0})(x^2 + w^-_{2,1} x + w^-_{2,0}) \\
&= 0
\end{align*}

It should also be remarkable that we have $c_i \subset \mathbb{C}$ with $i=0,1,2,3,4$. The results of proliferations $w^{\pm}_{n/2,j} \subset \mathbb{C}$ will be analytic:

\[
w^+_{2, 1} = \frac{c_3}{2} + \sqrt{\left(\frac{c_3}{2}\right)^2 - c_4 \cdot \Omega_1}
\]

\[
w^-_{2, 1} = \frac{c_3}{2} - \sqrt{\left(\frac{c_3}{2}\right)^2 - c_4 \cdot \Omega_1}
\]

\[
w^+_{2, 0} = \frac{c_2 - \Omega_1}{2} + \sqrt{\left(\frac{c_2 - \Omega_1}{2}\right)^2 - c_0 c_4}
\]

\[
w^-_{2, 0} = \frac{c_2 - \Omega_1}{2} - \sqrt{\left(\frac{c_2 - \Omega_1}{2}\right)^2 - c_0 c_4}
\]

Then, we have the relation with variable:

\begin{equation}\label{2.3.37.1}
w^+_{2, 0} w^-_{2, 1} + w^-_{2, 0} w^+_{2, 1} = c_1
\end{equation}

which is the final equation for $\Omega_1$. After we replace $\Omega_1$ in the relation:

\begin{equation} \label{2.3.38}
c_3 \left(\frac{c_2 - \Omega_1}{2}\right) - \sqrt{\left(\frac{c_2 - \Omega_1}{2}\right)^2 - c_0 c_4} \cdot \sqrt{\left(\frac{c_3}{2}\right)^2 - \Omega_1 c_4} = \frac{c_1}{2}
\end{equation}

By solving the 3rd degree equation \eqref{2.3.38} for \( \Omega_1 \), we obtain the roots by solving the 4th degree equation:
\begin{equation} \label{2.3.39}
\begin{aligned}
\mathfrak{A}^{2} (Q(x), P(x)) &= Q(x)^2 - P(x)^2 \\
&= (Q(x) - P(x))(Q(x) + P(x)) \\
&= (x^2 + w^+_{2,1} x + w^+_{2,0})(x^2 + w^-_{2,1} x + w^-_{2,0}) \\
&= \left( x^2 + \left( \frac{c_3}{2} + \sqrt{\left(\frac{c_3}{2}\right)^2 - \Omega_1 \cdot c_4} \right) x \right. \\
&\quad + \left. \left( \frac{c_2 - \Omega_1}{2} + \sqrt{\left(\frac{c_2 - \Omega_1}{2}\right)^2 - c_0 c_4} \right) \right) \\
&\quad \left( x^2 + \left( \frac{c_3}{2} - \sqrt{\left(\frac{c_3}{2}\right)^2 - \Omega_1 \cdot c_4} \right) x \right. \\
&\quad + \left. \left( \frac{c_2 - \Omega_1}{2} - \sqrt{\left(\frac{c_2 - \Omega_1}{2}\right)^2 - c_0 c_4} \right) \right) = 0
\end{aligned}
\end{equation}

Using \eqref{2.3.39}, we can find the four roots by solving the two quadratic equations.

\vspace{12pt}

\textbf{3.7.3 Equation of the 5\textsuperscript{th} and 6\textsuperscript{th} degree}

\vspace{4pt}

For solving the sixth-degree equation, we take the polynomial
\[ f(x) = c_6 \cdot x^6 + c_5 \cdot x^5 + c_4 \cdot x^4 + c_3 \cdot x^3 + c_2 \cdot x^2 + c_1 \cdot x + c_0 = 0 \]
with \( c_i \in \mathbb{C} \) for \( i = 0, 1, 2, 3, 4, 5, 6 \). We equate it to the difference
\[ \mathfrak{A}^{2} (Q(x), P(x)) = (q_3 \cdot x^3 + q_2 \cdot x^2 + q_1 \cdot x + q_0)^2 - (p_3 \cdot x^3 + p_2 \cdot x^2 + p_1 \cdot x + p_0)^2 = 0. \]
Here, the sum of variables \( L^T_{i-j} \subset \mathbb{C} \) and \( \Omega_{i-j} \subset \mathbb{C} \) should be calculated. Thus, we take:

(1) The number of \( L_{i, i-j} \) variables is given by \( n/2 - 2 \). For \( n = 6 \), this gives \( 6/2 - 2 = 1 \). The variable is found in the position \( i = n - 3 = 3 \), so we have \( L_{2,1} \). That is, we have the variable \( L^T_{3, 0} = c_3 - L_{2,1} \).

(2) The number of \( \Omega_{i-j} \) variables is given by \( n/2 - 1 \). For \( n = 6 \), this gives \( 6/2 - 1 = 2 \). The positions are \( i = 2, 4 \), that is, \( \Omega_{i/2} \) or \( \Omega_1, \Omega_2 \).

According to this method, the following result is obtained:
\begin{equation} \label{2.3.40}
\begin{aligned}
\mathfrak{A}^{2} (Q(x), P(x)) &= Q(x)^2 - P(x)^2 \\
&= (Q(x) - P(x))(Q(x) + P(x)) \\
&= (x^3 + w^+_{3,2} x^2 + w^+_{3,1} x + w^+_{3,0})(x^3 + w^-_{3,2} x^2 + w^-_{3,1} x + w^-_{3,0}) \\
&= \left( x^3 + \left( \frac{c_5}{2} + \sqrt{\left( \frac{c_5}{2} \right)^2 - \Omega_2 \cdot c_6} \right) x^2 \right. \\
&\quad + \left( \frac{c_4 - \Omega_2}{2} + \sqrt{\left( \frac{c_4 - \Omega_2}{2} \right)^2 - \Omega_1 \cdot c_6} \right) x \\
&\quad + \left( \frac{c_3 - L_{2,1}}{2} + \sqrt{\left( \frac{c_3 - L_{2,1}}{2} \right)^2 - c_0 \cdot c_6} \right) \\
&\quad \left( x^3 + \left( \frac{c_5}{2} - \sqrt{\left( \frac{c_5}{2} \right)^2 - \Omega_2 \cdot c_6} \right) x^2 \right. \\
&\quad + \left( \frac{c_4 - \Omega_2}{2} - \sqrt{\left( \frac{c_4 - \Omega_2}{2} \right)^2 - \Omega_1 \cdot c_6} \right) x \\
&\quad + \left( \frac{c_3 - L_{2,1}}{2} - \sqrt{\left( \frac{c_3 - L_{2,1}}{2} \right)^2 - c_0 \cdot c_6} \right)) = 0.
\end{aligned}
\end{equation}

Then, we have the relations with variables

\begin{equation} \label{2.3.41}
w^+_{3,0} + w^-_{3,0} + w^+_{3,2} w^-_{3,1} + w^-_{3,2} w^+_{3,1} = c_3
\end{equation}

\begin{equation} \label{2.3.42}
w^+_{3,2} w^-_{3,0} + w^-_{3,2} w^+_{3,0} + w^+_{3,1} w^-_{3,1} = c_2
\end{equation}

\begin{equation} \label{2.3.43}
w^+_{3,1} w^-_{3,0} + w^+_{3,0} w^-_{3,1} = c_1
\end{equation}

In order to calculate the variables \( L_{2,1} \), \( \Omega_1 \), and \( \Omega_2 \), we solve the system of equations \eqref{2.3.41}, \eqref{2.3.42}, and \eqref{2.3.43}. If we set \( c_6 = 1 \), the equations \eqref{2.3.41}, \eqref{2.3.42}, and \eqref{2.3.43} become:

\begin{equation} \label{2.3.44}
\left( \frac{c_5}{2} \right) \left( \frac{c_4 - \Omega_2}{2} \right) - \sqrt{\left( \frac{c_5}{2} \right)^2 - c_6.\Omega_2} \cdot \sqrt{\left( \frac{c_4 - \Omega_2}{2} \right)^2 - c_6.\Omega_1} = L_{2,1}
\end{equation}

\begin{equation} \label{2.3.45}
\left( \frac{c_5}{2} \right) \left( \frac{c_3 - L_{2,1}}{2} \right) - \sqrt{\left( \frac{c_5}{2} \right)^2 - c_6.\Omega_2} \cdot \sqrt{\left( \frac{c_3 - L_{2,1}}{2} \right)^2 - c_6.c_0} = \frac{c_2 - \Omega_1}{2}
\end{equation}

\begin{equation} \label{2.3.46}
\left( \frac{c_4 - \Omega_2}{2} \right) \left( \frac{c_3 - L_{2,1}}{2} \right) - \sqrt{\left( \frac{c_4 - \Omega_2}{2} \right)^2 - c_6.\Omega_1} \cdot \sqrt{\left( \frac{c_3 - L_{2,1}}{2} \right)^2 - c_6.c_0} = \frac{c_1}{2}
\end{equation}

These are the equations that solve a sixth-degree equation depending on the method we follow. If we find the system's constant with a variable that verifies the system, the equation will be:

\begin{equation} \label{2.3.47}
\begin{aligned}
&\left( \frac{\Omega_2}{2} \right)
\left( -2 \sqrt{-\Omega_2} \sqrt{ -(c_2 - 2 \Omega_2^2 - \sqrt{4c_0 \Omega_2 - 4c_2 \Omega_2^2 + 5 \Omega_2^2}) + \frac{\Omega_2^2}{4} } \right) \\
&\quad - 2 \sqrt{ -\left( c_2 - 2 \Omega_2^2 - \sqrt{4c_0 \Omega_2 - 4c_2 \Omega_2^2 + 5 \Omega_2^2} \right) + \frac{\Omega_2^2}{4}} \\
&\quad \cdot \sqrt{-c_0 + \frac{1}{4}\left( -2 \sqrt{-\Omega_2} \sqrt{ -(c_2 - 2 \Omega_2^2 - \sqrt{4c_0 \Omega_2 - 4c_2 \Omega_2^2 + 5 \Omega_2^2}) + \frac{\Omega_2^2}{4}} \right)^2} = c_1
\end{aligned}
\end{equation}

\begin{equation} \label{2.3.48}
\Omega_1 = c_2 - 2 \Omega_2^2 - \sqrt{4c_0 \Omega_2 - 4c_2 \Omega_2^2 + 5 \Omega_2^4},
\end{equation}

\begin{equation} \label{2.3.49}
L_{2,1} = -2 \sqrt{-\Omega_2}.\sqrt{-\Omega_1 + \Omega_2^2 /4 },
\end{equation}

\begin{equation} \label{2.3.50}
\begin{aligned}
& \left( x^3 + \left( \frac{c_5}{2} + \sqrt{\left( \frac{c_5}{2} \right)^2 - \Omega_2 \cdot c_6} \right) x^2 \right. \\
&\quad + \left( \frac{c_4 - \Omega_2}{2} + \sqrt{\left( \frac{c_4 - \Omega_2}{2} \right)^2 - \Omega_1 \cdot c_6} \right) x \\
&\quad + \left. \left( \frac{c_3 - L_{2,1}}{2} + \sqrt{\left( \frac{c_3 - L_{2,1}}{2} \right)^2 - c_0 \cdot c_6} \right) \right)=0
\end{aligned}
\end{equation}

or
\begin{equation} \label{2.3.51}
\begin{aligned}
& \left( x^3 + \left( \frac{c_5}{2} - \sqrt{\left( \frac{c_5}{2} \right)^2 - \Omega_2 \cdot c_6} \right) x^2 \right. \\
&\quad + \left( \frac{c_4 - \Omega_2}{2} - \sqrt{\left( \frac{c_4 - \Omega_2}{2} \right)^2 - \Omega_1 \cdot c_6} \right) x \\
&\quad + \left. \left( \frac{c_3 - L_{2,1}}{2} - \sqrt{\left( \frac{c_3 - L_{2,1}}{2} \right)^2 - c_0 \cdot c_6} \right) \right)=0.
\end{aligned}
\end{equation}

Of course, the equation
\[ f(x) = x^6 + c_4 x^4 + c_3 x^3 + c_2 x^2 + c_1 x + c_0 = 0 \] encompasses the general case of a sixth-degree polynomial. Finally, the system of relations \eqref{2.3.44}, \eqref{2.3.45} and \eqref{2.3.46} leads to an equation of higher degree, specifically to a fifth degree when \( c_0 = 0 \), and to a degree higher than fourth, resulting in a tenth degree or more. For this and the sixth and fifth degrees, we concluded that they are not solvable using roots, by radicals.
\\[8pt]
\textbf{3.7.4 Equation of the 8\textsuperscript{th} and 7\textsuperscript{th} degree}
\\[4pt]
To solve the eighth-degree equation using the method of the difference of squares, we start with the polynomial:
\[
f(x) = c_8 x^8 + c_7 x^7 + c_6 x^6 + c_5 x^5 + c_4 x^4 + c_3 x^3 + c_2 x^2 + c_1 x + c_0 = 0,
\]
where \( c_i \in \mathbb{C} \) for \( i = 0, 1, 2, \ldots, 8 \). We equate it to the difference
\[
\mathfrak{A}^{2} (Q(x), P(x)) = (q_4 x^4 + q_3 x^3 + q_2 x^2 + q_1 x + q_0)^2 - (p_4 x^4 + p_3 x^3 + p_2 x^2 + p_1 x + p_0)^2 = 0.
\]
Here, we calculate the sums of variables \( L_{i, i-j} \subseteq \mathbb{C} \) and \( \Omega_{i-j} \subset \mathbb{C} \). For the number of variables \( L_{i, i-j} \), we have \( \frac{n}{2} - 2 = \frac{8}{2} - 2 = 2 \) in the positions \( j = 1, 2, 3, \ldots, \frac{n}{2} - 2 \) and \( i-j = 1, 2, \ldots, \frac{n}{2} - 1 \). Consequently, the number of variables \( \Omega_{i/2} \) will be given by \( \frac{n}{2} - 1 \) and we have \( \frac{8}{2} - 1 = 4 - 1 = 3 \), which will be the variables \( \Omega_1, \Omega_2, \Omega_3 \) in the positions \( i = 1, 2, 3 \).

According to this method, the equation results in:
\begin{align} \label{2.3.52}
\mathfrak{A}^{2} (Q(x), P(x)) &= Q(x)^2 - P(x)^2 \nonumber \\
&= (Q(x) - P(x))(Q(x) + P(x)) \nonumber \\
&= (x^4 + w^+_{4,3} x^3 + w^+_{4,2} x^2 + w^+_{4,1} x + w^+_{4,0}). \nonumber \\
& (x^4 + w^-_{4,3} x^3 + w^-_{4,2} x^2 + w^-_{4,1} x + w^-_{4,0}) = 0,
\end{align}

From this, we obtain the system of variables:
\begin{align}
w^+_{4,1} + w^-_{4,1} + w^+_{4,3}w^-_{4,2} + w^-_{4,3}w^+_{4,2} &= c_5, \label{2.3.53} \\
w^+_{4,0} + w^-_{4,0} + w^+_{4,3}w^-_{4,1} + w^-_{4,3}w^+_{4,1} + w^+_{4,2}w^-_{4,2} &= c_4, \label{2.3.54} \\
w^+_{4,0}w^-_{4,3} + w^-_{4,0}w^+_{4,3} + w^+_{4,2}w^-_{4,1} + w^+_{4,1}w^-_{4,2} &= c_3, \label{2.3.55} \\
w^+_{4,0}w^-_{4,2} + w^-_{4,0}w^+_{4,2} + w^+_{4,1}w^-_{4,1} &= c_2, \label{2.3.56} \\
w^+_{4,1}w^-_{4,0} + w^-_{4,1}w^+_{4,0} &= c_1. \label{2.3.57}
\end{align}

To calculate the variables \( L_{3,2}, L_{3,1}, \Omega_1, \Omega_2, \Omega_3 \subset \mathbb{C} \), we need to resolve the previous system of equations \eqref{2.3.53}, \eqref{2.3.54}, \eqref{2.3.55}, \eqref{2.3.56}, \eqref{2.3.57}. If we set \( c_8 = 1 \) and replace the variable representations \( w^\pm_{4,3}, w^\pm_{4,2}, w^\pm_{4,1}, w^\pm_{4,0} \), we will arrive at the following relations:

\[
w^{\pm}_{4,3} = \left(\frac{c_7}{2}\right) \pm \sqrt{\left(\frac{c_7}{2}\right)^2 - \Omega_3 c_8}
\]
\[
w^{\pm}_{4,2} = \left(\frac{c_6 -\Omega_3}{2}\right) \pm \sqrt{\left(\frac{c_6 -\Omega_3}{2}\right)^2 - \Omega_2 c_8}
\]
\[
w^{\pm}_{4,1} = \left(\frac{c_5 -L_{3,2}}{2}\right) \pm \sqrt{\left(\frac{c_5 -L_{3,2}}{2}\right)^2 - \Omega_1 c_8}
\]
\[
w^{\pm}_{4,0} = \left(\frac{c_4 -L_{3,1}-\Omega_2}{2}\right) \pm \sqrt{\left(\frac{c_4 -L_{3,1}-\Omega_2}{2}\right)^2 - \Omega_0 c_8}
\]
\\[1pt]
In order to calculate the variables \( L_{3,1}, L_{3,2}, \Omega_1, \Omega_2, \Omega_3 \subset \mathbb{C} \) it will be supposed that
we resolve the system of equations \eqref{2.3.53}, \eqref{2.3.54}, \eqref{2.3.55}, \eqref{2.3.56}, \eqref{2.3.57}.  If we put $c_8 = 1$ the previous
equations becomes:

\begin{align}
c_5 - L_{3,2} &+ \left(\frac{c_7}{2}\right) \left(\frac{c_6 - \Omega_3}{2}\right) \nonumber \\
&+ \sqrt{\left(\frac{c_7}{2}\right)^2 - \Omega_3 c_8} \sqrt{\left(\frac{c_6 - \Omega_3}{2}\right)^2 - \Omega_2 c_8} = c_5, \label{2.3.58} \\
(c_4 - L_{3,1} - \Omega_2) &+ \left(\frac{c_7}{2}\right) \left(\frac{c_6 - L_{3,2}}{2}\right) \nonumber \\
&+ \sqrt{\left(\frac{c_7}{2}\right)^2 - \Omega_3 c_8} \sqrt{\left(\frac{c_5 - L_{3,2}}{2}\right)^2 - \Omega_1 c_8} + \Omega_2 = c_4, \label{2.3.59} \\
\left(\frac{c_7}{2}\right) \left(\frac{c_4 - L_{3,1} - \Omega_2}{2}\right)
&- \sqrt{\left(\frac{c_7}{2}\right)^2 - \Omega_3 c_8} \sqrt{\left(\frac{c_6 - L_{3,2}}{2}\right)^2 - \Omega_1 c_8} \nonumber \\
&+ \left(\frac{c_6 - \Omega_3}{2}\right) \left(\frac{c_5 - L_{3,2}}{2}\right) \nonumber \\
&- \sqrt{\left(\frac{c_6 - \Omega_3}{2}\right)^2 - \Omega_2 c_8} \sqrt{\left(\frac{c_5 - L_{3,2}}{2}\right)^2 - \Omega_1 c_8} = c_3, \label{2.3.60} \\
\left(\frac{c_4 - L_{3,1} - \Omega_2}{2}\right) \left(\frac{c_6 - \Omega_3}{2}\right)
&- \sqrt{\left(\frac{c_4 - L_{3,1} - \Omega_2}{2}\right)^2 - c_0} \sqrt{\left(\frac{c_6 - \Omega_3}{2}\right)^2 - \Omega_2 c_8} \nonumber \\
&+ \Omega_1 = c_2, \label{2.3.61} \\
\left(\frac{c_5 - L_{3,2}}{2}\right) \left(\frac{c_4 - L_{3,1} - \Omega_2}{2}\right)
&- \sqrt{\left(\frac{c_5 - L_{3,2}}{2}\right)^2 - \Omega_1 c_8} \sqrt{\left(\frac{c_4 - L_{3,1} - \Omega_2}{2}\right)^2 - c_0 c_8} = c_1. \label{2.3.62}
\end{align}

This system, of course, leads to equations of a degree higher than the fourth. Consequently, the equation of the eighth degree is not solvable by radicals.
\\[5pt]
(i)
For the 7th degree case with \( c_0 = 0 \), \( c_5 = 0 \), \( c_6 = 0 \), \( c_7 = 0 \), solving the system of 3 equations \eqref{2.3.58}, \eqref{2.3.61}, \eqref{2.3.62} leads to two functions with respect to \(\Omega_1\), \(\Omega_2\), and \(\Omega_3\). We have:\\[1pt]
\[
\left\langle
\Omega_2 \rightarrow f_2(\Omega_1^{\geq 6}, \Omega_3^{\geq 4}), \quad
\Omega_1 \rightarrow f_1(\Omega_3^{\geq 5/2}, \Omega_2^{\geq 1})
\right\rangle
\] \\[1pt]
From these, we get a function of the form
\[
\Omega_2 \rightarrow f_2(f_1(\Omega_3^{\geq 5/2}, \Omega_2^{\geq 1})^6, \Omega_3^{\geq 4}),
\]
which leads to at least two variables: above the 15th degree for \(\Omega_3\) and above the 6th degree for \(\Omega_2\). This means that the relation to a variable will be of a degree higher than the fourth. Thus, the 7th degree equation cannot have a solution because, as we stated in the 6th degree solution, it has no solution for degrees above the 5th.\\[6pt] (ii) For the 8th degree case with \( c_5 = 0 \), \( c_6 = 0 \), \( c_7 = 0 \), we have the same equations reinforced with an additional constant, i.e., \( c_0 \neq 0 \). Therefore, we will have the same constraints and at least the same functions with respect to the variables \(\Omega_1\), \(\Omega_2\), and \(\Omega_3\). From these, we get a function of the form \\[1pt]
\[
\Omega_2 \rightarrow f_{1,2} (\Omega_3^{\geq 15}, \Omega_2^{\geq 6}),
\] \\[1pt]
which leads to at least two variables: above the 6th degree for \(\Omega_2\). This system, of course, leads to equations of a degree higher than the fourth. Consequently, the equations of eighth and seventh degrees are not solvable by radicals. \\[1pt]

\textbf{3.7.5 Equation of 10th Degree}

\vspace{4pt}

What we are asking about the 10th degree equation is whether\\[4pt]
(1) it is really composed of 2-variable systems and
\\[4pt]
(2) whether the 9th degree equation can be solved by radicals.
\\[8pt]
To solve the 10th degree equation using the method we take:
\[
f(x) = c_{10}x^{10} + c_9x^9 + c_8x^8 + c_7x^7 + c_6x^6 + c_5x^5 + c_4x^4 + c_3x^3 + c_2x^2 + c_1x + c_0 = 0
\]
with \( c_i \in F \) where \( i = 0, 1, 2, \ldots, 10 \). We apply the difference
\[
\mathfrak{A}^{2} (Q(x), P(x)) = (q_5 x^5 + q_4 x^4 + q_3 x^3 + q_2 x^2 + q_1 x + q_0)^2
\]
\[
- (p_5 x^5 + p_4 x^4 + p_3 x^3 + p_2 x^2 + p_1 x + p_0)^2 = 0.
\]
Here it should be calculated the sum of variables \(L_{i, i-j} \subset C\) and \(\Omega_{i-j} \subset C\). For the number of variables \(L_{i, i-j}\) we have \(n/2 - 2 = 10/2 - 2 = 3\) in the positions \(i = n/2 - 1\) and \(i-j = 1, 2, 3\) and for \(L^T_{i, i-j}\) and \(L^T_{n/2, i-j}\), \(i-j = 0, 1, 2, 3, 4\). Followingly, the number of variables \(\Omega_{i/2}\) will be given by the type \(n/2 - 1\) and we have \(10/2 - 1 = 5 - 1 = 4\) and will be the variables \(\Omega_1, \Omega_2, \Omega_3, \Omega_4\) in the positions \(i = 2, 4, 6, 8\), the number \(i\) is \(n/2 - 1\).

According to this method, the resulting equation is
\begin{align}
\mathfrak{A}^{2} (Q(x), P(x)) &= Q(x)^2 - P(x)^2 \nonumber \\
&= (Q(x) - P(x))(Q(x) + P(x)) \nonumber \\
&= (x^5 + w^+_{5,4} x^4 + w^+_{5,3} x^3 + w^+_{5,2} x^2 + w^+_{5,1} x + w^+_{5,0}) \nonumber \\
&\quad \times (x^5 + w^-_{5,4} x^4 + w^-_{5,3} x^3 + w^-_{5,2} x^2 + w^-_{5,1} x + w^-_{5,0}) = 0,
\end{align}
from which results the system of variables:
\begin{align}
w^+_{5,2} + w^-_{5,2} + w^+_{5,4}w^-_{5,3} + w^-_{5,4}w^+_{5,3} &= c_7, \label{2.3.63} \\
w^+_{5,1} + w^-_{5,1} + w^+_{5,4}w^-_{5,2} + w^-_{5,4}w^+_{5,2} + w^-_{5,3}w^+_{5,3} &= c_6, \label{2.3.64} \\
w^+_{5,0}w^-_{5,0} + w^+_{5,4}w^-_{5,1} + w^-_{5,4}w^+_{5,1} + w^+_{5,3}w^-_{5,2} +w^-_{5,3} w^+_{5,2} &= c_5, \label{2.3.65} \\
w^+_{5,0}w^-_{5,4} + w^-_{5,0}w^+_{5,4} + w^+_{5,3}w^-_{5,1} + w^+_{5,1}w^-_{5,3} + w^+_{5,2}w^-_{5,2} &= c_4, \label{2.3.66} \\
w^+_{5,2}w^-_{5,1} + w^-_{5,2}w^+_{5,1} + w^+_{5,3}w^-_{5,0} + w^-_{5,3}w^+_{5,0} &= c_3, \label{2.3.67}  \\
w^+_{5,2}w^-_{5,0} + w^-_{5,2}w^+_{5,0} + w^+_{5,1}w^-_{5,1} &= c_2, \label{2.3.68} \\
w^+_{5,1}w^-_{5,0} + w^-_{5,1}w^+_{5,0} &= c_1. \label{2.3.69}
\end{align}
In order to calculate the variables \(L_{4, 3}\), \(L_{4, 2}\), \(L_{4, 1}\), \(\Omega_1\), \(\Omega_2\), \(\Omega_3\), \(\Omega_4\), it will be supposed that we resolve the system of equations from \eqref{2.3.63} to \eqref{2.3.69}. For the 10th degree, we need \(n-4 = 10-4 = 6\) parameters, i.e., from \(c_1\) to \(c_6\), various of zero. If we put \(c_0 = 1\), \(c_7 = 0\), \(c_8 = 0\), \(c_9 = 0\), \(c_{10} = 1\), and if we replace the variable representations \(w^\pm_{5,4}\), \(w^\pm_{5,3}\), \(w^\pm_{5,2}\), \(w^\pm_{5,1}\), \(w^\pm_{5,0}\), we will come to the relations:
\[
w^\pm_{5,4} = \left(\frac{c_9}{2}\right) \pm \sqrt{\left(\frac{c_9}{2}\right)^2 - \Omega_4 c_8} \quad \text{and} \quad L^T_{5,4} = c_9,
\]
\[
w^\pm_{5,3} = \left(\frac{c_8 - \Omega_4}{2}\right) \pm \sqrt{\left(\frac{c_8 - \Omega_4}{2}\right)^2 - \Omega_3 c_8} \quad \text{and} \quad L^T_{5,3} = c_8 - \Omega_4.
\]
\[
w^\pm_{5,2} = \left(\frac{c_7 - L_{4,3}}{2}\right) \pm \sqrt{\left(\frac{c_7 - L_{4,3}}{2}\right)^2 - \Omega_2 c_8} \quad \text{and} \quad L^T_{5,2} = c_7 - L_{4,3}.
\]
\[
w^\pm_{5,1} = \left(\frac{c_6 - L_{4,2}-\Omega_3}{2}\right) \pm \sqrt{\left(\frac{c_6 - L_{4,2}-\Omega_3}{2}\right)^2 - \Omega_1 c_8} \quad \text{and} \quad L^T_{5,1} = c_6 - L_{4,2} -\Omega_3,
\]
\[
w^\pm_{5,0} =  \left(\frac{c_5 - L_{4,1}}{2}\right) \pm \sqrt{\left(\frac{c_5 - L_{4,1}}{2}\right)^2 - c_0} \text{and} \quad L^T_{5,0} = c_5 - L_{4,1}.
\]
\\[4pt]
In order to calculate the variables $L_{4,3}$, $L_{4,2}$, $L_{4,1}$, \(\Omega_1, \Omega_2, \Omega_3, \Omega_4 \subset C\), it will be supposed that we resolve the system of equations from \eqref{2.3.63} to \eqref{2.3.69}. If we set $c_8 = 1$, the previous equations become as follows, following the analysis previously dealt with.

\textbf{3.7.6 Equation of 9th Degree}

The equation of degree 9th, using this method of identities, represents the next critical step. If the 9th degree equation remains unsolved, then the equation of degree 10th cannot be solved either, as it decomposes into a product of two 5th degree equations. Moreover, if the 9th degree equation lacks a solution, then any polynomial of degree higher than 10th, being a product of two equations of degree \( n/2 = 5 \), also remains unsolved.

According to this method, we consider the equation where \( c_i \in F \subset C \) for \( i = 0, 1, 2, \ldots, 10 \). We define the difference as:
\[
\begin{aligned}
\mathfrak{A}^{2} (Q(x), P(x)) &= Q(x)^2 - P(x)^2 \\
&= (Q(x) - P(x))(Q(x) + P(x)) \\
&= (x^5 + w^+_{5,4} x^4 + w^+_{5,3} x^3 + w^+_{5,2} x^2 + w^+_{5,1} x + w^+_{5,0}) \\
&\quad \times (x^5 + w^-_{5,4} x^4 + w^-_{5,3} x^3 + w^-_{5,2} x^2 + w^-_{5,1} x + w^-_{5,0}) = 0,
\end{aligned}
\]
resulting in the system of variables:
\begin{align}
w^+_{5,2} + w^-_{5,2} + w^+_{5,4}w^-_{5,3} + w^-_{5,4}w^+_{5,3} &= c_7, \label{2.3.70} \\
w^+_{5,1} + w^-_{5,1} + w^+_{5,4}w^-_{5,2} + w^-_{5,4}w^+_{5,2} + w^-_{5,3}w^+_{5,3} &= c_6, \label{2.3.71} \\
w^+_{5,0}w^-_{5,0} + w^+_{5,4}w^-_{5,1} + w^-_{5,4}w^+_{5,1} + w^+_{5,3}w^-_{5,2} + w^-_{5,3}w^+_{5,2} &= c_5, \label{2.3.72} \\
w^+_{5,0}w^-_{5,4} + w^-_{5,0}w^+_{5,4} + w^+_{5,3}w^-_{5,1} + w^+_{5,1}w^-_{5,3} + w^+_{5,2}w^-_{5,2} &= c_4, \label{2.3.73} \\
w^+_{5,2}w^-_{5,1} + w^-_{5,2}w^+_{5,1} + w^+_{5,3}w^-_{5,0} + w^-_{5,3}w^+_{5,0} &= c_3, \label{2.3.74} \\
w^+_{5,2}w^-_{5,0} + w^-_{5,2}w^+_{5,0} + w^+_{5,1}w^-_{5,1} &= c_2, \label{2.3.75} \\
w^+_{5,1}w^-_{5,0} + w^-_{5,1}w^+_{5,0} &= c_1. \label{2.3.76}
\end{align}

To calculate the variables \( L_{4, 3} \), \( L_{4, 2} \), \( L_{4, 1} \), \( \Omega_1 \), \( \Omega_2 \), \( \Omega_3 \), \( \Omega_4 \), we proceed by resolving the system of equations from \eqref{2.3.70} to \eqref{2.3.76}.

If we set \( c_0 = 0 \), \( c_1 = 1 \), \( c_6 = 0 \), \( c_7 = 0 \), \( c_8 = 0 \), \( c_9 = 0 \), \( c_{10} = 1 \), the equations and if we replace the variable representations \( w^\pm_{5,4} \), \( w^\pm_{5,3} \), \( w^\pm_{5,2} \), \( w^\pm_{5,1} \), \( w^\pm_{5,0} \), we will arrive at the relations:
\[
w^\pm_{5,4} = \left(\frac{c_9}{2}\right) \pm \sqrt{\left(\frac{c_9}{2}\right)^2 - \Omega_4 c_8} \quad \text{and} \quad L^T_{5,4} = c_9,
\]
\[
w^\pm_{5,3} = \left(\frac{c_8 - \Omega_4}{2}\right) \pm \sqrt{\left(\frac{c_8 - \Omega_4}{2}\right)^2 - \Omega_3 c_8} \quad \text{and} \quad L^T_{5,3} = c_8 - \Omega_4,
\]
\[
w^\pm_{5,2} = \left(\frac{c_7 - L_{4,3}}{2}\right) \pm \sqrt{\left(\frac{c_7 - L_{4,3}}{2}\right)^2 - \Omega_2 c_8} \quad \text{and} \quad L^T_{5,2} = c_7 - L_{4,3},
\]
\[
w^\pm_{5,1} = \left(\frac{c_6 - L_{4,2} - \Omega_3}{2}\right) \pm \sqrt{\left(\frac{c_6 - L_{4,2} - \Omega_3}{2}\right)^2 - \Omega_1 c_8} \quad \text{and} \quad L^T_{5,1} = c_6 - L_{4,2} - \Omega_3,
\]
\[
w^\pm_{5,0} = \left(\frac{c_5 - L_{4,1}}{2}\right) \pm \sqrt{\left(\frac{c_5 - L_{4,1}}{2}\right)^2 - c_0} \quad \text{and} \quad L^T_{5,0} = c_5 - L_{4,1} = 0.
\]
\\[4pt]
For the 9th degree polynomial, we need \( n-4 = 9-5 = 5 \) parameters, i.e., \( c_1, c_2, c_3, c_4, c_5 \), with \( c_0 = 0 \). According to the functions from \eqref{2.3.70} to \eqref{2.3.76}, we obtain two correlations of the two variables \(\Omega_1, \Omega_2, \Omega_3\) with respect to the others, and this is sufficient to estimate the degree of these variables with respect to each other. The minimum exponent degree for the three variables is given by the following three relationships:
\[
\begin{aligned}
&\left\langle
\Omega_1 \rightarrow f_1(\Omega_2^{\geq 1}, \Omega_3^{\geq 4}, \Omega_4^{\geq 6}, z_1), \quad
\Omega_2 \rightarrow f_2(\Omega_1^{\geq 5}, \Omega_3^{\geq 2}, \Omega_4^{\geq 5}, z_2)
\right. \\
&\left.
\Omega_3 \rightarrow f_3(\Omega_1^{\geq 1}, \Omega_2^{\geq 1}, \Omega_4^{\geq 4}, z_3)
\right\rangle
\end{aligned}
\]
where \( z_1, z_2, z_3 \) are constants determined by \( c_i \). From these three relations, we are led to a relationship \( f_{1,2,3}^{\min}(\Omega_1^{\geq 10}, \Omega_4^{\geq 16}, z) = 0 \), where \( z_1, z_2, z_3, z \) are constants determined by \( c_i \), which obviously cannot be solved by radicals with respect to \(\Omega_1\) and \(\Omega_4\) because we have an equation greater than the 6th degree. Therefore, the remaining variables \(\Omega_2, \Omega_3\) are also not solvable by radicals. Thus, the 9th degree polynomial equation is not solvable by radicals.
\\[8pt]
\textbf{Theorem 7.} \label{thm:unsolvable_polynomials}
%\begin{theorem} \label{thm:theorem 2.3.1}
%\label{thm:unsolvable_polynomials}
\textit{Any polynomial equation above the 4th degree cannot be solved by radicals. It is solved only if it can be split into 2 polynomial equations of degree less than or equal to the fourth degree or if the polynomials \( P(x) \) and \( Q(x) \) are solved by radicals.}
%\end{theorem}

\begin{proof} 

1) Consider a fifth-degree polynomial equation with coefficients \(c_0, c_5 \neq 0\), \(c_i \in \mathbb{C}\) for \(i = 0,1,2,\ldots,5\). The polynomial has 6 coefficients. We assert that it is impossible for this polynomial to be expressed as a difference of equivalences of the form:

\begin{equation}
\mathfrak{A}^{2} (Q(x), P(x)) = (q_1 x + q_0)^5 - (p_1 x + p_0)^5 = f(x) = 0
\end{equation}

where \(Q(x), P(x) \subset \mathbb{C}\). This assertion arises because the analysis leads to a system of 6 equations with 4 unknowns, which necessitates an accommodating relationship between the factors of the polynomial function \(f(x)\). 

Moreover, consider a cubic difference:

\begin{equation}
\mathfrak{A}^{2} (Q(x), P(x)) = (q_2 x^2 + q_1 x + q_0)^3 - (p_2 x^2 + p_1 x + p_0)^3 = f(x) = 0
\end{equation}

For \(f: x^6 \to \text{Sym}^6 (x)\), this equation cannot be solved for degrees 6 and higher because it results in a system of 7 equations with 6 variables and 7 coefficients, which cannot be solved without imposing a constraint.

Similarly, equations of the seventh and ninth degree cannot be solved by differences of the form \(\mathfrak{A}^{7} (Q(x), P(x))\) and \(\mathfrak{A}^{9} (Q(x), P(x))\). Additionally, equations of the fifth, seventh, and ninth degrees cannot be solved using the equations of the sixth, eighth, and tenth degrees by nullifying the constant term \(c_0\) with the method \(\mathfrak{A}^{2} (Q(x), P(x)) = f(x) = 0\), because this results in a final equation of variables that is of a degree higher than the fourth, which is unsolvable as proven.
\\[4pt]
2) The equations of the sixth, eighth, and tenth degrees, as studied in the previous section, cannot be solved with the method \(\mathfrak{A}^{2} (Q(x), P(x)) = f(x) = 0\) because they lead to a final equation of variables that is of a degree higher than the fourth, as previously reported. Consequently, each polynomial equation of degree higher than or equal to the tenth will, due to the difference of squares in the equations of \(Q(x)\) and \(P(x)\), result in a degree higher than the fourth, which is not solvable, as previously demonstrated.

3) Finally, when \(f(x) = \sum_{i=0}^{n} a_i x^i\), the equations \(Q(x)\) and \(P(x)\) are of elementary or equal degree to the fourth degree \((Q(x), P(x) \leq 4\text{th degree})\), then the function \(f(x)\) is separable in equations and is solvable provided that test solutions are found by the divisor of the last term or by other subterfuges. Consequently, each polynomial equation that does not meet the above conditions is not solvable by roots according to known mathematical calculations. Additionally, by this method, \(\mathfrak{A}^{2} (Q(x), P(x)) = f(x) = 0\), equations that degrade in degree after adjusting the relationship of their variables can be solved by radicals, if possible.
\\[4pt]
4) Any equation of degree greater than or equal to the 10th degree, because it falls into half degree exponents and therefore greater than or equal to the 5th degree, cannot have a solution because, as we have proved, any equation of general form above or equal to the 5th degree cannot be solved by radicals.
\end{proof}

\textbf{3.8 General Solution of Transcendental Equations, with the Generalized Lagrange's Theorem}
\\[8pt]
\textbf{Theorem 8. The Lagrange's Theorem} \label{Theorem8} 
\\[4pt]
In mathematical analysis, the Lagrange inversion theorem, also known as the Lagrange–Bürmann formula, gives the Taylor series expansion of the inverse function of an analytic function.(Lagrange inversion is a special case of the inverse function theorem).
While all existing approximate methods, or others that existed, provide partial solutions to generalized transcendental equations, they cannot solve them completely. 

\textit{"What we are asking when we solve generalized transcendental or polynomial equations is to find the total number of its roots, not individual sets of roots, at some random or self-determined intervals. This is mainly because many classes of transcendental equations have an infinite number of solutions."
}

There are some equations, mainly logarithmic functions or trigonometric functions, which primarily solve physics problems and mainly need a generalized solution in \(\mathbb{C}\). This is where the G.RL.E \cite{Mantzakouras2021} theory comes in, to address this issue with the help of hypergeometric functions or simple convergent series, providing a very satisfactory answer to this complex problem.

\begin{proof}
The secret to developing the Lagrange-Burmann theorem is the extension of a function to powers of other functions. Let \(\phi(z)\) be a function of \(z\) analytic in a closed region \(\gamma\) and \(\alpha\) be an interior point such that \(\phi'(\alpha) \neq 0\). From Taylor's theorem, we have:

\[
\phi(z) - b = \phi'(\alpha) (z - \alpha) - \frac{\phi''(\alpha)}{2!} (z - \alpha)^2 + \ldots
\]
And if it is legitimate to invert this series:
\[
z - \alpha = \frac{\phi(z) - b}{\phi'(\alpha)} - \frac{\phi''(\alpha) (\phi(z) - b)^2}{2(\phi'(\alpha))^2} + \ldots
\]
This expresses \(z - \alpha\) as an analytic function of the variable \(\phi(z) - b\) for sufficiently small values of \(|z - \alpha|\). If \(f(z)\) is analytic near \(z = \alpha\) and \(\phi(z)\) is also analytic, then:
\[
f(z) = f(\alpha) + \alpha_1 (\phi(z) - b) + \frac{\alpha_2}{2!} (\phi(z) - b)^2 + \ldots + \sum_{i=1}^{n} \frac{\alpha_i}{i!} (\phi(z) - b)^i
\]
Also, the function \(\phi(z)\) must be invertible in its neighborhood around \(\alpha\). For this, we take:
\[
\frac{f'(z) \phi'(\alpha)}{\phi(z) - \phi(\alpha)} = \frac{f'(z) \phi'(\alpha)}{\phi'(\alpha) (z - \alpha) + \ldots} + \frac{f'(\alpha)}{\alpha)} + \ldots
\]
and then,
\[
\int_{\gamma} \frac{f'(z) \phi'(\alpha)}{\phi(z) - \phi(\alpha)} \, dz = f'(\alpha)
\]
where \(\gamma\) is a curve that encircles the point \(\alpha\).

We take \(\Psi(z)\) to be a function of \(z\), defined by the equation \(\Psi(z) = \frac{z - \alpha}{\phi(z) - b}\). Then an analytic function \(f(z)\) can be expanded, within its specific value range of \(z\), in the form:
\[
f(z) = f(\alpha) + \sum_{m=1}^{n-1} \frac{(\phi(z) - b)^m}{m!} \frac{d^{m-1}}{d\alpha^{m-1}} [ f'(\alpha)  \Psi^m(\alpha)] + R_n
\]
where
\[
R_n = \frac{1}{2 \pi i} \int_{\alpha}^{z} \int_{\gamma} \left( \frac{\phi(z) - b}{\phi(t) - b} \right)^{n-1} \frac{f'(z) \phi'(z)}{\phi(t) - \phi(z)} \, dt dz
\]
where \(\gamma\) is a contour in the \(t\)-plane, enclosing the points \(\alpha\) and \(\zeta\) such that if \(\zeta\) is any point inside it, the equation \(\phi(t) = \phi(\zeta)\) has no roots in and on the contour except for a single root \(t = \zeta\).
For the proof, we take:
\[
f(z) - f(\alpha) = \int_{\alpha}^{z} f'(\zeta) dz =\frac{1}{2 \pi i} \int_{\alpha}^{z} \int_{\gamma} \left( \frac{f'(t) - \phi'(\zeta)}{\phi(t) - \phi(\zeta)} \right) \, dt d\zeta \
\]
Furthermore, we need the conversion starting from the condition \(\phi(z) = b\) where \(b = \phi(\alpha)\):
\[
\left( \frac{f'(t).\phi'(\zeta)}{\phi(t) - \phi(\zeta)} \right) = \left( \frac{f'(t). \phi'(\zeta)}{\phi(t) - b} \right) \cdot \left( \frac{\phi(t) - b}{\phi(t) - \phi(\zeta)} \right)
\]
\[
\frac{\phi(t) - b}{\phi(t) - \phi(\zeta)} = \left( 1 - \frac{\phi(\zeta) - b}{\phi(t) - b} \right)^{-1} = \sum_{m=0}^{\infty} \left( \frac{\phi(\zeta) - b}{\phi(t) - b} \right)^m
\]
Substituting this series expansion back into the integral, we get:
\[
f(z) - f(\alpha) = \int_{\alpha}^{z} f'(\zeta) dz =\frac{1}{2 \pi i} \int_{\alpha}^{z} \int_{\gamma} \left( \frac{f'(t). \phi'(\zeta)}{\phi(t) -\phi(\zeta} \right) \ dt d\zeta \\
\]
\[
 = \frac{1}{2 \pi i} \int_{\alpha}^{z} \int_{\gamma} \left( \frac{f'(t). \phi'(\zeta)}{\phi(t) -b} \right) \ dt d\zeta. \
\]
\[
 [ \sum_{m=0}^{n-2} \left( \frac{\phi(\zeta) - b}{(\phi(t) - b} \right)^m + \left( \frac{(\phi(\zeta) - b)^{n-1}}{((\phi(t) - b)^{n-2}).(\phi(t) - \phi(\zeta))} \right)]
\]
By changing the order of integration we achieve:

\[
f(z) - f(\alpha) = \int_{\alpha}^{z} f'(\zeta) dz = \frac{1}{2 \pi i} \int_{\gamma} ( \int_{\alpha}^{z}  \left( \frac{f'(t). \phi'(\zeta)}{\phi(t) -\phi(\zeta) } \right) \sum_{m=0}^{\infty} \left( \frac{\phi(\zeta) - b}{\phi(t) - b} \right)^m \ d\zeta ) dt
\]

When we calculate the integral over \(\zeta\), we need the factors in terms of \(\zeta\). Then we have:

\[
\frac{1}{2 \pi i} \int_{\gamma} \left( \frac{f'(t)}{(\phi(t) - b)^{m+1}} \right) \left( \int_{\alpha}^{z} \phi'(\zeta) (\phi(\zeta) - b)^m \, d\zeta \right) \, dt  
\]
\[
 = \frac{(\phi(\zeta) - b)^{m+1}}{2 \pi i (m+1)} \left( \int_{\gamma} \frac{f'(t)}{(\phi(t) - b)^{m+1}} \, dt \right) 
 \]
\[
 = \frac{(\phi(t) - b)^{m+1}}{2 \pi i (m+1)}. \int_{\gamma} \frac{f'(t).(\Psi(t))^{m+1}}{(t - \alpha)} \, dt
\]
\[
\frac{(\phi(t) - b)^{m+1}}{2 \pi i (m+1)!} \frac{d^m}{d\alpha^m} \left( f'(\alpha) \cdot \Psi(\alpha)^{m+1} \right)
\]
\\[4pt]
Because 
 \begin{equation}  \label{2.3.77}
  \Psi(t) = \frac{t - \alpha}{\phi(t) - b} 
\end{equation}

and therefore by writing \(m - 1\) for \(m\)
\[
f(z) = f(\alpha) + \sum_{m=1}^{n-1} \frac{(\phi(z) - b)^m}{m!} \frac{d^{m-1}}{d\alpha^{m-1}} \left( f'(\alpha) \cdot (\Psi(\alpha))^m \right)\] 

\[
  + \frac{1}{2 \pi i} \int_{\alpha}^{z} \int_{\gamma} \left( \frac{\phi(\zeta) - b}{\phi(t) - b} \right)^{n-1} \frac{f'(t) \phi ' (\zeta)}{\phi(t) - \phi(\zeta)} \, dt \, d\zeta,
\]

and because the last integral tends to zero ($0$) when $n \to \infty$

\begin{equation}  \label{2.3.78}
f(z) = f(\alpha) + \sum_{m=1}^{n-1} \frac{(\phi(z) - b)^m}{m!} \frac{d^{m-1}}{d\alpha^{m-1}} \left( f'(\alpha) \cdot (\Psi(\alpha))^m \right)\ 
\end{equation}

In more general form, if $\Psi(z) = \frac{z - \alpha}{\phi(z) - b} = h(z)$, \\
\(
f(z) = f(\alpha) + \sum_{m=1}^{n-1} \frac{s^m}{m!}. \frac{d^{m-1}}{d\alpha^{m-1}} \left( f'(\alpha) \cdot (h(\alpha))^m \right)\), 
\\[8pt]
and $s = \phi(z)$, specifically,
\begin{equation} \label{2.3.79} 
z = \alpha + \sum_{m=1}^{\infty} \frac{s^m}{m!}  \frac{d^{m-1}}{d\alpha^{m-1}} \left( (h(\alpha))^m \right)
\end{equation} 
So, if the series \eqref{2.3.79} can be calculated, then we can find the root of the equation $z = \alpha + s \cdot h(z)$. The series \eqref{2.3.78} are called Burmann's series. Burmann developed these in 1799 by generalizing Lagrange's series from 1770.
\end{proof}

{\bf Examples.} We present some examples of this section: \\[2pt]
1) Calculate the root of the equation $z = \alpha + \frac{s}{z}$. We take $h(z) = \frac{1}{z}$, therefore

\begin{equation}
z = \alpha + \sum_{m=1}^{\infty} \frac{s^m}{m!} \frac{d^{m-1}}{d\alpha^{m-1}} \left( f'(\alpha) \cdot (h(\alpha))^m \right) = \sum_{m=1}^{\infty} \frac{s^m}{m!} \frac{d^{m-1}}{d\alpha^{m-1}} \left( (h(\alpha))^m \right) \\
\end{equation} 

\begin{equation} \label{2.3.80} 
  = \alpha + \sum_{m=1}^{\infty} (-1)^{m-1} \frac{(2m-2)! \cdot s^m}{(m-1)! \cdot \alpha^{2m-1}}
\end{equation} 

This series converges for $|s| < \frac{|\alpha|^2}{2}$. The second degree equation has two roots: \\[6pt]
$z_1 = \frac{\alpha}{2} \left( 1 + \sqrt{1 + \frac{4s^2}{\alpha^2}} \right)$ and $z_2 = \frac{\alpha}{2} \left( 1 - \sqrt{1 + \frac{4s^2}{\alpha^2}} \right)$.
\\[6pt]
2) Estimation of the roots of a hyperbolic equation. If it is one of the roots of the equation $x = 1 + y \cdot x^{\alpha}$ which tends to $1$ when $y \to 0$, show that:
\[ \log(x) = y + \left( \frac{2\alpha - 1}{2} \right) y^2 + \left( \frac{(3\alpha - 1).(3\alpha - 2)}{3!} \right) y^3 + \dots \]
Expansion makes sense when $|y| < \left| (\alpha - 1)^{\alpha - 1} \cdot \alpha^{-\alpha} \right|$, in the form of McClintock.

3) Estimate the roots of the polynomial in the general form. Consider, for example, the polynomial
\[ f(z) = \alpha_0 + \alpha_1(z - c) + \alpha_2(z - c)^2 + \ldots + \alpha_k(z - c)^k. \]
The equation \( f(z) = 0 \) can be expressed in the form \( z = \alpha + s \cdot h(z) \) where \( s = -\frac{1}{\alpha_1} \) and
\[ h(z) = \alpha_0 + \alpha_1(z - c) + \alpha_2(z - c)^2 + \ldots + \alpha_k(z - c)^k. \]
The Lagrange series for the equation is of the form
\[ z = c + \sum_{m=1}^{\infty} \frac{s^m}{m!} \frac{d^{m-1}}{d\alpha^{m-1}} \left( (h(\alpha))^m \right) \text{ at } z = c. \]
In our case, we have
\[ h^{(m)}(z) = \sum_{n_1+n_2+\ldots+n_k}^{\infty} \frac{\alpha_0^{n_1} \alpha_1^{n_2} \ldots \alpha_k^{n_k}}{n_1! n_2! \ldots n_k!} n!(z - c)^{n-1}. \]
And we take the sum given by $n_1 +n_2 +...+n_k= n$
\[ \frac{d^{n-1}}{dz^{n-1}} \left( (h(z))^n \right) \bigg|_{z=c} = \sum_{n_1+n_2+\ldots+n_k}^{\infty} \frac{(n-1)!}{n_1! n_2! \ldots n_k!} \alpha_0^{n_1} \alpha_2^{n_2} \ldots \alpha_k^{n_k}, \]
where \( n \) is the polynomial degree and \( 2n_2 + \ldots + k n_k = n-1 \) and \( n-1 = 2n_2 + \ldots + k n_k \), \( n_1 = n_2 + 2n_3 + \ldots + (k-1)n_k + 1 \).
After \( s = -\frac{1}{\alpha_1} \), we succeed:
\[ z = c - \frac{\alpha_0}{\alpha_1} \sum_{n_1+n_2+\ldots+n_k=n}^{\infty} \frac{(n-1)!}{n_1! n_2! \ldots n_k!} \left( \frac{\alpha_0 \cdot \alpha_2}{(-\alpha_1)^2} \right)^{n_2} \ldots \left( \frac{\alpha_0^{k-1} \cdot \alpha_k}{-\alpha_1} \right)^{n_k}, \]

where \( z \) is a root of the equation \( f(z) = 0 \) of degree \( n \).
\\[6pt]
4) Solving polynomial equations with the Bring-Jerrard method. Without losing generality, we can find at least one root of the
\begin{equation} \label{2.3.81}  
z^s + a z - q = 0 \quad (s = 2, 3, 4, \ldots)
\end{equation}
by setting \( z^s = \zeta \) $\iff$  \( z_{\text{in}} = \zeta^{1/s} \exp\left( \frac{2k\pi i}{s} \right) = f(\zeta) \). We easily find that \eqref{2.3.81} becomes
\begin{equation} \label{2.3.82}
_k z_{\text{root}} = \zeta^{1/s} \exp\left( \frac{2k\pi i}{s} \right) - \Phi(\zeta)
\end{equation}
where
\begin{equation} \label{2.3.83}
\Phi(\zeta) = \alpha \zeta
\end{equation}
and
\begin{equation} \label{2.3.84}
f(\zeta) = \zeta^{1/s} \exp\left( \frac{2k\pi i}{s} \right).
\end{equation}
The Lagrange Theorem states that for any analytic function in a nearby region of equation \eqref{2.3.82}, then
\[
_k z_{\text{root}} = f(\zeta) + \sum_{n=1}^{\infty} \frac{(-1)^n}{n!} \frac{d^{n-1}}{d\alpha^{n-1}} \left( f'(\zeta) \left( \Phi(\zeta) \right)^n \right) \bigg|_{\zeta \rightarrow q}, \quad k = 0, \ldots, [s-1].
\]
With a known integral we get
\[
\frac{d^{n-1}}{d\zeta^{n-1}} \left( f'(\zeta) \cdot (\Phi(\zeta))^n \right) = \frac{\alpha^n \exp\left( \frac{2k\pi i n}{s} \right) \zeta^{\frac{1+n-ns}{s}} \Gamma\left( \frac{1+n}{s} \right)}{s \Gamma\left( \frac{1+n+s-ns}{s} \right)}
\]
and we come to a root of the form
 
\[
\frac{d^{n-1}}{d\zeta^{n-1}} \left( f'(\zeta) \cdot (\Phi(\zeta))^n \right) = \frac{\alpha^n \exp\left( \frac{2k\pi i n}{s} \right) \zeta^{\frac{1+n-ns}{s}} \Gamma\left( \frac{1+n}{s} \right)}{s \Gamma\left( \frac{1+n+s-ns}{s} \right)}
\]
and we come to a root of the form
\[
_k z_{\text{root}} = \zeta^{1/s} \exp\left( \frac{2k\pi i}{s} \right) + \sum_{n=1}^{\infty} \frac{(-\alpha)^n}{n!} \exp\left( \frac{2k\pi i n}{s} \right) \frac{\zeta^{\frac{1+n-ns}{s}} \Gamma\left( \frac{1+n}{s} \right)}{s \Gamma\left( \frac{1+n+s-ns}{s} \right)}
\]
with \(\zeta \rightarrow q \).

Example for \( s = 5 \) and \( k = 0, \ldots, [s-1] \) we have:
\begin{align*}
_k z_{\text{root}} &= \exp\left( \frac{2k\pi i}{5} \right).{z^{1/5}} \\
&  + \frac{1}{5} \exp\left( \frac{2k\pi i}{5} \right) z^{1/5} \Bigg( -5 +5\, \text{HypergeometricPFQ}\left[ \left\{-\frac{1}{20}, \frac{1}{5}, \frac{9}{20}, \frac{7}{10}\right\}, \right.\\
& \left. \left\{\frac{2}{5}, \frac{3}{5}, \frac{4}{5}\right\}, -256 \exp(2ik\pi) \frac{\alpha^5}{3125z^4} \right] \\
& - \exp\left( \frac{2k\pi i}{5} \right) \alpha.z^{-4/5} \, \text{HypergeometricPFQ}\left[ \left\{\frac{3}{20}, \frac{2}{5}, \frac{13}{20}, \frac{9}{10}\right\}, \right.\\
& \left. \left\{\frac{3}{5}, \frac{4}{5}, \frac{6}{5}\right\}, -\frac{256 \exp(2ik\pi) \alpha^5}{3125z^4} \right] \\
& - \exp\left( \frac{4k\pi i}{5} \right) \frac{1}{5z^{8/5}} \alpha^2 \, \text{HypergeometricPFQ}\left[ \left\{\frac{7}{20}, \frac{3}{5}, \frac{17}{20}, \frac{11}{10}\right\}, \right.\\
& \left. \left\{\frac{4}{5}, \frac{6}{5}, \frac{7}{5}\right\}, -\frac{256 \exp(2ik\pi) \alpha^5}{3125z^4} \right] \\
& - \exp\left( \frac{6k\pi i}{5} \right) \frac{1}{25z^{12/5}} \alpha^3 \, \text{HypergeometricPFQ}\left[ \left\{\frac{11}{20}, \frac{4}{5}, \frac{21}{20}, \frac{13}{10}\right\}, \right.\\
& \left. \left\{\frac{6}{5}, \frac{7}{5}, \frac{8}{5}\right\}, -\frac{256 \exp(2ik\pi) \alpha^5}{3125z^4} \right] \Bigg)
\end{align*}

Where \( z \to q \) and \(\text{HypergeometricPFQ}\) is the generalized hypergeometric function \cite{wolfram} ${}_pF_q \left( \{a_1, \ldots, a_p\}; \{b_1, \ldots, b_q\}; z \right)$. With another way using the Gauss Theorem, an infinite series is decomposed into infinite series of hypergeometric functions.

\[
\Psi(j) = \left( \frac{\omega. q}{s-1} \right)^{\frac{j.s}{s-1}} 
\prod_{k=0}^{s-1} \Gamma\left( \frac{s.j}{s.(s-1)} + 1/s + k/s \right) \frac{1}{\Gamma\left( \frac{j}{N-1} + 1 \right) \cdot \prod_{k=0}^{s-2} \Gamma\left( \frac{j + k + 2}{s-1} \right)}
\]

And finally we get:

\[
\begin{aligned}
x_{\text{root}} = & \omega^{-1} - \frac{q}{(s-1)^2} \sqrt{\frac{s}{2 \pi (s-1)}} \sum_{i=0}^{s-2} \Psi(j)_s \cdot F^{s+1}. \\
& [\left( \frac{j.s}{s(s-1)} + \frac{1}{s} \right), \left( \frac{j.s}{s(s-1)} + \frac{2}{s} \right), \ldots, \left( \frac{j.s}{s(s-1)} + \frac{s}{s} \right); \\
& \left( \frac{j+2}{s-1} \right), \left( \frac{j+3}{s-1} \right), \ldots, \left( \frac{j+s}{s-1} \right); \\
& \left( \left( \frac{q \cdot \omega}{s-1} \right)^{s-1} \cdot s^s \right)]
\end{aligned}
\]

Where $\omega = \exp\left(\frac{2\pi i}{s-1}\right)$. A root of an equation can be expressed as a sum of at most $s - 1$ hypergeometric functions. Applying the Bring-Jerrard method to the quintic equation, we define the following functions:
\\[6pt]
 $F_1(q) = F_2(q)$, 

$F_2(q) = {}_4F_3\left(\frac{1}{5}, \frac{2}{5}, \frac{3}{5}, \frac{4}{5}; \frac{1}{2}, \frac{3}{4}, \frac{5}{4}; \frac{3125 q^4}{256}\right)$, 

$F_3(q) = {}_4F_3\left(\frac{9}{20}, \frac{13}{20}, \frac{17}{20}, \frac{21}{20}; \frac{3}{4}, \frac{5}{4}, \frac{3}{2}; \frac{3125 q^4}{256}\right)$, 

$F_4(q) = {}_4F_3\left(\frac{7}{10}, \frac{9}{10}, \frac{11}{10}, \frac{13}{10}; \frac{5}{4}, \frac{3}{2}, \frac{7}{4}; \frac{3125 q^4}{256}\right)$. 
\\[6pt]
These are hypergeometric functions listed above. The roots of the quintic equation are:
\begin{align*}
x_1 &= -q^4 \cdot F_1(q), \\
x_2 &= -F_1(q) + \frac{q}{4} \cdot F_2(q) + \frac{5}{32} \cdot (q^2) \cdot F_3(q) + \frac{5}{32} \cdot (q^3) \cdot F_3(q), \\
x_3 &= -i \cdot F_1(q) + \frac{q}{4} \cdot F_2(q) - i \cdot \frac{5}{32} \cdot (q^2) \cdot F_3(q) + \frac{5}{32} \cdot (q^3) \cdot F_3(q), \\
x_4 &= -i \cdot F_1(q) + \frac{q}{4} \cdot F_2(q) - i \cdot \frac{5}{32} \cdot (q^2) \cdot F_3(q) - \frac{5}{32} \cdot (q^3) \cdot F_3(q), \\
x_5 &= -i \cdot F_1(q) + \frac{q}{4} \cdot F_2(q) + i \cdot \frac{5}{32} \cdot (q^2) \cdot F_3(q) - \frac{5}{32} \cdot (q^3) \cdot F_3(q).
\end{align*}
This is the same result that was achieved with the method of differential solvers developed by James Cockle and Robert Harley in 1860.

\vspace{10pt}
\noindent \textbf{3.9 Generalized existence and root number theorem of a random transcendental equation.}

\vspace{6pt}
\noindent \textbf{3.9.1. Definition} Type of function means one of the 5 general forms of functions, which appear in mathematics namely: Exponential Function, Logarithmic Function, Trigonometric Function, Dynamic Function, Exponential Dynamic Function.

\vspace{6pt}
\noindent \textbf{3.9.2. Definitions of primary functions}
\\[6pt]
\noindent (i). A primary simple transcendental equation is any equation of the form: $\Psi(z) = q(z) + t = 0$, with $t \in \mathbb{C}$, which has its roots in $\mathbb{C}$.

\noindent (ii). A primary composite transcendental equation is any equation of the form: $\Psi_2(z) = q(z) + m \cdot p(z) + t = 0$, which has its roots in $\mathbb{C}$. Furthermore, the coefficients $m$, $t$ are different from $0$, and also take values from $\mathbb{C}$. In general, the functions $q(z)$, $p(z)$ are different types, with values in $\mathbb{C}$.
\\[6pt]
\textbf{Theorem 9. }\label{theo:random}
%\begin{theorem}\label{theo:random}
  \textit{Every random transcendental equation of the form:}
\[
p_{\lambda}(z) = \frac{1}{m_{\lambda}} \sum_{\substack{k=1 \\ k \ne \lambda}}^{n} m_k \cdot p_k(z) + \frac{t}{m_{\lambda}} = 0, \quad k=1,2,\ldots,n, \quad \lambda \ne k 
\]
\begin{equation}
\label{eq:transcendental_eq}
\end{equation}
\textit{has a plurality of roots, the union of the fields of roots \(L_1, L_2, \ldots, L_n\) corresponding to each transformation \(p_{\lambda}(z)\) and divided by \(L = \bigcup_{i=1}^{n} L_i\), which fields are derived from the individual functions \(p_{\lambda}(z)\) and Lagrange's Theorem, taking a cyclically transposed position within the transcendental equation. This can be expressed as taking \(n\) forms with the following equations:}

\[
m_1 \cdot p_1(z) + \sum_{k=2}^{n} m_k \cdot p_k(z) + t = 0
\]
\begin{equation}
\label{eq:form1}
\end{equation}

\[
m_2 \cdot p_2(z) + \sum_{\substack{k=1 \\ k \ne 2}}^{n} m_k \cdot p_k(z) + t = 0
\]
\begin{equation}
\label{eq:form2}
\end{equation}

\[
\vdots
\]

\[
m_n \cdot p_n(z) + \sum_{k=1}^{n-1} m_k \cdot p_k(z) + t = 0
\]
\begin{equation}
\label{eq:formn}
\end{equation}

\textit{These are derived as shown by the generalized transcendental equation:}

\[
\sum_{k=1}^{n} m_k \cdot p_k(z) + t = 0 \quad \text{with} \quad k=1, 2, \ldots, n
\]

\textit{if and only if:}

\textit{i) The functions are analytic functions \(p_i(z)\) at all points within and on contour \(C\). At the same time they are different as kind functions or of different form or power, in general.}

\textit{ii) The coefficients \(m_i, t \ne 0, \, i \in \mathbb{N}_+\) take values in \(\mathbb{C}\), is at most \(n\) in number, with at least 2 coefficients \(m_i\) must be non-zero.}

\textit{iii) The fields of roots \(L_1, L_2, \ldots, L_n\) of \(\Psi_1, \Psi_2, \ldots, \Psi_n\) forms of equal equations:}
\[
\Psi_{\lambda}(z) = m_{\lambda} \cdot p_{\lambda}(z) + \sum_{\substack{k=1 \\ k \ne \lambda}}^{n} m_k \cdot p_k(z) + t = 0, \quad 1 \le \lambda \le n
\]
\textit{are located according to the Theorem of Lagrange and belong to \(\mathbb{C}\).}

\textit{iv) The number of sub-fields of the roots of the equation is \(n\), and hence the total field of the roots of the equation:}
\[
\sum_{i=1}^{n} m_i \cdot p_i(z) + t = 0, \quad \lambda = 1,2,\ldots,n
\]
is \(L = \bigcup_{i=1}^{n} L_i\).  
 
%\end{theorem}

\begin{proof}
The pairs of ordered pairs of the function

\begin{equation} \label{2.3.90}
 \Psi_{\lambda}(z) = m_{\lambda} \cdot p_{\lambda}(z) + \sum_{\substack{k=1 \\ k \ne \lambda}}^{n} m_k \cdot p_k(z) + t = 0 
\end{equation}

come in the form of:
\begin{equation} \label{2.3.91}
m_{\lambda} \cdot p_{\lambda}(z) = -\sum_{\substack{k=1 \\ k \ne \lambda}}^{n} m_k \cdot p_k(z) - t = 0 
\end{equation}

and then in the final form:
\begin{equation} \label{2.3.92}
p_{\lambda}(z)=\left(-\frac{1}{m_{\lambda}}\right)\sum_{\substack{k=1 \\ k \ne \lambda}}^{n} m_k \cdot p_k(z) - \frac{t}{m_{\lambda}} = 0 
\end{equation} 
This is the decisive form for finding the individual fields of its roots. It proves that \(n\) is the number of iterative transfers.
\\[6pt]
\textit{\textbf{Local proof:}}
Suppose we have a transfer \(\kappa\) with the property "There are \(m_1, m_2 \ldots,\) distinguished elements of the set \(\{1, 2, \ldots, n\}\) with the matching \(\alpha_1 \pi = \alpha_2, \alpha_2 \pi = \alpha_3, \ldots, \alpha_1\) then \(\pi\) is called a circle of length \(k\) and is written:
\[
(m_1. m_2. \ldots. m_k) = (m_2. m_3. \ldots. m_k. m_1) = (m_3. m_4. \ldots. m_k. m_1. m_2) = \ldots = (m_k. m_1. \ldots. m_{k-1})
\]
So it will have a number of transpositions \(k\). Similarly, if \(k=n\), then the number of transpositions is \(n\).

According to the developed form \eqref{2.3.92}, we obtain the analysis in the following forms:
\[
p_1(z) + \sum_{k=2}^{n} \left(\frac{m_k}{m_1}\right)p_k(z) + \frac{t}{m_1} = 0
\]
\[
p_2(z) = \sum_{\substack{k=1 \\ k \ne 2}}^{n} \left(\frac{m_k}{m_2}\right)p_k(z) + \frac{t}{m_2} = 0
\]
\[
\vdots
\]
\[
p_{\lambda}(z) + \sum_{\substack{k=1 \\ k \ne \lambda}}^{n} \left(\frac{m_k}{m_{\lambda}}\right)p_k(z) + \frac{t}{m_{\lambda}} = 0
\]
\[
p_n(z) + \sum_{k=1}^{n-1} \left(\frac{m_k}{m_n}\right)p_k(z) + \frac{t}{m_n} = 0
\]

Now we take \( p_{\lambda}(z) \), \( f_{\lambda}(z) \), and \( \phi_{\lambda}(z) \) to be analytic functions of \( z \) on and inside the contour \( C \) (a closed area \( C \)), surrounding a point \( -\frac{t}{m_{\lambda}} \). We take \( -\frac{1}{m_{\lambda}} \) such that the inequality
\[
\left| -\frac{1}{m_{\lambda}} \sum_{\substack{k=1 \\ k \ne \lambda}}^{n} m_k \cdot p_k(z) \right| < \left| z - \left(-\frac{t}{m_{\lambda}}\right) \right|
\]
is satisfied at all points of complex numbers on the perimeter of the contour \( C \) with \( m_1 = 1 \).

Taking therefore \(p_\lambda(z)\) = \( {}_\lambda\zeta_{t_{\lambda}}, \, n \geq \lambda \geq 1, \, t \geq 1, \, \{n, t, \lambda \in \mathbb{N}\} \) and by doing the inverse of the function, we get:
\[
f_({}_\lambda\zeta_{t_{\lambda}}) = p_{\lambda}^{-1}( {}_\lambda\zeta_{t_{\lambda}})
\]
and from the generalized transcendental equation we will have:
\[
p_{\lambda}(z) = \left( -\frac{1}{m_{\lambda}} \right) \sum_{\substack{k=1 \\ k \ne \lambda}}^{n} m_k \cdot p_k(z) - \frac{t}{m_{\lambda}}, \quad m_1 = 1
\]

Then by simple correlation, it applies:
\[
{}_\lambda\zeta_{t_{\lambda}} = p_{\lambda}(z) = -\left( \frac{1}{m_{\lambda}} \right) \sum_{\substack{k=1 \\ k \ne \lambda}}^{n} m_k \cdot p_k(z) - \frac{t}{m_{\lambda}}, \quad i \geq 1, \, \{i, \lambda, k \in \mathbb{N}\}
\]

and 
\[
\phi_{\lambda} {}(_\lambda\zeta_{i_{\lambda}}) = \sum_{\substack{i=1 \\ i \ne \lambda}}^{n} m_i \cdot p_i \left( p_{\lambda}^{-1}(\zeta) \right), \quad i \geq 1, \, i \ne \lambda, \, i \leq n
\]

Then after the replacement:
\[
f( {}_\lambda\zeta_{t_{\lambda}}) = p_k^{-1} ( {}_\lambda\zeta_{t_{\lambda}})
\]
we have the equation:
\[
{}_\lambda\zeta_{t_{\lambda}} = \left( -\frac{1}{m_{\lambda}} \right) \sum_{\substack{i=1 \\ i \ne \lambda}}^{n} m_i \cdot p_i ({}_\lambda\zeta_{i_{\lambda}}) - \frac{t}{m_{\lambda}}
\]
which is considered as an equation of \( \zeta \), that has a root (or more than one if \( q() \) is an Exponential, Logarithmic, Trigonometric, Dynamic Function (\( n > 1 \)), Exponential Dynamic Function) within the contour \( C \). Any function of \( \zeta \) that is analytic within and on \( C \) can be developed into a power series with a similar \( C \), can be developed into a power series with an identical variable \( {}_\lambda\zeta_{t_{\lambda}} \), i.e., \( {}_\lambda\zeta_{t_{0}} = w \to \left\{ -\frac{t}{m_{\lambda}} \right\} \) which is an initial value, and then we can find all the partial solutions per field by form \( {}_\lambda\zeta_{t_{\lambda}}, \, n \geq i, \, \lambda \geq 1, \, \{ \lambda, i \in \mathbb{N} \} \) and according to the general relation.

\begin{equation}  \label{2.3.92.1}
{}_\lambda\zeta_{i_{\lambda}} = f_{\lambda}(w) \bigg|_{w \to -\frac{t}{m_{\lambda}}}  + \sum_{s=1}^{\infty} \frac{(-1/m_{\lambda})^s}{s!} \left. \frac{d^{s-1}}{d w^{s-1}} \left( f'_{\lambda}(w) \cdot \phi^s_{\lambda}(w) \right) \right|_{w \to -\frac{t}{m_{\lambda}}}, \quad 1 \leq \lambda \leq n.
\end{equation}
Remark: we can also have a more clearly ideal solution if $\left| m_r \right| < \left| m_i \right|$ for every $i \geq 1$, $i \ne \lambda$, $i \leq n$ or $m_r = \min\{\left|m_1\right|, \left|m_2\right|, \ldots\}$. Thus, $\left|\frac{m_r}{m_i}\right| < 1$ because we will have a given convergence in the sum and then we finally get for the root the relation:
\begin{equation}  \label{2.3.93}
{}_\lambda\zeta_{i_{\lambda}} = f_{\lambda}(w) \bigg|_{w \to -\frac{t}{m_{\lambda}}} + \sum_{s=1}^{\infty} \frac{(-m_r/m_{\lambda})^s}{s!} \left. \frac{d^{s-1}}{d w^{s-1}} \left( f'_{\lambda}(w) \cdot \phi^s_{\lambda}(w) \right) \right|_{w \to -\frac{t}{m_{\lambda}}}, \quad \lambda \in \mathbb{N}, \, 1 \leq \lambda \leq n, \, \lambda \ne r.
\end{equation}
Relations \eqref{2.3.92}, \eqref{2.3.93} are the relations for finding each root of the primary composite transcendental equation:
\[
\sum_{i=1}^{n} \left( m_i \cdot p_i(z) \right) + t = 0, \quad i \leq n, \, i \in \mathbb{N},
\]
with its total number of roots being \( K = \sum_{i=1}^{n} k_i \) with \( k_i \) the partial integer of the sets of roots, total number \( n \), and thus \( K \) equal to the number of roots of the transcendental equation. In particular, with the fields of the roots of the equation in the form:
\[
m_{\lambda} \cdot p_{\lambda}(z) + \sum_{\substack{k=1 \\ k \ne \lambda}}^{n} \left( m_k \cdot p_k(z) \right) + t = 0,
\]
which we call \( L_k \), concerns only this form, i.e.
\[
p_{\lambda}(z) + \sum_{\substack{k=1 \\ k \ne \lambda}}^{n} \left(\frac{m_k}{m_{\lambda}}\right)p_k(z) + \frac{t}{m_{\lambda}} = 0.
\]
Now, for the generalization of the cases, because \( i, \lambda \), they take values from \( 1 \) to \( n \), the number of partial fields of roots will be \( L \), and as a consequence we will have that the total field of roots of the transcendental equation:
\[
\sum_{k=1}^{n} \left( m_k \cdot p_k(z) \right) + t = 0, \quad m_1 = 1,
\]
will be:
\[
L = \bigcup_{i=1}^{n} L_i.
\]

\end{proof}

\textbf{3.10 A General form of 2 arguments of 7th degree is solvable only by sub-conditions using hyper-geometric series or with methods of approximations} \label{des:general1}
\\[4pt]
By the definition "to prove whether there is a solution for all 7th degree equations using algebraic (variant: continuous) functions of two arguments" means that if we talk about the 7th degree equation of the form \(x^7 + c x^3 + b x^2 + a x + 1\) from  the Kolmogorov and after Arnold we will have the following forms as an example \(x^7 + x^3 + b x^2 + a x + 1\) or \(a x^7 + b x^3 + x^2 + x + 1\) but not in the trinomial form such as the form \(x^7 + a x + b = 0\) which is solved by hypergeometric series or by the method of approximation. We proved that a generalized form of degree 7 has the form \(x^7 + a x^3 + b x^2 + a x + 1\) which is not solvable with roots, so we look for various functions such as hypergeometric series to solve this form and more specifically with 2 arguments as form \(x^7 + x^3 + b x^2 + a x + 1\) if it can be solved with the hypergeometric series as mentioned before. We can consider 2 of the 4 parameters as constants and calculate the other 2 according to our computational requirements to get a good approximation, as we will see below. It is of course necessary to calculate the roots of the trinomial in general form.
\\[4pt]
\textbf{3.10.1 Solving the polynomial trinomial form \(x^5 - \alpha x^b - q = 0\), \(1 \leq b \leq 4\).}
\\[2pt]
The only process with which we can compute the roots of a trinomial in both direct approximation and in type form are using hypergeometric series, as originally used by the Bring-Jerrard Vrin form for the 5th degree. Here we use the complete Burman Lagrangian series in a different new form and will preserve it for each analysis with the existing known hypergeometric series and in higher degree polynomial equations. We are looking to find a general process for the trinomial of any exponent. Of course, we were talking about exponents inside the \(\mathbb{R}\). Starting from on page 32 we will extend the trinomial to be in 2 exponents.
\begin{equation} \tag{3.2.104} \label{2.3.94}
z^s - a z^b - q = 0 \quad \{s, b\} \in \mathbb{N}
\end{equation}
By setting \(z^s = \zeta\) if and only if \(z_i = \zeta^{1/s} \exp\left(\frac{2k\pi i}{s}\right) = f(\zeta)\), we easily find that \eqref{2.3.94} becomes:

\begin{equation} \tag{3.2.105} \label{2.3.95}
z_{\text{root}, k} = \zeta^{1/s} \exp\left(\frac{2k\pi i}{s}\right) - \alpha \Phi(\zeta)
\end{equation}
with
\begin{equation} \tag{3.2.106} \label{2.3.96}
\Phi(\zeta) = \zeta^b
\end{equation}
\begin{equation} \tag{3.2.107}  \label{2.3.97}
f(\zeta) = \zeta^{1/s} \exp\left(\frac{2k\pi i}{s}\right)
\end{equation}
Lagrange's theorem states that for any analytic function in a nearby region of equation \eqref{2.3.95} it will hold:
\begin{equation} \tag{3.2.108} \label{2.3.98}
z_{\text{root}, k} = f(\zeta) + \sum_{n=1}^{\infty} \frac{(-1)^n}{n!} \left. \frac{d^{n-1}}{d \alpha^{n-1}} \left( f'(\zeta) \cdot \Phi^n(\zeta) \right) \right|_{\zeta \to q}, \quad k = 0, \ldots, [s-1]
\end{equation}
According to the known integral we get:
\[
\frac{d^{n-1}}{d \zeta^{n-1}} \left( f'(\zeta) \cdot \Phi^n(\zeta) \right) = \frac{\alpha^n \exp\left(\frac{2k\pi.i.b.n}{s}\right) \zeta^{\left(\frac{1+b. n - n.s}{s}\right)} \Gamma\left(\frac{1+b n}{s}\right)}{s \Gamma\left(\frac{1+b.n + s - n.s}{s}\right)}
\]
and we come to a root of the form:
\begin{equation} \tag{3.2.109} \label{2.3.99}
z_{\text{root}, k} = \zeta^{1/s} \exp\left(\frac{2k\pi i}{s}\right) + \sum_{n=1}^{\infty} \frac{(-\alpha)^n}{n!} \frac{\exp\left(\frac{2k\pi i b n}{s}\right) \zeta^{\left(\frac{1 + bn - ns}{s}\right)} \Gamma\left(\frac{1 + bn}{s}\right)}{s \Gamma\left(\frac{1 + bn + s - ns}{s}\right)}
\end{equation}

i). Example for \(s = 5\) and \(b = 2\) and \(k = 0, \ldots, [s-1]\) and \(z \to q\) we have the roots

\begin{equation} \tag{3.2.110} \label{2.3.100}
\begin{aligned}
z_{\text{root}, k} &= \exp\left(\frac{2 i k \pi}{5}\right) q^{1/5} \\
&\quad + \frac{1}{625 q^{11/5}} \exp\left(\frac{2 i k \pi}{5}\right) \left(-625 q^{12/5} + 625 q^{12/5} \cdot {}_4F_3 \left[ \begin{array}{c}
- \frac{1}{15}, \frac{1}{10}, \frac{4}{15}, \frac{3}{5} \\
\frac{1}{5}, \frac{2}{5}, \frac{4}{5} \end{array} \middle| \frac{108 \exp(4 i k \pi) \alpha^5}{3125 q^3} \right] \right) \\
&\quad - 125 \exp\left(\frac{4 i k \pi}{5}\right) \alpha q^{9/5} \cdot {}_4F_3 \left[ \begin{array}{c}
\frac{2}{15}, \frac{3}{10}, \frac{7}{15}, \frac{4}{5} \\
\frac{2}{5}, \frac{3}{5}, \frac{6}{5} \end{array} \middle| \frac{108 \exp(4 i k \pi) \alpha^5}{3125 q^3} \right] \\
&\quad + 5 \exp\left(\frac{12 i k \pi}{5}\right) \alpha^3 q^{3/5} \cdot {}_4F_3 \left[ \begin{array}{c}
\frac{8}{15}, \frac{7}{10}, \frac{13}{15}, \frac{6}{5} \\
\frac{4}{5}, \frac{7}{5}, \frac{8}{5} \end{array} \middle| \frac{108 \exp(4 i k \pi) \alpha^5}{3125 q^3} \right] \\
&\quad + \exp\left(\frac{16 i k \pi}{5}\right) \alpha^4 \cdot {}_4F_3 \left[ \begin{array}{c}
\frac{11}{15}, \frac{9}{10}, \frac{16}{15}, \frac{7}{5} \\
\frac{6}{5}, \frac{8}{5}, \frac{9}{5} \end{array} \middle| \frac{108 \exp(4 i k \pi) \alpha^5}{3125 q^3} \right]
\end{aligned}
\end{equation}

ii) Example for \( s=5 \) and \( b=3 \) and \( k=0,\ldots,[s-1] \) and \( z \rightarrow q \) we have the roots:
\begin{equation} \tag{3.2.111} \label{2.3.101}
\begin{aligned}
z_{\text{root}, k} = &\exp\left(\frac{2 i k \pi}{5}\right) q^{1/5} \\
&- \frac{1}{625 z^{7/5}} \exp\left(\frac{2 i k \pi}{5}\right) \left(625 q^{8/5} - 625 q^{8/5} \cdot {}_4F_3 \left[ \begin{array}{c}
- \frac{1}{10}, \frac{1}{15}, \frac{2}{5}, \frac{11}{15} \\
\frac{1}{5}, \frac{3}{5}, \frac{4}{5} \end{array} \middle| \frac{-108 \exp(6 i k \pi) \alpha^5}{3125 q^2} \right] \right) \\
&+ 125 \exp\left(\frac{6 i k \pi}{5}\right) \alpha q^{6/5} \cdot {}_4F_3 \left[ \begin{array}{c}
\frac{1}{10}, \frac{4}{15}, \frac{3}{5}, \frac{14}{15} \\
\frac{2}{5}, \frac{4}{5}, \frac{6}{5}, \end{array} \middle| \frac{-108 \exp(6 i k \pi) \alpha^5}{3125 q^2} \right] \\
&- 25 \exp\left(\frac{12 i k \pi}{5}\right) \alpha^2 q^{4/5} \cdot {}_4F_3 \left[ \begin{array}{c}
\frac{3}{10}, \frac{7}{15}, \frac{4}{5}, \frac{17}{15} \\
\frac{3}{5}, \frac{6}{5}, \frac{7}{5} \end{array} \middle| \frac{-108 \exp(6 i k \pi) \alpha^5}{3125 q^2} \right] \\
&+ 2 \exp\left(\frac{24 i k \pi}{5}\right) \alpha^4 \cdot {}_4F_3 \left[ \begin{array}{c}
\frac{7}{10}, \frac{13}{15}, \frac{6}{5}, \frac{23}{15} \\
\frac{7}{5}, \frac{8}{5}, \frac{9}{5} \end{array} \middle| \frac{-108 \exp(6 i k \pi) \alpha^5}{3125 q^2} \right]
\end{aligned}
\end{equation}

From relation \eqref{2.3.99}, we can derive any relation that gives its roots of 5th degree in trinomial form for \( 1 \leq b \leq 4 \).

\textbf{3.10.2 Solving the polynomial trinomial form \( x^6 - \alpha x^b - q = 0 \), \( 1 \leq b \leq 5 \)}. \label{desc:new} From the general relation \eqref{2.3.99}, 
\[
z_{\text{root}, k} = \zeta^{1/s} \exp\left(\frac{2k\pi i}{s}\right) + \sum_{n=1}^{\infty} \frac{(-\alpha)^n}{n!} \frac{\exp\left(\frac{2k\pi i b n}{s}\right) \zeta^{\left(\frac{1 + bn - ns}{s}\right)} \Gamma\left(\frac{1 + bn}{s}\right)}{s \Gamma\left(\frac{1 + bn + s - ns}{s}\right)}
\] 
we can easily calculate the hypergeometric function that will solve and find the roots of the 6th degree trinomial for each exponent \( 1 \leq b \leq 5 \) if \( s=6 \). 
\newline \indent $^{*}$ Corresponding author\texttt{ nikmatza@gmail.com}. Relations 3.2.110,3.2.111 were developed with the
program mathematica 14.1, and are therefore considered undoubtedly correct. 

i) For example, if \( s=6 \), \( b=1 \), and \( k=0,\ldots,[s-1] \) and \( z \rightarrow q = b \), 
\begin{equation} \tag{3.2.112} \label{2.3.102}
\begin{aligned}
z_{\text{root}, k} = &\exp\left(\frac{i k \pi}{3}\right) q^{1/6} \\
&+ \frac{1}{31104} \exp\left(\frac{i k \pi}{3}\right) q^{1/6} \left(-31104 + 31104. \cdot
{}_5F_3 \left[ \begin{array}{c}
- \frac{1}{30}, \frac{1}{6}, \frac{11}{30}, \frac{17}{30}, \frac{23}{30} \\
\frac{1}{3}, \frac{1}{2}, \frac{2}{3}, \frac{5}{6} \end{array} \middle| \frac{-3125 \alpha^6 \exp(2 i k \pi)}{46656 q^5} \right] \right) \\
&- 5184 q^{-5/6} \alpha \exp\left(\frac{i k \pi}{3}\right) \cdot {}_5F_3 \left[ \begin{array}{c}
\frac{2}{15}, \frac{1}{3}, \frac{8}{15}, \frac{11}{15}, \frac{14}{15} \\
\frac{1}{2}, \frac{2}{3}, \frac{5}{6}, \frac{7}{6} \end{array} \middle| \frac{-3125 \exp(2 i k \pi) \alpha^6}{46656 q^5} \right] \\
&- 1296 q^{-5/3} \alpha^2 \exp\left(\frac{2 i k \pi}{3}\right) \cdot {}_5F_3 \left[ \begin{array}{c}
\frac{3}{10}, \frac{1}{2}, \frac{7}{10}, \frac{9}{10}, \frac{11}{10} \\
\frac{2}{3}, \frac{5}{6}, \frac{7}{6}, \frac{4}{3} \end{array} \middle| \frac{-3125 \exp(2 i k \pi) \alpha^6}{46656 q^5} \right] \\
&- 384 q^{-5/2} \alpha^3 \exp\left(i k \pi\right) \cdot {}_5F_3 \left[ \begin{array}{c}
\frac{7}{15}, \frac{2}{3}, \frac{13}{15}, \frac{16}{15}, \frac{19}{15} \\
\frac{5}{6}, \frac{7}{6}, \frac{4}{3}, \frac{3}{2} \end{array} \middle| \frac{-3125 \exp(2 i k \pi) \alpha^6}{46656 q^5} \right] \\
&- 91 q^{-10/3} \alpha^4 \exp\left(\frac{4 i k \pi}{3}\right) \cdot {}_5F_3 \left[ \begin{array}{c}
\frac{19}{30}, \frac{5}{6}, \frac{31}{30}, \frac{37}{30}, \frac{43}{30} \\
\frac{7}{6}, \frac{4}{3}, \frac{3}{2}, \frac{5}{3} \end{array} \middle| \frac{-3125 \exp(2 i k \pi) \alpha^6}{46656 q^5} \right] \\
\end{aligned}
\end{equation}

ii. In case if $s=6$ and $b=5$ and $k=0,\ldots,[s-1]$ and $z \rightarrow q = b$ we have the root
\\[2pt]
\begin{equation} \tag{3.2.113} \label{2.3.103}
\begin{aligned}
k_{z_{\text{root}}} = & e^{\frac{ik\pi}{3}} q^{1/6} + \frac{1}{93312 q^{2/3} \Gamma\left(\frac{7}{6}\right)} \left( -15552 e^{\frac{ik\pi}{3}} q^{5/6} \Gamma\left(\frac{1}{6}\right) + 15552 ae^{2i a \pi 3} q^{2/3} \Gamma\left(\frac{7}{6}\right) \right) \\
& + 93312 e^{\frac{ik\pi}{3}} q^{5/6} \Gamma\left(\frac{7}{6}\right) \cdot {}_5F_4 \left[ \begin{array}{c} \frac{-1}{6}, \frac{1}{30}, \frac{7}{30}, \frac{13}{30}, \frac{19}{30} \\ \frac{1}{6}, \frac{1}{3}, \frac{1}{2}, \frac{2}{3} \end{array} \middle| \frac{-3125 a^6 e^{10ik\pi}}{46656 q} \right] \\
& + 6480 a^2 e^{\frac{11ik\pi}{3}} \sqrt{q} \Gamma\left(\frac{7}{6}\right) \cdot {}_5F_4 \left[ \begin{array}{c} \frac{1}{6}, \frac{11}{30}, \frac{17}{30}, \frac{23}{30}, \frac{29}{30} \\ \frac{1}{2}, \frac{2}{3}, \frac{5}{6}, \frac{4}{3} \end{array} \middle| \frac{-3125 a^6 e^{10ik\pi}}{46656 q} \right] \\
& + 2880 a^3 e^{\frac{16ik\pi}{3}} q^{1/3} \Gamma\left(\frac{7}{6}\right) \cdot {}_5F_4 \left[ \begin{array}{c} \frac{1}{3}, \frac{8}{15}, \frac{11}{15}, \frac{14}{15}, \frac{17}{15} \\ \frac{2}{3}, \frac{5}{6}, \frac{7}{6}, \frac{3}{2} \end{array} \middle| \frac{-3125 a^6 e^{10ik\pi}}{46656 q} \right] \\
& + 1215 a^4 e^{\frac{7ik\pi}{3}} q^{1/6} \Gamma\left(\frac{7}{6}\right) \cdot {}_5F_4 \left[ \begin{array}{c} \frac{1}{2}, \frac{7}{10}, \frac{9}{10}, \frac{11}{10}, \frac{13}{10} \\ \frac{5}{6}, \frac{7}{6}, \frac{4}{3}, \frac{5}{3} \end{array} \middle| \frac{-3125 a^6 e^{10ik\pi}}{46656 q} \right] \\
& + 448 a^5 e^{\frac{26ik\pi}{3}} \Gamma\left(\frac{7}{6}\right) \cdot {}_5F_4 \left[ \begin{array}{c} \frac{2}{3}, \frac{13}{15}, \frac{16}{15}, \frac{19}{15}, \frac{22}{15} \\ \frac{7}{6}, \frac{4}{3}, \frac{3}{2}, \frac{11}{6} \end{array} \middle| \frac{-3125 a^6 e^{10ik\pi}}{46656 q} \right]
\end{aligned}
\end{equation}
\noindent

\textbf{3.10.3 Solving the polynomial trinomial form \( x^7 - \alpha x^b - q = 0 \), \( 1 \leq b \leq 6 \)}. Similarly, from the general relation, we obtain two typical cases for the seventh degree: \\
i) For the case where \( s=7 \), \( b=1 \), and \( k=0, \ldots, [s-1] \), and \( z \to q=b \). We have the root:
\newline \indent $^{*}$ Corresponding author\texttt{ nikmatza@gmail.com}. Relations 3.2.112,3.2.113, 3.2.114 were developed with the
program mathematica 14.1, and are therefore considered undoubtedly correct.

\vspace{-2\baselineskip} % Adjust vertical space before equation
\begin{equation} \tag{3.2.114} \label{2.3.104}
\begin{aligned}
z_{\text{root}, k} = & \exp\left(\frac{2 i k \pi}{7}\right) q^{1/7} + \frac{1}{16807} \exp\left(\frac{2 i k \pi}{7}\right) q^{1/7} \\
& \quad \times \left(-16807 + 16807 \cdot {}_6F_5 \left[ \begin{array}{c}
- \frac{1}{42}, \frac{1}{7}, \frac{13}{42}, \frac{10}{21}, \frac{9}{14}, \frac{17}{21} \\
\frac{2}{7}, \frac{3}{7}, \frac{4}{7}, \frac{5}{7}, \frac{6}{7} \end{array} \middle| \frac{46656 \alpha^7 \exp(2 i k \pi)}{823543 q^6} \right] \right) \\
& + 2401 q^{-6/7} \alpha \exp\left(\frac{2 i k \pi}{7}\right) \\
& \quad \times {}_6F_5 \left[ \begin{array}{c}
\frac{5}{42}, \frac{2}{7}, \frac{19}{42}, \frac{13}{21}, \frac{11}{14}, \frac{20}{21} \\
\frac{3}{7}, \frac{4}{7}, \frac{5}{7}, \frac{6}{7}, \frac{8}{7} \end{array} \middle| \frac{46656 \alpha^7 \exp(2 i k \pi)}{823543 q^6} \right] \\
& - 686 q^{-12/7} \alpha^2 \exp\left(\frac{4 i k \pi}{7}\right) \\
& \quad \times {}_6F_5 \left[ \begin{array}{c}
\frac{11}{42}, \frac{3}{7}, \frac{25}{42}, \frac{16}{21}, \frac{13}{14}, \frac{23}{21} \\
\frac{4}{7}, \frac{5}{7}, \frac{6}{7}, \frac{8}{7}, \frac{9}{7} \end{array} \middle| \frac{46656 \alpha^7 \exp(2 i k \pi)}{823543 q^6} \right] \\
& + 245 q^{-18/7} \alpha^3 \exp\left(\frac{6 i k \pi}{7}\right) \\
& \quad \times {}_6F_5 \left[ \begin{array}{c}
\frac{17}{42}, \frac{4}{7}, \frac{31}{42}, \frac{19}{21}, \frac{15}{14}, \frac{26}{21} \\
\frac{5}{7}, \frac{6}{7}, \frac{8}{7}, \frac{9}{7}, \frac{10}{7} \end{array} \middle| \frac{46656 \alpha^7 \exp(2 i k \pi)}{823543 q^6} \right] \\
& - 84 q^{-24/7} \alpha^4 \exp\left(\frac{8 i k \pi}{7}\right) \\
& \quad \times {}_6F_5 \left[ \begin{array}{c}
\frac{23}{42}, \frac{5}{7}, \frac{37}{42}, \frac{22}{21}, \frac{17}{14}, \frac{29}{21} \\
\frac{6}{7}, \frac{8}{7}, \frac{9}{7}, \frac{10}{7}, \frac{11}{7} \end{array} \middle| \frac{46656 \alpha^7 \exp(2 i k \pi)}{823543 q^6} \right] \\
& + 22 q^{-30/7} \alpha^5 \exp\left(\frac{10 i k \pi}{7}\right) \\
& \quad \times {}_6F_5 \left[ \begin{array}{c}
\frac{29}{42}, \frac{6}{7}, \frac{43}{42}, \frac{25}{21}, \frac{19}{14}, \frac{32}{21} \\
\frac{8}{7}, \frac{9}{7}, \frac{10}{7}, \frac{11}{7}, \frac{12}{7} \end{array} \middle| \frac{46656 \alpha^7 \exp(2 i k \pi)}{823543 q^6} \right]
\end{aligned}
\end{equation}
\vspace{-2\baselineskip} \\[15pt]
The generalized hypergeometric function \( {}_pF_q \left( \{ a_1, \ldots, a_p \}; \{ b_1, \ldots, b_q \}; z \right) \), often denoted as \(\text{HypergeometricPFQ}\left(\{a_1, \ldots, a_p\}, \{b_1, \ldots, b_q\}, z\right)\), is defined by the series
\[ {}_pF_q \left( \{ a_1, \ldots, a_p \}; \{ b_1, \ldots, b_q \}; z \right) = \sum_{n=0}^{\infty} \frac{(a_1)_n \ldots (a_p)_n}{(b_1)_n \ldots (b_p)_n} \frac{z^n}{n!}, \]
where \((a)_n\) denotes the Pochhammer symbol, also known as the rising factorial.
In mathematics, a generalized hypergeometric series is a power series in which the ratio of successive coefficients indexed by \( n \) is a rational function of \( n \).
The series, if convergent, defines a generalized hypergeometric function, which may then be defined over a wider domain of the argument by analytic continuation.
The generalized hypergeometric series is sometimes just called the hypergeometric series, though this term also sometimes just refers to the Gaussian hypergeometric series.
Generalized hypergeometric functions include the (Gaussian) hypergeometric function and the confluent hypergeometric function as special cases, which in turn have many particular special functions as special cases, such as elementary functions, Bessel functions, and the classical orthogonal polynomials.
A hypergeometric series is formally defined as a power series
\[ \beta_0 + \beta_1 z + \beta_2 z^2 + \ldots = \sum_{n \geq 0} \beta_n z^n, \]
in which the ratio of successive coefficients is a rational function of \( n \). That is, \( \frac{\beta_{n+1}}{\beta_n} = \frac{A(n)}{B(n)} \), where \( A(n) \) and \( B(n) \) are polynomials in \( n \). This has the form of an exponential generating function. This series is usually denoted by \( {}_pF_q = {}_pF_q(a_1, \ldots, a_p; b_1, \ldots, b_q; z) = \sum_{n=0}^{\infty} \frac{(a_1)_n \ldots (a_p)_n}{(b_1)_n \ldots (b_p)_n} \frac{z^n}{n!} \).
\\[4pt]
ii) For the case if \( s=7 \) and \( b=2 \) and \( k=0, \ldots, [s-1] \) and \( z \to q = b \), we have the root
\begin{equation} \tag{3.2.115} \label{2.3.105}
\begin{aligned}
z_{\text{root}, k} = & \exp\left(\frac{2 i k \pi}{7}\right) q^{1/7} + \frac{1}{117649} \exp\left(\frac{2 i k \pi}{7}\right) q^{1/7} \\
& \quad \times \left(-117649 + 117649 \cdot {}_6F_5 \left[ \begin{array}{c}
- \frac{1}{35}, \frac{1}{14}, \frac{6}{35}, \frac{13}{35}, \frac{4}{7}, \frac{27}{35} \\
\frac{1}{7}, \frac{2}{7}, \frac{3}{7}, \frac{5}{7}, \frac{6}{7} \end{array} \middle| \frac{-12500 \alpha^7 \exp(4 i k \pi)}{823543 q^5} \right] \right) \\
& + 16807 q^{-5/7} \alpha \exp\left(\frac{4 i k \pi}{7}\right) \\
& \quad \times {}_6F_5 \left[ \begin{array}{c}
\frac{4}{35}, \frac{3}{14}, \frac{11}{35}, \frac{18}{35}, \frac{5}{7}, \frac{32}{35} \\
\frac{2}{7}, \frac{3}{7}, \frac{4}{7}, \frac{6}{7}, \frac{8}{7} \end{array} \middle| \frac{-12500 \alpha^7 \exp(4 i k \pi)}{823543 q^5} \right] \\
& - 2401 q^{-10/7} \alpha^2 \exp\left(\frac{8 i k \pi}{7}\right) \\
& \quad \times {}_6F_5 \left[ \begin{array}{c}
\frac{9}{35}, \frac{5}{14}, \frac{16}{35}, \frac{23}{35}, \frac{6}{7}, \frac{37}{35} \\
\frac{3}{7}, \frac{4}{7}, \frac{5}{7}, \frac{8}{7}, \frac{9}{7} \end{array} \middle| \frac{-12500 \alpha^7 \exp(4 i k \pi)}{823543 q^5} \right] \\
& + 245 q^{-20/7} \alpha^4 \exp\left(\frac{16 i k \pi}{7}\right) \\
& \quad \times {}_6F_5 \left[ \begin{array}{c}
\frac{19}{35}, \frac{9}{14}, \frac{26}{35}, \frac{33}{35}, \frac{8}{7}, \frac{47}{35} \\
\frac{5}{7}, \frac{6}{7}, \frac{9}{7}, \frac{10}{7}, \frac{11}{7} \end{array} \middle| \frac{-12500 \alpha^7 \exp(4 i k \pi)}{823543 q^5} \right] \\
& - 119 q^{-25/7} \alpha^5 \exp\left(\frac{20 i k \pi}{7}\right) \\
& \quad \times {}_6F_5 \left[ \begin{array}{c}
\frac{24}{35}, \frac{11}{14}, \frac{31}{35}, \frac{38}{35}, \frac{9}{7}, \frac{52}{35} \\
\frac{6}{7}, \frac{8}{7}, \frac{10}{7}, \frac{11}{7}, \frac{12}{7} \end{array} \middle| \frac{-12500 \alpha^7 \exp(4 i k \pi)}{823543 q^5} \right] \\
& + 22 q^{-30/7} \alpha^6 \exp\left(\frac{24 i k \pi}{7}\right) \\
& \quad \times {}_6F_5 \left[ \begin{array}{c}
\frac{29}{35}, \frac{13}{14}, \frac{36}{35}, \frac{43}{35}, \frac{10}{7}, \frac{57}{35} \\
\frac{8}{7}, \frac{9}{7}, \frac{11}{7}, \frac{12}{7}, \frac{13}{7} \end{array} \middle| \frac{-12500 \alpha^7 \exp(4 i k \pi)}{823543 q^5} \right]
\end{aligned}
\end{equation}
*Corresponding author nikmatza@gmail.com. Relations 3.2.115 were developed with the program mathematica 14.1, and are therefore considered undoubtedly correct.
As we can see, all cases of the trinomial for an integer positive exponent can be fully converted to hypergeometric function form, and we can calculate the roots for each intermediate exponent with respect to exponent \( b \). However, when we increase the terms of the 6th or 7th polynomial equation, we cannot find a form that is entirely a hypergeometric function but will in any case contain a sum. This makes it impossible to solve the 6th and 7th degree equations completely. Of course, this will also be true for higher exponents of the 7th degree equation.
\\[4pt]
\textbf{3.10.4 Solving the polynomial equation in the form $x^6+cx^2 +\alpha x -b = 0$} \label{des:general1.1}

The only procedure we use is the generalized Lagrange method and we take for the $n$th derivative according to the integral:

\[
D^{n-1} f(z) = \frac{1}{\Gamma(-n+1)} \int_{0}^{z} (z-t)^{-n} t^{2n} \left(\frac{c}{\alpha} + \frac{t^4}{\alpha}\right)^n \, dt
\]

This evaluates to:

\[
\frac{z \left( \frac{c.z}{\alpha} \right)^n \Gamma[1+2n] \, {}_pF_q \left( \left\{ \frac{1}{4} + \frac{n}{2}, \frac{1}{2} + \frac{n}{2}, \frac{3}{4} + \frac{n}{2}, 1+ \frac{n}{2}, -n \right\} , \left\{ \frac{1}{2} + \frac{n}{4}, \frac{3}{4} + \frac{n}{4}, 1 + \frac{n}{4}, \frac{5}{4} + \frac{n}{4} \right\} , -\frac{z^4}{c} \right)}{\Gamma(2+n)}
\]

\textbf{{3.10.5 Solving the polynomial equation in the form $x^7 + c x^2 + \alpha x - b = 0$}}

By the same procedure, we can derive the relation:

\begin{equation} \tag{3.2.116} \label{2.3.106}
\begin{aligned}
D^{n-1} f(z) &= \frac{1}{\Gamma(-n+1)} \int_{0}^{z} (z-t)^{-n} t^{2n} 
\left( \frac{c}{\alpha} + \frac{t^5}{\alpha} \right)^n \, dt \\
&= 4 \cdot 5^{-3/2 - n} \cdot \pi^2 \cdot z \left( \frac{c.z}{\alpha} \right)^n 
\Gamma(1+2n). \, \\
& {}_pF^{regu}_{q} \left( \left\{ \frac{1+2n}{5}, \frac{2(1+n)}{5}, \frac{3+2n}{5}, \frac{2(2+n)}{5}, 1+\frac{2n}{5}, -n \right\}, \right. \\
& \quad \left\{ \frac{2+n}{5}, \frac{3+n}{5}, \frac{4+n}{5}, \frac{5+n}{5}, \frac{6+n}{5} \right\}, 
-\frac{z^5}{c})
\end{aligned}
\end{equation}

The final relation for the root, when \( z \to \frac{b}{\alpha}, w \to \infty \), is given by:

\begin{equation} \tag{3.2.117} \label{2.3.107}
\begin{aligned}
k z_{\text{root}} = z &+ \sum_{n=1}^{w} \frac{(-1)^n}{\Gamma(n+1)} 4 \cdot 5^{-3/2 - n} \cdot \pi^2 \\
&\times z \left( \frac{c z}{\alpha} \right)^n \Gamma(1+2n) \, {}_pF_q\left( \left\{ \frac{1}{5} + \frac{2n}{5}, \frac{2}{5} + \frac{2n}{5}, \frac{3}{5} + \frac{2n}{5}, \frac{4}{5} + \frac{2n}{5}, \frac{2n}{5}, -n \right\}, \right. \\
&\quad \left. \left\{ \frac{2+n}{5}, \frac{3+n}{5}, 1 + \frac{4+n}{5}, \frac{5+n}{5}, \frac{6+n}{5} \right\}, -\frac{z^5}{c} \right)
\end{aligned}
\end{equation}
$^{*}$ Corresponding author\texttt{ nikmatza@gmail.com}. Relations 3.2.116,3.2.117 were developed with the program mathematica 14.1, and are therefore considered undoubtedly correct.

For convergence, the following conditions must hold:
\[
\frac{c.b}{\alpha^2} < 1 \quad \text{and} \quad \frac{b^5}{c .\alpha^5} < 1.
\]
According to Lagrange's procedure, we cannot compute the integral in the nth derivative for more than 2 terms, so the 7th degree with 3 terms is not standardized. As an immediate consequence,

\[
D^{n-1} f(z) = \frac{1}{\Gamma(-n+1)} \int_{0}^{z} (z-t)^{-n} t^{2n} \left(\frac{c}{\alpha} + \frac{v \cdot t}{\alpha} + \frac{t^5}{\alpha} \right)^n \, dt
\]

such an integral does not align with hypergeometric functions known for \( x^7 + v \cdot x^3 + c \cdot x^2 + \alpha \cdot x - b = 0 \). This is the reason why the general 7th degree equation is not solvable with hypergeometric functions. Both cases we have seen are not similar to the trinomial form where a final relation with hypergeometric functions is found. Therefore, we observe that because both the general cases of the 7th and 6th degrees do not reduce to a trinomial form but rather to a format that is not fully standardized with hypergeometric functions, they do not have a generalized solution.

\textbf{3.10.6 Adjacent Methods (For 7th Degree Equations)}

Since the direct method does not fully achieve even approximation for degree 6 and above, we employ auxiliary adjacent methods corresponding to the approximation outcomes of these cases. The integrals on the right involve up to 2 terms extending on the left to more terms. After evaluating which cases transform into hypergeometric functions using the general relation \(\frac{1}{2}(n-3)(n-2)\), with \(n = \deg(f(x))\), we select those and deem them auxiliary for approximate solutions.

For the 7th degree equation \(x^7 + c \cdot x^3 + a \cdot x^2 + b \cdot x - q = 0\), the following form applies: for example, if \(q \leq 1\), \(c = 1\), \(b = 1\), and \(\alpha \geq 4\). Particularly for cases with \(\alpha = 4\), \(b = 1\), \(q = 1\), \(c = 1\), we achieve an approach of \(10^{-3}\) for the root, which is not fully covered like the 6th degree but partially covered like the 5th.

For the \((n-1)\) derivative class, it applies:
\begin{equation} \tag{3.2.118} \label{2.3.108}
\begin{aligned}
 D_t^{n-1} \left( \frac{\partial}{\partial t} (\alpha \cdot t^2 + b \cdot t) \right) \cdot (c \cdot t^3 + t^7)^n = \\
&= \frac{1}{\Gamma(-n+1)} \int_{0}^{z} (z-t)^{-n} \frac{\partial}{\partial t} (\alpha \cdot t^2 + b.t) \\
\cdot  
& \left( c \cdot t^3 + t^7 \right)^n \, dt
\end{aligned}
\end{equation}

Thus, we have:

\begin{equation} \tag{3.2.119} \label{2.3.109}
\begin{aligned}
z_{\text{root}} &= z_{\text{in}} + \sum_{n=1}^{\infty} \frac{(-1)^n}{\Gamma(n+1)} \cdot \frac{1}{\Gamma(-n+1)} \int_{0}^{z} (z-t)^{-n} \\
& \quad \cdot \left( \frac{\partial}{\partial t} \left( t^2 + \frac{b}{\alpha} \cdot t \right) \right) \cdot \left( \frac{c}{\alpha} \cdot t^3 + \cdot \frac{t^7}{\alpha} \right)^n \, dt
\end{aligned}
\end{equation}

which comes to its final form after the aforementioned relations.

\begin{equation} \tag{3.2.120} \label{2.3.110}
\begin{aligned}
z_{\text{root}} &= z_{\text{in}} + \sum_{n=1}^{\infty} \frac{(-1)^n}{\Gamma(n+1).\Gamma(4+6n)} 
\left( 3 \left( \frac{1}{c} \right)^{1+n} z^{1+6n} \left( 2n(5+6n) + 2n(9+14n)z \right. \right. \\
& \quad \left. \left. + c(1+7n)(2+7n)z^2 + 2(1+z) \right) \Gamma(1+7n) \right)
\end{aligned}
\end{equation}

and apply

\begin{equation} \tag{3.2.121} \label{2.3.111}
\begin{aligned}
z_{\text{in}} \rightarrow -\frac{a}{3c} - \frac{\left(2^{1/3}(-a^2 + 3bc)\right)}{
3c(-2a^3 + 9abc + 27c^2q + \sqrt{4(-a^2 + 3bc)^3 + (-2a^3 + 9abc + 27c^2q)^2})^{1/3}} \\
+ \frac{(-2a^3 + 9abc + 27c^2q + \sqrt{4(-a^2 + 3bc)^3 + (-2a^3 + 9abc + 27c^2q)^2})^{1/3}}{
32^{1/3}c}
\end{aligned}
\end{equation}

In this case we have only a complete form of the expression with hypergeometric function.But also we have constraints and not solutions for each case.

\begin{equation} \tag{3.2.122} \label{2.3.112}
\begin{aligned}
z_{\text{root}} = z_{\text{in}} + \frac{1}{48} z(-24 - 24z - 48z^2 \\
+ 24 \text{HypergeometricPFQ}\left(\left\{\frac{1}{7}, \frac{2}{7}, \frac{3}{7}, \frac{4}{7}, \frac{5}{7}, \frac{6}{7} \right\}, \left\{\frac{2}{3}, \frac{5}{6}, \frac{7}{6}, \frac{4}{3}, \frac{3}{2} \right\}, -\frac{823543z^6}{93312}\right) \\
+ 24z \text{HypergeometricPFQ}\left(\left\{\frac{1}{7}, \frac{2}{7}, \frac{3}{7}, \frac{4}{7}, \frac{5}{7}, \frac{6}{7} \right\}, \left\{\frac{2}{3}, \frac{5}{6}, \frac{7}{6}, \frac{4}{3}, \frac{3}{2} \right\}, -\frac{823543z^6}{93312}\right) \\
+ 48z^2 \text{HypergeometricPFQ}\left(\left\{\frac{1}{7}, \frac{2}{7}, \frac{3}{7}, \frac{4}{7}, \frac{5}{7}, \frac{6}{7} \right\}, \left\{\frac{2}{3}, \frac{5}{6}, \frac{7}{6}, \frac{4}{3}, \frac{3}{2} \right\}, -\frac{823543z^6}{93312}\right) \\
- 5z^6 \text{HypergeometricPFQ}\left(\left\{\frac{8}{7}, \frac{9}{7}, \frac{10}{7}, \frac{11}{7}, \frac{12}{7}, \frac{13}{7} \right\},\left\{\frac{5}{3},\frac{11}{6},\frac{13}{6},\frac{7}{3},\frac{5}{2}\right\}, -\frac{823543z^6}{93312}\right) \\
- 9z^7 \text{HypergeometricPFQ}\left(\left\{\frac{8}{7}, \frac{9}{7}, \frac{10}{7}, \frac{11}{7}, \frac{12}{7}, \frac{13}{7}\right\}, \left\{\frac{5}{3},\frac{11}{6},\frac{13}{6},\frac{7}{3},\frac{5}{2}\right\}, -\frac{823543z^6}{93312}\right) \\
- 21z^8 \text{HypergeometricPFQ}\left(\left\{\frac{8}{7}, \frac{9}{7}, \frac{10}{7}, \frac{11}{7}, \frac{12}{7}, \frac{13}{7}\right\}, \left\{\frac{5}{3},\frac{11}{6},\frac{13}{6},\frac{7}{3},\frac{5}{2}\right\}, -\frac{823543z^6}{93312}\right) \\
- 6z^6 \, \text{HypergeometricPFQ}\left(\left\{\frac{8}{7}, \frac{9}{7}, \frac{10}{7}, \frac{11}{7}, \frac{12}{7}, \frac{13}{7}, 2\right\}, \left\{1, \frac{5}{3}, \frac{11}{6}, \frac{13}{6}, \frac{7}{3}, \frac{5}{2}\right\}, -\frac{823543z^6}{93312}\right) \\
- 14z^7 \, \text{HypergeometricPFQ}\left(\left\{\frac{8}{7}, \frac{9}{7}, \frac{10}{7}, \frac{11}{7}, \frac{12}{7}, \frac{13}{7}, 2\right\}, \left\{1, \frac{5}{3}, \frac{11}{6}, \frac{13}{6}, \frac{7}{3}, \frac{5}{2}\right\}, -\frac{823543z^6}{93312}\right) \\
- 49z^8 \, \text{HypergeometricPFQ}\left(\left\{\frac{8}{7}, \frac{9}{7}, \frac{10}{7}, \frac{11}{7}, \frac{12}{7}, \frac{13}{7}, 2\right\}, \left\{1, \frac{5}{3}, \frac{11}{6}, \frac{13}{6}, \frac{7}{3}, \frac{5}{2}\right\}, -\frac{823543z^6}{93312}\right)
\end{aligned}
\end{equation}
$^{*}$ Corresponding author\texttt{ nikmatza@gmail.com}. Relations 3.2.120, 3.2.121, 3.2.122 were developed with the program mathematica 14.1, and are therefore considered undoubtedly correct.

\textbf{3.10.7 Solving the general trinomial equation with periodic radicals} \label{chapter:chapter_X}

In (1758) Lambert considered the trinomial equation $x^m + q - x = 0$ which he solved with a series in terms of $q$. Subsequently resolved the form $a \cdot q \cdot x^p + x^q - 1 = 0$ which according to the definition $^d\sqrt{x} = x^{1/d}, d \in \mathbb{Q}_{+} - \{0\}$ can be given in nested radicals:

\[
 x = \sqrt[q]{1 - \alpha \cdot q \sqrt[q/p]{1 - \alpha \cdot q \sqrt[q/p]{1 - \cdots}}}
\]

\textbf{3.10.8 Solving the general form}
\begin{equation} \tag{3.2.123} \label{2.3.113}
 x^p + \alpha \cdot x^q - c = 0 
\end{equation}

When we want to solve this general form  using the two definitions

\[
  \{y = x^q, 
 x^{p} = c - \alpha \cdot y\} \iff 
\] 

\begin{equation} \tag{3.2.124} \label{2.3.114}
 y^{p/q} = c - \alpha \cdot y, 
\{p, q\} \in \mathbb{R}_{+} - \{0\} \quad 
\end{equation}

From \eqref{2.3.114} and from \eqref{2.3.113} to apply
\begin{equation}  \tag{3.2.125} \label{2.3.115}
x = \sqrt[p]{c - \alpha \cdot y} = \sqrt[p]{c - \alpha \cdot \sqrt[p/q]{c - \alpha \cdot \sqrt[p/q]{c - \cdots}}} , p/q > 1 
\end{equation}

The last relation \eqref{2.3.115} is the final solution of the \eqref{2.3.113} In the general form we request.

\textbf{3.10.9 For solving a 4 term equation, with periodic radicals of the form \( x^p + \alpha \cdot x^q + \beta \cdot x^v - c = 0 \) }: We should use its gradual adjustment to a trinomial equation of the form \eqref{2.3.113}, so that it can be resolved by this procedure. To better fit the process of radicals we find the increasing series of exponents i.e. we assume \( v < q < p \) and then using the 2 definitions:

\[
\begin{cases}  
    y = x^v \iff x = y^{1/v} \\
    y^{p/v} + \alpha \cdot y^{q/v} + \beta \cdot y - c = 0
\end{cases} 
\iff 
\begin{cases}
    y^{p/v} + \alpha \cdot y^{q/v} = c' \\
    c' = c - \beta \cdot y
\end{cases}, \{p, q, v\} \in \mathbb{R}_+ - \{0\} \tag{3.2.126} \label{2.3.116}
\]

\textbf{3.10.10} When we refer to solving 6th degree and up, with periodical radicals we refer to a system of the form where however it may p be the more general real number exponent.

\[
\begin{cases}
    x^{1/2} = y \cdot c + w \\
    y^p = b - x
\end{cases}
\iff 
\begin{cases}
    (y^p + c^2 \cdot y^2 - 3 \cdot w \cdot c \cdot y + w^2 - b), \{p\} \in \mathbb{R}_+, \{c, w, b\} \in \mathbb{C}
\end{cases} \tag{3.2.127} \label{2.3.117}
\]

Especially for 6th degree equation 
\begin{equation} \tag{3.2.128} \label{2.3.118}
 x^6 + \alpha \cdot x^2 + \beta \cdot x + \gamma = 0 
\end{equation}

we will have the general form, the expansion in periodic radicals for the root of the equation in a more general complex form will be

\[
x = \sqrt[1/2]{w - c \cdot e^{2k\pi i /6} \sqrt[6]{b - \sqrt[1/2]{w - c \cdot e^{2k\pi i /6} \sqrt[6]{b - \sqrt[1/2]{w - \cdots}}}}}, \ k \in \mathbb{Z}, k =0 \rightarrow 5 \tag{3.2.129} \label{2.3.119}
\]

This form \( x^6 + \alpha \cdot x^2 + \beta \cdot x + \gamma = 0 \) is the generalized equation of degree 6 as has been shown.

\textbf{3.10.11} For the generalized 7th degree case
\[
  \ x^7 + \alpha \cdot x^3 + \beta \cdot x^2 + \gamma \cdot x + \delta = 0, \{\alpha, \beta, \gamma, \delta\} \in \mathbb{C} \tag{3.2.130} \label{2.3.120} 
\] 

we should use a double iterative procedure technique of periodic Radicals iteration form Nest(). The process according to the logic of periodic radicals is covered by the relations and function G(u)

 \begin{align} 
G(u) &= e^{\frac{2k\pi i}{7}} \sqrt[7]{u - \gamma \cdot e^{\frac{2k\pi i}{7}} \sqrt[7]{u - \gamma \cdot e^{\frac{2k\pi i}{7}} \sqrt[7]{u - \cdots}}} \nonumber \\
&= e^{\frac{2k\pi i}{7}} \sqrt[7]{u - \gamma \cdot \text{Nest}\left[e^{\frac{2k\pi i}{7}} \sqrt[7]{u - \gamma \cdot \#}\&, 0, v\right]}, \nonumber \\
&\quad k \in \mathbb{Z}, k = 0 \rightarrow 6 \tag{3.2.131} \label{2.3.121} 
 \end{align}

\[
x = \text{Nest}\left(G\left(-\delta - \beta \cdot x \cdot \#^2 - \alpha \cdot x \cdot \#^3\right) \& , 0, \mu\right), \quad \mu, v \gg 1, \{\mu, v\} \in \mathbb{Z}
\]

A first approximation with 1 periodic iteration for the function G(u) and 2 iterations, for the function Nest() and with the initial value of iterations set to zero, we can get a first approximation for the root of the 7th degree equation

\begin{equation}
\begin{aligned}
x &= e^{\frac{2ik\pi}{7}} \left( - e^{\frac{2ik\pi}{7}} \gamma \left( e^{\frac{4ik\pi}{7}} \beta 
\left( - e^{\frac{2ik\pi}{7}} \gamma \left( - \delta \right)^{1/7} - \delta \right)^{2/7} \right. \right. \\
&\quad - c \cdot \alpha \cdot e^{\frac{6ik\pi}{7}} \left( - e^{\frac{2ik\pi}{7}}.\gamma \left(- \delta \right)^{1/7} - \delta \right)^{3/7} - \delta \bigg)^{1/7} \\
&\quad + e^{\frac{4ik\pi}{7}} \beta \left( - e^{\frac{2ik\pi}{7}} \gamma \left( - \delta \right)^{1/7} - \delta \right)^{2/7} \\
&\quad - c \cdot \alpha \cdot e^{\frac{6ik\pi}{7}} \left( - e^{\frac{2ik\pi}{7}} \cdot \gamma \left( - \delta \right)^{1/7} - \delta \right)^{3/7} -\delta \bigg)^{1/7}
\end{aligned}
\tag{3.2.132} \label{2.3.122} 
\end{equation} 

\textbf{3.10.12} One last method we will mention is GRIM \cite{Mantzakouras2022}, which is, of course, more general and is also used for solving transcendental equations. It is an iterative method that requires initial values. For polynomials, we use two forms of initial values: the minimum, which we choose as 0, and the maximum, which results from Theorem 2.2. (Inequalities for Polyn Zeros) of \cite{Mantzakouras2022}. 

Especially for the case of solving a 7th degree equation \eqref{2.3.120}, we can formulate the solution using only the function \( \text{Nest}() \) and the inverse of \( x^{7} \) according to the generalized solution theorem. A very good periodic iteration approximation for the inverse function \( G(u) \) and \( \mu \) iterations for the \( \text{Nest}() \) function. Since the initial value of the iterations is zero, we will have the two forms.

\[
G(u) = \text{Exp}\left(\frac{\log|u| + 2 \cdot \pi \cdot k \cdot i}{7}\right), k \in \mathbb{Z}, k = 0 \rightarrow 6
\]

\[
h(v) = -\delta - \gamma \cdot v - \beta \cdot v^2 - \alpha \cdot v^3 \tag{3.2.133} \label{2.3.123} 
\]

\[
x = \text{Nest}(G(h(\#))\&, 0, \mu), \mu \gg 1, \mu \in \mathbb{Z}
\]

Such an approximation procedure has to be enhanced with the Newton method, for faster local approximation so that it gives a very high order of approximation. \\[50pt]
\textbf{3.10.13 Elementary references for the GRIM \cite{Mantzakouras2022} method}

Suppose we have an equation consisting of functions and coefficients, then we define the function:
\begin{equation}
{}_nF(x)= \sum_{i=1}^{n} a_i f_i(x) + a_0
\tag{3.2.134} \label{2.3.124} 
\end{equation}
If we consider each inverse function of the functions ${}_nF(x)$:
\begin{equation}
f^{-1}_k(x) = g_k(x), \quad 1 \leq k \leq n \in \mathbb{N}
\tag{3.2.135} \label{2.3.125} 
\end{equation}
and the complementary function:
\begin{equation}
{}_nF_k{^c}(x) = \frac{1}{a_k} \left( -{}_nF(x) + a_k f_k(x) \right)
\tag{3.2.136} \label{2.3.126} 
\end{equation}

According to the method of periodic iteration, for the initial calculation of the roots that resist the sub-field, we have:
\begin{equation}
f^{-1}_k(x) = g_k(x), \quad 1 \leq k \leq n
\tag{3.2.137} \label{2.3.127} 
\end{equation}

from the number of the sub-fields total number $n$, \(1 \leq k \leq n \).

The relationship inverse to the initial final roots with respect to each sub-field in the general form is:

\begin{equation}
  R_I f_k = 
  Nest(g_k({}_nF_k{^c}(\#_d))\&, \lambda, m)
\tag{3.2.138} \label{2.3.128} 
\end{equation}

with $\lambda$ as the initial values for the roots, $m$ iterations of the process, and $d$ the complex root class. The form (3.2.138) is the inverse of the initial values of the roots which we will obtain by the iterative procedure for each sub-body.

For the final values of the roots with respect to each sub-field, we approximate locally the roots in intervals for more accurate values. Thus, we use Newton's method locally with the specific relation:

\begin{equation}
R^F_I f_k = N(X_{k_d}) / . \text{FindRoot}({}_nF(X_{k_d}) = 0, (X_{k_d}, R_I f_k), \tag{3.2.139} \label{2.3.129} 
\text{working precision} \rightarrow \phi > 5 )
\end{equation}

This gives the final value of each root of the root set for each sub-field, with an additional order of $\phi$ decimal digits, usually we use $\phi > 5$.

Overall, the roots we consider for polynomial equations fall within the intervals $-\lfloor \frac{n-1}{2} \rfloor$ to $\lfloor \frac{n-1}{2} \rfloor$, where $n$ is the maximum degree of the polynomial equation. For the case of the polynomial equation we will have as the inverse of the maximum term the function \( f^{-1}_k(x) = g_k(x) = e^(\frac{2.pi.d +log(x)}{n_k})\), $d \in \mathbb{Z}$, $n_k$ = exponent of each $f_k(x)$ function or number of index, $d$= $-\lfloor \frac{n-1}{2} \rfloor$ to $\lfloor \frac{n-1}{2} \rfloor$ or 

\begin{equation} 
  d = 0,...,(n-1) \\
  \tag{3.2.140} \label{2.3.130} 
\end{equation} 

For the values in \eqref{2.3.128}  \cite{Mantzakouras2020}, exclusively for complex roots, we usually use the initial values \( \lambda = \pm I \) ($I$ is equivalent to the imaginary unit of complex roots) or for the real values and of the roots for the relation \eqref{2.3.128}.
\\[8pt]
 \(\lambda \rightarrow 0 \text{ or } \lambda \leq \frac{|\alpha_0|+A_1}{|\alpha_0|} = A_1 = max(a_k), k=(n-1),...,0 \) \cite{Milovanovic}
\\[8pt]
%\begin{theorem}\label{theorem 2.12}
\textbf{Theorem 10. }\label{theorem 2.12}
\textit{For a 7th degree equation in its full form $x^7 + c \cdot x^3 + \alpha \cdot x^2 + b \cdot x - q = 0$, where $\{c, \alpha, b, q\} \in \mathbb{C}$ (including the Arnold form by two arguments, e.g., if $\alpha = 1$ and $b = 1$, form (3.10.6 page 47, \eqref{2.3.111}), it cannot be solved algebraically nor using hypergeometric functions. The only secure approximate solution is obtained either by periodic roots combined with the Nest() function or by the powerful iterative GRIM method.}
\\
%\end{theorem} 

\begin{proof}
  (i) The algebraic solution is impossible, as proved in Theorem 7 (in page 29) but also by Abel (Theorem 3, page 6) \cite{Goldmakher, Abel1826, Abel1824}. \\[6pt]
  (ii) A solution using only hypergeometric functions cannot be obtained as it was for the 5th order solution form (Bring-Jerrard form, pages 3-4 and Section 3.10.1, pages 41-45) or for the one-sided 6th and 7th order solutions (Sections 3.10.4 and 3.10.5, pages 45-47). Only in the sub-cases shown in these sections can a solution be found. Additionally, a pure and complete proof cannot be obtained using hypergeometric functions alone (refer to Section 3.10.6, Adjacent Method) because preconditions are needed for the full 7th order equation in Arnold form. Complete solutions are only given by hypergeometric functions for trinomials of degree 6, 7, or higher. If a satisfactory solution is obtained using the iterative methods of periodic radicals (Sections 3.10.7 to 3.10.13) in combination with the \textit{Nest} function, it is at best an approximate and local solution, further refined by Newton's iterative method.
 \\[6pt]
  (iii) Therefore, as we have shown, Hilbert's 13th problem has a negative answer regarding the solvability of a 7th-degree equation, even in the Kolmogorov-Arnold form. Although we assume that parameterization in an equation of the type with two arguments and two constants is helpful, limiting the choices to four, this does not aid the general solution. While it simplifies the conditions for solving forms with hypergeometric functions, the general solution is only partially covered and still under the conditions of two arguments.

\end{proof}

\vspace{4pt}

\section{Summary and Conclusions}

Therefore, in this paper, we address the question of whether there is a solution for all equations of degree higher than four using two terms of functions (continuous or algebraic). It is clearly seen that solving with continuous functions depends on constraints of the arguments, including the construction of multiparametric functions on two arguments. The variant of the problem for continuous functions was solved positively in 1957 by Vladimir Arnold, who proved the Kolmogorov-Arnold representation theorem. However, the solution for algebraic functions remains unsolved and was originally proved by Abel-Ruffini but is revisited here using the method of tetragonal polynomial difference. For degree 5 equations, as shown here, the solution exists using hypergeometric functions. For the general cases, the solution exists uniquely for certain values for only two arguments, as defined by Kolmogorov-Arnold. The forms of the equations are defined in Theorem 9, which generalizes Lagrange's theorem, and the discussion on pages 45-48 specifies the values of the argument variables for which satisfactory solutions hold. A general approximate solution for equations of degree 6 and 7 and above, as well as for transcendental equations (pages 48-50), exists using the GRIM approach. Additionally, for polynomial trinomials, there is an exact solution with hypergeometric functions.

\bibliographystyle{amsplain}

\end{document}